\documentclass[preprint,3p,nopreprintline]{elsarticle}

\usepackage{amssymb}

\usepackage{adjustbox}
\usepackage{afterpage}
\usepackage{amsmath}
\usepackage{bm}
\usepackage[overload]{empheq}
\usepackage{float}
\usepackage{hyperref}
\usepackage{lipsum}
\usepackage{longtable}
\usepackage{multicol}
\usepackage{multirow}
\usepackage{physics}
\usepackage{siunitx}
\usepackage{subfig}
\usepackage{tabularx}
\usepackage{todonotes}
\everymath{\displaystyle}

\DeclareMathOperator*{\argmin}{arg\,min\,}

\begin{document}
	
	\begin{frontmatter}
		
		
		
		\title{Physics-Informed Neural Networks for the High-Resolution Reconstruction of Flow Measurement Indicators in Fluid Dynamics}
		

        \author[1]{Irena Radišić\corref{cor1}}
		\author[0]{Raffaele Tirotta}
		\author[1,2]{Alberto Zingaro}
		\author[1]{Stefano Pagani}		
		\author[1]{Luca Dede'}
        
        \cortext[cor1]{Corresponding author: irena.radisic@polimi.it}
   
		\affiliation[0]{organization={Politecnico di Milano},
			addressline={Piazza Leonardo Da Vinci 32}, 
			city={Milan},
			postcode={20133}, 
			country={Italy}}
   
		\affiliation[1]{organization={MOX, Dipartimento di Matematica, Politecnico di Milano},
			addressline={Piazza Leonardo Da Vinci 32}, 
			city={Milan},
			postcode={20133}, 
			country={Italy}}
   
        \affiliation[2]{organization={ELEM Biotech S.L.},
			addressline={Via Laietana 26}, 
			city={Barcelona},
			postcode={08003},
			country={Spain}}

		\begin{abstract}
            Accurate, spatially resolved flow field measurements are essential for the reliable assessment of hemodynamic quantities in cardiovascular research and clinical practice. Experimental techniques, such as 4D flow MRI, PIV, or Doppler ultrasound, often yield data that are sparse, noisy, or under-resolved, particularly near vessel walls and in regions of complex flow. This limits the fidelity of distributed or derived hemodynamic indicators such as the wall shear stress and the clinical utility of such measurements. To address these challenges, we propose a physics-informed neural network (PINN) framework that integrates the incompressible Navier–Stokes equations with velocity measurements coming from experimental flow field data. By embedding physical laws into data, PINN enhances the reconstruction of velocity fields, enables the estimation of unmeasured quantities such as pressure and wall shear stress, and improves the spatial resolution of hemodynamic indicators. We show the effectiveness of our approach using both in silico and experimental data. First, we apply our method to the FDA nozzle benchmark, leveraging both control particle image velocimetry (PIV) measurements and computational fluid dynamics (CFD) simulations. Next, we apply our method to the more complex case of blood flow in an aneurysm model, exploiting in vitro 4D flow MRI data. In both cases, the synergy between data-driven learning and physics-based regularization yields results that align more closely with ground truth observations than standard CFD or pure data-driven approaches. Our findings highlight the potential of PINNs to improve the fidelity of under-resolved flow field measurements and yield spatially resolved hemodynamic indicators.

		\end{abstract}
		
		\begin{graphicalabstract}
			\begin{figure}[t]
				\centering
				\includegraphics[width = \textwidth]{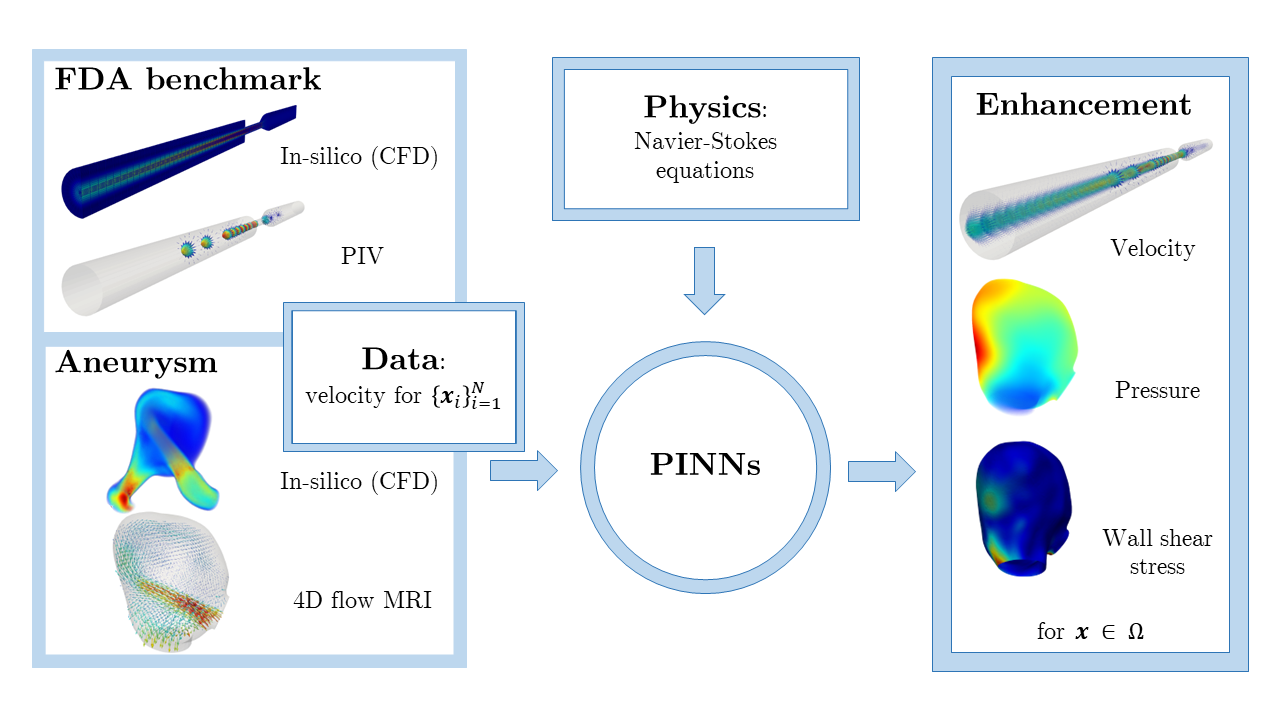}
				\caption{Graphical abstract. We use physics-informed neural networks (PINNs) to integrate physical knowledge through the Navier–Stokes equations with spatially under-resolved blood flow velocity measurements. This enables us to derive additional hemodynamic insights and achieve a spatially resolved solution.}
			\label{fig:summary}
			\end{figure}
            We use physics-informed neural networks (PINNs) to integrate physical knowledge through the Navier–Stokes equations with spatially under-resolved blood flow velocity measurements. This enables us to derive additional hemodynamic insights and achieve a spatially resolved solution.
		\end{graphicalabstract}
		
		
		
	\end{frontmatter}
	
	
 
	\section{Introduction}
	\label{sec:introduction}
	4D flow MRI is a clinically established and non-invasive imaging technique for the in vivo velocity measurement of cardiac and vascular flows~\cite{stankovic20144d,miyazaki2017validation,sughimoto2016effects}. It consists of a 3D time-varying description of the three components of the blood velocity vector field based on phase-contrast MRI~\cite{markl20124d}. The reconstructed space-time description of the blood velocity provides hemodynamics biomarkers for diagnosis~\cite{casas20164d}, including the wall shear stress, which is used to assess the risk of rupture in a blood vessel or to assist the diagnosis of important cardiac pathologies like stenosis or aneurysms~\cite{baek2009wall,caro2009discovery}. However, the accuracy in reconstructing biomarker indicators is constrained by the limited spatio-temporal resolution and signal-to-noise ratio~\cite{miyazaki2017validation}. Indeed, 4D flow MRI is not able to capture the small scales or flow-field features, such as recirculation regions or transitions to turbulence, as pointed out in~\cite{zingaro2021hemodynamics, harloff20093d}. Moreover, data acquisition is particularly dispersed and uncertain on the boundary~\cite{ngo2019four, pravdivtseva2022influence}, reducing the reliability in the estimation of vessel stresses and quantifying other spatially-resolved biomarkers~\cite{zingaro2022}. Another source of uncertainty arises from the time resolution of the measurements. Using low time-resolution acquisition techniques leads to the loss of information on instantaneous events in the flow field, such as the sudden formation of vortex cores in short time scale, which would support the evaluation of the rupture risk of cerebral aneurysms~\cite{futami2019identification, le2013vortex}.
    \par
    On the other hand, the use of computational fluid dynamics (CFD) and simulation-based approaches for patient specific geometries yields fluid flow information with arbitrary space and time resolution. Standard high fidelity methods, such as the finite element method (FEM) or the finite volume method (see, e.g., \cite{quarteroni2017numerical}) can be very time consuming~\cite{totorean_patient-specific_2022}. Moreover, these methods often fail to account for the inherent variability of the data, and the personalization of these models requires multiple stages and is particularly laborious. Indeed, a typical simulation pipeline consists of image segmentation, meshing, discretization and parameter tuning, and in the case of model personalization some of these steps may even be iterated.
    \par
	The goal of this work is to develop a new approach based on scientific machine learning~\cite{goodfellow2016deep,Tassi2022,quarteroni_combining_2025} to enhance 4D flow MRI post-processing in a computationally efficient and accurate manner. Specifically, we aim to reconstruct high-resolution velocity fields from 4D flow MRI data by leveraging physical constraints in a physics-informed neural networks (PINNs) framework~\cite{toscano_pinns_2025}. PINNs are an established technique for data-model integration in scenarios characterized by shortage of data~\cite{raissi2017physics} and, additionally, avoid geometry segmentation and mesh generation due to their mesh-free formulation.
    To compensate the limitations of 4D flow MRI, the training of the neural network leverage physical knowledge encoded into the Navier-Stokes equations, which are enforced in the loss function by computing a residual on a dataset of collocation points. This physical regularization is able to improve the reconstruction of the velocity field and its derived quantities, like the wall shear stress, while also allowing the estimation of quantities that are hardly measurable in-vivo, like the pressure distribution.
    \par
	Several approaches have already been developed to recover the pressure field and the wall shear stress from 4D flow MRI velocity measurements with an appropriate level of accuracy, as well as other hemodynamics biomarkers, see e.g. \cite{marsden25}. Zhang et al.~\cite{zhang20194d} introduced a weighted least-squares reconstruction method for the spatial integration of pressure gradients. The main advantage of this approach is that the pressure estimation is obtained non-invasively. However, the pressure gradients are recovered a posteriori from the velocity field, so the accuracy on the pressure field is strongly dependent on the accuracy of the velocity field. In ~\cite{ferdian2022wssnet}, 4D flow MRI data and CFD simulations are used to train a convolutional neural network in order to assess hemodynamic quantities, such as wall shear stress. Other than reconstructing quantities not present in the original data, the inclusion of the physics not only augments sparse data, but also has beneficial effects on the regularization of noisy data. For example, in~\cite{zhang_divergence-free_2021}, the authors use a divergence-free constraint in order to effectively denoise data coming from 4D flow MRI measurements.
    \par
    PINNs are becoming widely used in computational hemodynamics~\cite{cai2022physics}, particularly alongside domain segmentation and 1D reduced order modeling. Important contributions are the ones for the reconstruction of blood pressure given by Kissas et al.~\cite{kissas2020machine}, and for the estimation of brain hemodynamics and arterial wall dilatation by Sarabian et al.~\cite{sarabian2022physics}. Both of these works are based on 4D flow MRI data. Moreover, outside of the field of computational hemodynamics, PINNs have already succesfully been used for the reconstruction of latent flow field quantities from imaging data~\cite{cai_flow_2021}. 
    \par
	In this paper, we propose a computational pipeline for high-resolution velocity, wall-shear stress and pressure reconstruction based on PINN. To validate the approach we consider two test cases: the FDA nozzle benchmark case, known as Computational Round Robin n.~1~\cite{wiki:xxx,stewart2012assessment}, and a realistic patient-derived aneurysm experimental model. We provide a graphical abstract of our work in Figure~\ref{fig:summary}. In the first test case, we generate an in silico dataset by carrying out CFD simulations of the benchmark problem provided by the FDA. To mimic 4D MRI data, we consider the velocity samples on some selected cross-sections, and we match the CFD results on the remaining part of the dataset. The procedure is then repeated on an experimental dataset given by the three velocity components obtained from PIV~\cite{hariharan2011multilaboratory, NCIHub43}. Next, we assess the flow field in an aneurysm. In this case, the data are given by open-access 4D flow MRI data for an in vitro patient-derived aneurysm model~\cite{pravdivtseva2022influence}. In particular, our method allows to compute the reconstructed velocity, pressure and wall shear stress wherever these quantities are not available.
	\par
	This paper is organized as follows. In Section~\ref{sec:methodology}, we introduce the methodology and our approach based on PINNs. Section~\ref{sec:FDAbench} is devoted to the FDA nozzle benchmark. In Section~\ref{sec:aneurysm}, we apply our methodology to an aneurysm model. We discuss the results of the two different test cases in Section~\ref{sec:discussion}, and the limitations of our approach in Section~\ref{sec:limitations}. Finally, in Section~\ref{sec:conclusions} we draw our conclusions and we discuss the results.
	

	\section{Materials \& Methods}
	\label{sec:methodology}
	PINNs are a class of artificial neural networks (ANNs) that are trained to simultaneously leverage data and physical knowledge. The physical knowledge is typically incorporated by PDEs in the loss function~\cite{raissi2017physics}. 
    In this work, we apply PINNs in the context of vascular hemodynamics, and since we assume that blood obeys the incompressible Navier-Stokes equations for a Newtonian fluid~\cite{quarteroni2019mathematical}, the PDE loss term is accordingly specified.
    \par
    Let $\Omega\subset\mathbb R^3$ be a bounded domain with boundary $\partial\Omega=\Gamma_\mathrm{inl}\cup\Gamma_\mathrm{wall}\cup\Gamma_\mathrm{out}$, with $\Gamma_\mathrm{inl},\Gamma_\mathrm{wall},\Gamma_\mathrm{out}$ disjoint subsets of $\partial\Omega$. We consider the velocity field $\bm{u}_\mathrm{ex}$ and pressure field $p_\mathrm{ex}$,
    \begin{equation*}
        \bm{u}_\mathrm{ex}:\Omega\times[0,T]\to\mathbb R^3, \qquad p_\mathrm{ex}:\Omega\times[0,T]\to\mathbb R,
    \end{equation*}
    exact solutions of the dimensionless steady-state Navier-Stokes equations in $\Omega$ with suitable boundary conditions, that is, for every time $t\in [0,T]$ it holds that:
    \begin{subequations}
    \label{eq:NS_strong}
    \begin{empheq}[left=\empheqlbrace]{align}
    & (\bm{u}_\mathrm{ex}\cdot \nabla)\bm{u}_\mathrm{ex} + \nabla p_\mathrm{ex}  - \frac{1}{\mathbb{R}\mathrm{e}}\Delta\bm{u}_\mathrm{ex} = \mathbf{0}, \quad &&\mbox{in} \ \Omega, \label{eq:NS_mom_cons} \\
    & \nabla \cdot \bm{u}_\mathrm{ex} = 0, \quad &&\mbox{in} \ \Omega, \label{eq:NS_mass_con} \\
    & \bm{u}_\mathrm{ex} = \bm{g}, \quad &&\mbox{on} \ \Gamma_\mathrm{inl}, \label{eq:NS_bc_inl} \\
    & \bm{u}_\mathrm{ex} = \bm{0}, \quad &&\mbox{on} \ \Gamma_\mathrm{wall}, \label{eq:NS_bc_wall} \\
    & \left(\frac{1}{\mathbb{R}\mathrm{e}}\nabla{\bm{u}_\mathrm{ex}} - p_\mathrm{ex}\bm I\right)\bm n = \bm{0}, \quad &&\mbox{on} \ \Gamma_\mathrm{out}. \label{eq:NS_bc_out}
    \end{empheq}
    \end{subequations}
    Eq.s~\eqref{eq:NS_mom_cons} and~\eqref{eq:NS_mass_con} are the momentum balance and mass conservation equations respectively, Eq.~\eqref{eq:NS_bc_inl} is the Dirichlet inflow boundary condition, Eq.~\eqref{eq:NS_bc_wall} is the Dirichlet no-slip boundary condition, whereas Eq.~\eqref{eq:NS_bc_out} is the homogeneous Neumann outflow boundary condition. Given a volumetric flow rate $Q$, we define the Reynolds number $\mathbb{R}\mathrm{e}$ over a characteristic section as:
    \begin{equation*}
     \mathbb{R}\mathrm{e} = \frac{\rho U D_t}{\mu},\qquad U = \frac{4 Q}{\pi D_\mathrm{t}^2},
    \end{equation*}
    where $\rho$ and $\mu$ are respectively the fluid density and dynamic viscosity, $U$ is a characteristic velocity, and $D_t$ is a characteristic diameter. We report the fluid flow properties parameters used throughout this work in Table~\ref{tab:fluid_flow_properties}.
    \begin{table}[h]
		\centering
			\begin{tabular}{ p{6em} p{6em} p{6em} p{6em} p{7em}}
				\hline
				 & Unit & FDA & ANE-CFD & ANE-4DFLOW \\
				\hline
				$Q$ & $\si[per-mode=reciprocal]{\meter\cubed\per\second}$ & $5.2062\cdot10^{-6}$ & $3.8\cdot10^{-6}$ & $3.8\cdot10^{-6}$  \\
				$D_t$ & $\si[per-mode=reciprocal]{\meter}$ & $0.004$ & $0.0035$& $0.013$ \\
				$\mu$ & $\si[per-mode=reciprocal]{\pascal\second}$ & $0.0035$ & $0.00089$& $0.00089$ \\
				$\rho$ & $\si[per-mode=reciprocal]{\kilo\gram\per\cubic\meter}$ & $1056$ & $997$ & $997$ \\
				$\mathbb{R}\mathrm{e}$ & - & $500$ & $1500$ & $1500$ \\
				\hline
			\end{tabular}
		\caption{Fluid flow properties values for the FDA nozzle benchmark model (FDA) in Section~\ref{sec:FDAbench} and for the aneurysm model (ANE) from Sections~\ref{sec:aneVeri} and ~\ref{sec:aneVali}. The acronyms are further explained in Section~\ref{sec:results}.}
		\label{tab:fluid_flow_properties}
	\end{table}
    
    \par
    While no time derivative of the velocity appears in Eq.~\eqref{eq:NS_mom_cons}, the solutions $\bm u_\mathrm{ex}$ and $p_\mathrm{ex}$ may still be time dependent due to the time dependence of the inflow boundary condition $\bm g$ in Eq.~\eqref{eq:NS_bc_inl}. Moreover, since Eq.~\eqref{eq:NS_mom_cons} is steady-state, the model is assumed to hold for every $t\in [0,T]$.
    \par
    In view of the construction of loss functions that reflect the physics expressed by system~\eqref{eq:NS_strong}, we proceed by defining the Navier-Stokes residuals, for every $t\in[0,T]$, as in \cite{raissi2017physics}, for any $\bm u:\Omega\times[0,T]\to\mathbb R^3$ and $p:\Omega\times[0,T]\to\mathbb R$ sufficiently smooth:
    \begin{subequations}
		\label{eq:NS_res}
		\begin{align}
			f_\mathrm{C}(\bm u,p) &= f_1(\bm u,p)  = \nabla\cdot \bm{u}, \quad &&\mbox{in} \ \Omega, \label{eq:mass_cons} \\
			\bm{f}_\mathrm{M}(\bm u,p) &= \left(\begin{array}{c}
				f_2(\bm u,p) \\
				f_3(\bm u,p) \\
				f_4(\bm u,p)
			\end{array}\right) = (\bm{u}\cdot\nabla)\bm{u} + \nabla p - \frac{1}{\mathbb{R}\mathrm{e}}\Delta\bm{u}, \quad &&\mbox{in} \ \Omega, \label{eq:mom}
		\end{align}
	\end{subequations}
    where in Eq.~\eqref{eq:mass_cons} is the mass conservation i.e. continuity equation residual $f_\mathrm{C}$, and in Eq.~\eqref{eq:mom} is the momentum balance residual $\bm{f}_\mathrm{M}$.
    \par
    Similarly to the residuals~\eqref{eq:NS_res} associated to the PDE, we define the residuals associated to the boundary conditions \eqref{eq:NS_bc_inl}, \eqref{eq:NS_bc_wall} and \eqref{eq:NS_bc_out} as:

    \begin{subequations}
    \label{eq:bc_residuals}
	\begin{align}
		\bm{f}_{{\mathrm{BC},1}}(\bm u,p) &= \bm u - \bm g, \quad &&\mbox{on} \ \Gamma_\mathrm{inl},
		\label{eq:NS_bc1}\\
		\bm{f}_{{\mathrm{BC},2}}(\bm u,p) &= \bm{u}, \quad &&\mbox{on} \ \Gamma_\mathrm{wall},
		\label{eq:NS_bc2}\\
		\bm{f}_{{\mathrm{BC},3}}(\bm u,p) &= \left(\frac{1}{\mathbb{R}\mathrm{e}}\nabla{\bm{u}_\mathrm{ex}} - p_\mathrm{ex}\bm I\right)\bm n, \quad &&\mbox{on} \ \Gamma_\mathrm{out}.
		\label{eq:NS_N}
	\end{align}
    \end{subequations}
    While the residual associated to the PDE can be evaluated in any point $(\bm x,t)$ of the spatio-temporal domain $\Omega\cross[0,T]$, the residuals associated with the boundary conditions are only defined on the corresponding subsets of the boundary of the spatio-temporal domain. Consequently, $\bm{f}_{{\mathrm{BC},1}}$ can only be evaluated in $\Gamma_\mathrm{inl}\times[0,T]$, $\bm{f}_{{\mathrm{BC},2}}$ in $\Gamma_\mathrm{wall}\times[0,T]$ and $\bm{f}_{{\mathrm{BC},3}}$ in $\Gamma_\mathrm{out}\times[0,T]$. The residuals are well defined for all $\bm u$ and $p$ sufficiently smooth. Moreover, for $\bm u_\mathrm{ex}$ and $p_\mathrm{ex}$ all the residuals are identically zero since this pair satisfies the strong Navier-Stokes equations~\eqref{eq:NS_strong}.
    \par
    Once defined the residuals associated to the PDE~\eqref{eq:NS_res} and the boundary conditions~\eqref{eq:bc_residuals}, we introduce the physics-based loss functions, a key characteristic of PINNs. We aim at finding ANNs $\bm u(\bm \theta)\, : \, \Omega\times[0,T]\to\mathbb R^3$ and $p(\bm \theta):\Omega\times[0,T]\to\mathbb R$, where $\bm{\theta}$ are the parameters of the ANN, such that they are approximations of the exact solutions of system~\eqref{eq:NS_strong} $\bm u_\mathrm{ex}$, $p_\mathrm{ex}$, that is such that
    \begin{equation*}
        \bm{u}(\bm{x}, t; \bm{\theta})  \simeq   \bm{u}_\mathrm{ex}(\bm{x}, t),\ p(\bm{x}, t; \bm{\theta})  \simeq   p_\mathrm{ex}(\bm{x}, t),\ \forall (\bm{x}, t)\in \Omega\times[0,T].
    \end{equation*}
    The parameters of the ANNs, corresponding to the vector $\bm \theta$, can be found by solving the following minimization problem:
	\begin{equation}
		\label{eq:param_minimization}
		\bm{\theta} = \underset{\Tilde{\bm{\theta}}\in\Theta}{\argmin}\mathcal{L}(\Tilde{\bm{\theta}}),
	\end{equation}
    where $\mathcal{L}(\Tilde{\bm{\theta}})$ is the total loss function evaluated on a particular parametrization of the ANN given by $\Tilde{\bm{\theta}}$, and $\Theta$ is the parameter space in which the minimization is performed. 
    \par
    The PINN method is based on defining the total loss function as:
    \begin{equation}
        \label{eq:loss}
        \mathcal{L}(\Tilde{\bm{\theta}}) = \mathcal{L}_{\mathrm{PDE}}(\Tilde{\bm{\theta}}) + \mathcal{L}_{\mathrm{BC}}(\Tilde{\bm{\theta}}) +
		\mathcal{L}_{\mathrm{train}}(\Tilde{\bm{\theta}}),
    \end{equation}
    where  $\mathcal L_\mathrm{PDE}$ and $\mathcal{L}_\mathrm{BC}$ are the loss functions associated to the physics, and $\mathcal{L}_{\mathrm{train}}$ is a suitably defined data-driven contribution to the total loss function.
    In order to define the physics-based loss functions $\mathcal L_\mathrm{PDE}$ and $\mathcal{L}_\mathrm{BC}$, associated to the residuals \eqref{eq:NS_res} and \eqref{eq:bc_residuals}, we introduce  discrete (finite-dimensional) subsets of the spatio-temporal domain, made of the so-called collocation points, ${\Gamma}_{\mathrm{BC},1}\subset{\Gamma}_\mathrm{inl}\cross[0,T]$ with $|{\Gamma}_{\mathrm{BC},1}|=N_{\mathrm{BC},1}$, ${\Gamma}_{\mathrm{BC},2}\subset{\Gamma}_\mathrm{wall}\cross[0,T]$ with $|{\Gamma}_{\mathrm{BC},2}|=N_{\mathrm{BC},2}$, ${\Gamma}_{\mathrm{BC},3}\subset{\Gamma}_\mathrm{out}\cross[0,T]$ with $|{\Gamma}_{\mathrm{BC},3}|=N_{\mathrm{BC},3}$, and ${\Omega}_\mathrm{PDE}\subset\Omega\cross[0,T]$ with $|\Omega_\mathrm{PDE}|=N_\mathrm{PDE}$.
    \par
    For any family of sufficiently smooth functions $\bm u(\bm\theta)$ and $p(\bm\theta)$ parametrized by $\bm\theta$, we define the loss function associated to the PDE evaluated at $\bm\theta$ as:
    \begin{equation}
		\label{eq:loss_PDE}
		\mathcal{L}_\mathrm{PDE}(\bm \theta) = \sum_{i=1}^{N_\mathrm{PDE}}\sum_{j=1}^{4}\, \frac{\lambda_{\mathrm{PDE},j}}{N_{\mathrm{PDE}}}\, \left( \left|f_j(\bm{u}_{i}(\bm{\theta}), p_{i}(\bm{\theta})\right)|^2 + \lambda_{\mathrm{PDE},\mathrm{MSLE}}\, \log\left|f_j(\bm{u}_{i}(\bm{\theta}), p_{i}(\bm{\theta}))+1\right|^2\right),
	\end{equation}
    where the subscript $i$ indicates that the functions $\bm u(\bm \theta)$ and $p(\bm \theta)$ are evaluated at the $i$-th PDE collocation point $(\bm x_i,t_i)\in\Omega_\mathrm{PDE}$, that is $\bm{u}_{i}(\bm{\theta})=\bm u(\bm x_i,t_i;\bm\theta)$ and $p_{i}(\bm{\theta})= p(\bm x_i,t_i;\bm\theta)$. We remark that the relative weights $\left\{\lambda_{\mathrm{PDE},j}\right\}_{j=1}^{4}$ assigned to each residual term are in general different from each other and that the choice of the collocation points is a priori arbitrary. These aspects allow for a larger flexibility when defining the loss function~\eqref{eq:loss_PDE}~\cite{bischof_multi-objective_2025}. Our choice of adding the mean squared logarithmic error (MSLE) term in~\eqref{eq:loss_PDE} in addition to the mean squared error (MSE) is meant to enable accurate flow approximation in regions where the physics residuals have small magnitudes. Similarly, we define the loss function associated to the boundary conditions as:
    \begin{equation}
			\mathcal{L}_{\mathrm{BC}}(\bm \theta) = \sum_{j=1}^{3}\sum_{i=1}^{N_{\mathrm{BC},j}}\, \frac{\lambda_{\mathrm{BC},j}}{N_{\mathrm{BC},j}} \, \left|\bm f_{\mathrm{BC},j}(\bm{u}_{j,i}(\bm{\theta}), p_{j,i}(\bm{\theta}))\right|^2, \label{eq:loss_BC}
    \end{equation}
    where the double subscript $j_i$ indicates that the functions $\bm u(\bm \theta)$ and $p(\bm \theta)$ are evaluated at the $i$-th collocation point of the $j$-th boundary set of collocation points $(\bm x_{j,i},t_{j,i})\in\Gamma_{\mathrm{BC},j}$. Once again, the weights $\left\{\lambda_{\mathrm{BC},j}\right\}_{j=1}^{3}$ associated to the different boundary conditions are in general different from each other.
    \par
    The definition of the physics-based loss functions~\eqref{eq:loss_PDE} and~\eqref{eq:loss_BC} is central to the formulation of the PINN method. In the following,  we will consider $\mathcal{L}_\mathrm{train}$ to be expressed as:
    \begin{equation}
        \label{eq:training_loss}
        \mathcal{L}_\mathrm{train}({\bm \theta}) = \mathcal{L}_{\bm u}({\bm \theta}) + \mathcal{L}_{\hat{\bm u}}({\bm \theta}),
    \end{equation}
    where $\mathcal{L}_{\bm u}$ is the loss term of the velocity with respect to the data, and $\mathcal{L}_{\hat{\bm u}}$ is the loss term with respect to the velocity direction unitary vector data. The data are high-fidelity known values of the solution, available for a discrete subset of the spatio-temporal domain. These data are typically obtained from empirical or experimental measurements, and the type of data depends on the problem at hand. For fluid dynamics problems, the data can come in the form of the fluid velocity from 4D flow MRI acquisitions \cite{kissas2022feasibility}, like the focus of this work,  or from PIV measurements~\cite{hariharan2011multilaboratory}. Once the data are available, the entire dataset is split into training and testing datasets~\cite{goodfellow2016deep, chollet2021deep}. The training dataset is used to define the loss function $\mathcal{L}_\mathrm{train}$ in Eq.~\eqref{eq:loss}. In particular, given the $N_\mathrm{train}$ training data points comprising of the velocity data $\left\{\bm u_{\mathrm{ref},i}\right\}_{i=1}^{N_\mathrm{train}}$ at $N_\mathrm{train}$  different points of the spatio-temporal domain, the data-driven loss terms in~\eqref{eq:training_loss} are defined as:
    \begin{subequations}
		\begin{align}
			\mathcal{L}_{\bm u}(\bm \theta) &= \begin{aligned}[t]\sum_{i=1}^{N_\mathrm{train}}\sum_{j=1}^3\, \frac{\lambda_{(\bm u)_j}}{N_\mathrm{train}}\, \left(|(\bm u_{i}(\bm{\theta}))_j - (\bm u_{\mathrm{ref},i})_j|^2 + \lambda_{\bm u,\mathrm{MSLE}}\, \left| \log\left( \frac{(\bm u_{i}(\bm{\theta}))_j+1}{(\bm u_{\mathrm{ref},i})_j+1} \right) \right|^2\right),
			\end{aligned} \label{eq:loss_vel}\\
			\mathcal{L}_{\hat{\bm u}}(\bm \theta) &= \sum_{i=1}^{N_{\mathrm{train}}}\, \frac{\lambda_{\hat{\bm u}}}{N_{\mathrm{train}}}\, \left(|\hat{\bm{u}}_i(\bm{\theta}) - \hat{\bm{u}}_{\mathrm{ref},i}|^2 + \left|\log\left(\frac{\hat{\bm{u}}_i(\bm \theta)+\bm{1}}{ \hat{\bm{u}}_{\mathrm{ref},i}+\bm{1}}\right)\right|^2\right). \label{eq:loss_dir}
		\end{align}
		\label{eq:loss_train_comps}
	\end{subequations}
   The terms in $\mathcal{L}_{{\bm u}}$~\eqref{eq:loss_vel} correspond to the MSE and MSLE of the velocity vector, where the subscript $j$ in~\eqref{eq:loss_vel} indicates the $j$-th Cartesian component of the velocity vector.
   On the other hand, by introducing the loss term $\mathcal{L}_{\hat{\bm u}}$~\eqref{eq:loss_dir} related to the velocity vector field direction $\hat{\bm u}$, with $\bm 1 \in \mathbb R^3$, we enhance the velocity direction information obtained from the training dataset, even when the velocity magnitude is small~\cite{regazzoni_universal_2022}. To avoid problems with zero-magnitude vectors, we compute the unit vector $\hat{\bm u}$ as:
    \begin{equation*}
		\hat{\bm{u}} = \frac{\bm{u}}{\max{\left\{|\bm{u}|_2, \varepsilon\right\}}}, \quad\text{with} \quad \varepsilon = 10^{-2}.
	\end{equation*}  
    We remark that we do not use pressure data for the ANN training, since we aim to reconstruct it Eq.~\eqref{eq:mom}; this choice is consistent with 4D flow MRI data, where only velocity measurements are available.
    \par
    We train the PINN, consisting of solving the minimization problem~\eqref{eq:param_minimization}, either by using a stochastic gradient descent algorithm, such as the ADAM optimizer~\cite{goodfellow2016deep, kingma2014adam}, for large datasets, or by using a quasi-Newton optimization algorithm such as L-BFGS-B~\cite{zhu1997algorithm} in the case of small datasets~\cite{raissi2017physics}.
    \par
    Finally, the accuracy of the ANN approximations of the velocity and the pressure is evaluated by using the remaining $N_\mathrm{test}$ data comprising the testing dataset. To this end, analogously to Eq.~\eqref{eq:loss_vel}, we define the velocity and pressure test loss functions, $\mathcal{L}^\mathrm{test}_{\bm{u}}$ and $\mathcal{L}^\mathrm{test}_p$ respectively:
    \begin{subequations}
		\begin{align}
			\mathcal{L}^\mathrm{test}_{\bm{u}}(\bm{\theta}) &= \frac{1}{N_{\mathrm{test}}} \sum_{i=1}^{N_{\mathrm{test}}} |\bm{u}_i(\bm{\theta}) - \bm{u}_{\mathrm{ref},i}|^2, \label{eq:loss_test_u}\\ \mathcal{L}^\mathrm{test}_{p}(\bm{\theta}) &= \frac{1}{N_{\mathrm{test}}} \sum_{i=1}^{N_{\mathrm{test}}} |p_i(\bm{\theta}) - p_{\mathrm{ref},i}|^2. \label{eq:loss_test_p}
		\end{align}
		\label{eq:loss_test}
	\end{subequations}
    The test data loss function $\mathcal{L}^\mathrm{test}_p$ is well-defined only when pressure data are available.
    \section{Results}
    \label{sec:results}
    In this section we numerically verify and validate our computational pipeline based on the PINN methodology presented in Section~\ref{sec:methodology} on two test cases. To this end, we employ both real measurement data and FEM generated in silico data.
    \par
    Section~\ref{sec:FDAbench} describes the FDA nozzle benchmark, whereas in Section~\ref{sec:veriFDA}, the FDA-CFD cases address the in silico verification of different PINN configurations by training and testing FEM data, whereas in Section~\ref{sec:valiFDA} the FDA-EXP cases are devoted to experimental validation by training and testing the PINNs on PIV measurement data. On the other hand, Section~\ref{sec:aneurysm} focuses on the fluid flow in an aneurysm geometry. In Section~\ref{sec:aneVeri} the ANE-CFD networks are trained on an in silico FEM generated dataset of fluid flow through an aneurysm geometry, and in Section~\ref{sec:aneVali} the ANE-4DMRI cases are trained using experimental in vitro data.
	
    \subsection{FDA nozzle benchmark}
	\label{sec:FDAbench}
\begin{figure}[t]
    {
    \centering
	\includegraphics[width=0.9\textwidth]{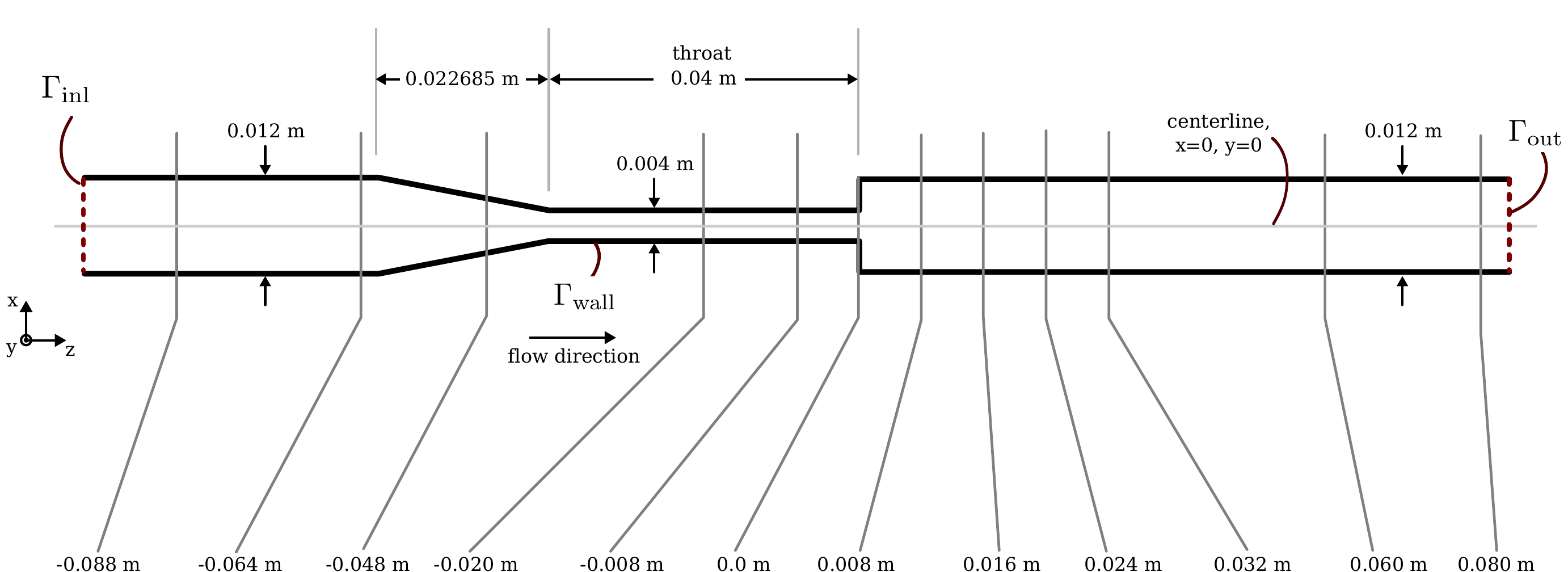}
	\caption[]{FDA nozzle domain geometry, characteristic dimensions with sketch of inlet boundary $\Gamma_\mathrm{inl}$, outlet boundary $\Gamma_\mathrm{out}$, and wall $\Gamma_\mathrm{wall}$. Vertical lines in gray are cross-sections where experimental data are available. }
    \label{fig:FDAgeo}
    }
	\end{figure}

	First, we apply the PINN method to the laminar case of the FDA nozzle benchmark~\cite{hariharan2011multilaboratory, NCIHub43}. We consider a laminar flow through an axially symmetric nozzle, whose geometry and dimensions are represented in Figure~\ref{fig:FDAgeo}. First, in Section~\ref{sec:veriFDA}, we verify the proposed PINN method by performing a CFD simulation of our problem and by using a subset of the generated in silico data to train and test the PINN. Then, in Section~\ref{sec:valiFDA}, we validate the method by employing real experimental data from PIV measurements~\cite{hariharan2011multilaboratory, NCIHub43}.
	
	\subsubsection{Numerical verification against in silico data}
	\label{sec:veriFDA}
	We discuss the verification of the PINN method for the FDA nozzle benchmark by training and testing against FEM generated in silico data. In Section~\ref{sec:veriSetFDA} we describe the setup for our problem, focusing first on data generation, then on the training of the PINNs. Following this, in  Section~\ref{sub:test1}, we present a sensitivity study of the FDA-CFD PINN output by considering different definitions of the loss function.
    \par
    In particular we study 7 different PINN configurations, named FDA-CFD and numbered 1 through 7, using 5 different tests. The exact configurations are reported in Table~\ref{tab:wstd}. The tests are designed as follows:
    \begin{itemize}
        \item In Test 1, we examine the effects of modifying the weights of the radial components in the velocity. This is done by comparing FDA-CFD-1 and 2, where FDA-CFD-2 has a higher weights on the radial components of the velocity in the data-driven term;
        \item In Test 2, we examine the effect of assigning different weight on terms coming from the PDE-based loss function on the velocity reconstruction accuracy. To this end, we compare FDA-CFD-2 and 3, by forcing the residuals in the latter to be of the same order of magnitude;
        \item In Test 3, by introducing in FDA-CFD-4 the MSLE, we improve the reconstruction of the centerline pressure;
        \item In Test 4, we test the centerline pressure and velocity accuracy reconstruction depending on the PINN size, by using twice the number of neurons per layer in FDA-CFD-5 with respect to FDA-CFD-4;
        \item In Test 5, we test the performance of FDA-CFD-6 and FDA-CFD-7 by including the velocity-direction dependent loss function in differently sized networks
    \end{itemize}
	\paragraph{Setup}
	\label{sec:veriSetFDA}
    To generate data for the training and testing of the PINNs, we first carry out CFD simulations. We perform the simulations using \texttt{life$^\texttt{X}$}~\cite{africa2022lifex,africa_lifex-cfd_2024}, a high-performance solver of multiphysics and multiscale differential models developed at the MOX laboratory of Politecnico di Milano. The software is based on  the \texttt{deal.II} finite element core~\cite{arndt_dealii_2021,africa_dealii_2024} and mainly targets cardiocirculatory applications.
    \par
	We consider a laminar flow in the nozzle characterized by a Reynolds number at the throat section of $\mathbb{R}\mathrm{e} = 500$. As the characteristic velocity $U$ we consider the averaged velocity over the throat cross-section. We set the fluid flow parameters, reported in Table~\ref{tab:fluid_flow_properties}, first column, in accordance with the FDA nozzle benchmark~\cite{stewart2012assessment}.
    \par
    In accordance with Eq.s~\eqref{eq:NS_bc_inl}, \eqref{eq:NS_bc_out} and \eqref{eq:NS_bc_wall}, we impose a parabolic velocity profile with a volumetric flow rate of $Q$ at the inlet section $\Gamma_{\mathrm{inl}}$, a homogeneous Neumann boundary condition at the outlet section $\Gamma_{\mathrm{out}}$, and a no-slip condition at the lateral wall of the nozzle $\Gamma_{\mathrm{wall}}$. As the final time of the simulation we set $T=2~\si{\second}$.
    \par
	For the space discretization we use trilinear continuous nodal finite elements (FE) for both velocity and pressure, stabilized with the Variational-Multiscale Large Eddy Simulation (VMS-LES) method~\cite{forti2015semi, zingaro2021hemodynamics}. We bring the solution to a steady-state regime by using a smooth ramp lasting $0.01~\si{\second}$ starting from a null velocity initial condition. For the time integration of the unsteady solver we use the Backward Euler method with a semi-implicit treatment of the nonlinear term \cite{forti2015semi}. We consider two hexahedral mesh refinement levels, a coarse level and a fine level. We report the simulation parameters in Table~\ref{table:FEM}. Simulations were ran on 20 cores on the \texttt{gigat} cluster of the Department of Mathematics, Politecnico di Milano~\cite{polimiHardware}. Once we have generated the in silico FEM data, we perform the ANN initialization, training and testing.
    \par
    The PINN inputs are the spatial coordinates $x, y, z$, whereas the outputs are the pressure $p$ and the velocity $\bm{u}$. We use three hidden layers and hyperbolic tangent activation functions~\cite{dung_choice_2023}. These dimensions were chosen on the basis of trial and error, and while it is possible to render them hyperparameters and tune them further using automatic parameter optimization, we skip this procedure due to its significant additional computational cost. In order to address the problem of vanishing and exploding gradients during training, we use batch-normalization before each layer~\cite{geron2019hands,ioffe2015batch}. We consider the total loss function $\mathcal L$ defined in Eq.~\eqref{eq:loss}. We recall that the pressure does not enter the definition of the total loss function $\mathcal{L}$, of Eq.~\eqref{eq:loss}. Instead, the pressure field will be directly reconstructed from the velocity using Eq.~\eqref{eq:mom}.
    In order to avoid saturation of the activation function~\cite{glorot2010understanding}, we normalize the features by making spatial coordinates and flow properties dimensionless. Moreover, normalizing the outputs guarantees that they will all have the same order of magnitude~\cite{raschka2015python}.
    \par
\begin{figure}[t]
		\centering
		\subfloat[\label{fig:trainPINNs}]{
			\includegraphics[width=0.4\textwidth]{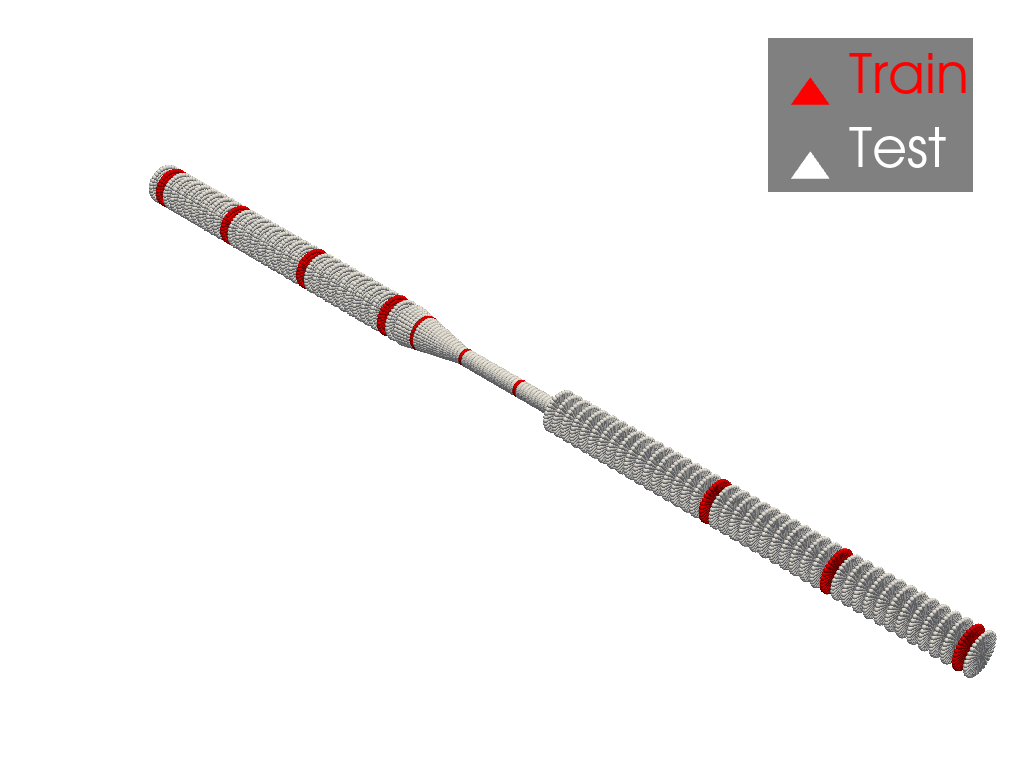}
		}
		\quad
		\subfloat[
		\label{fig:col_bc}]{
			\includegraphics[width=0.4\textwidth]{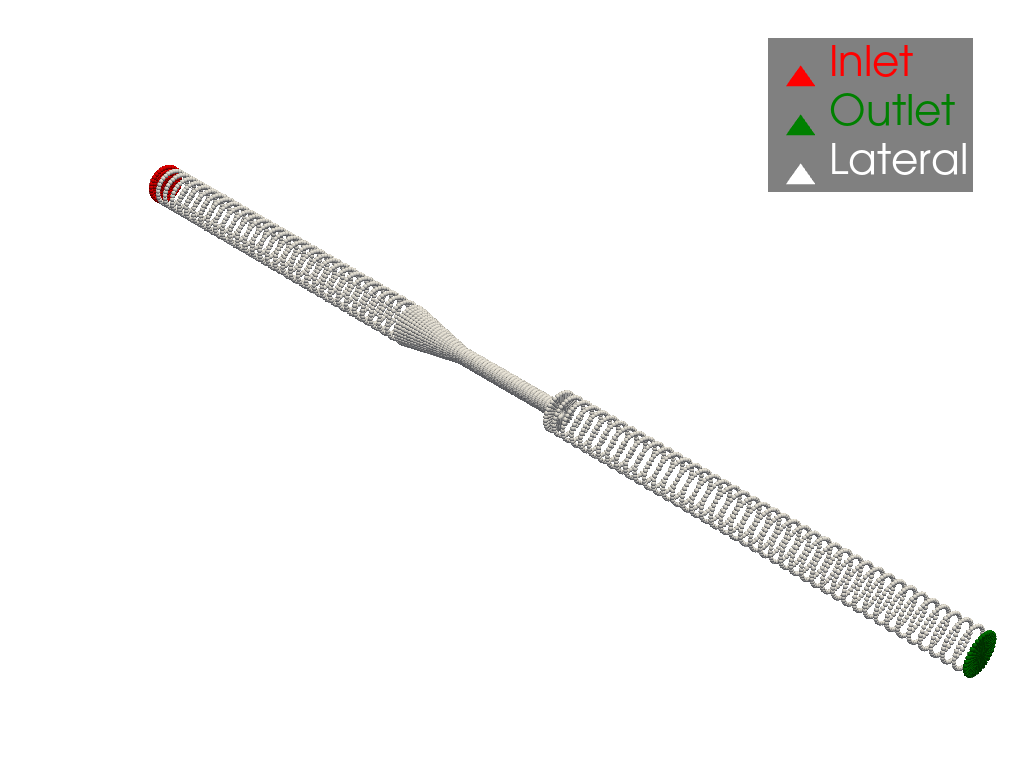}
		}
		
		\subfloat[
		\label{fig:col_less}]{
			\includegraphics[width=0.4\textwidth]{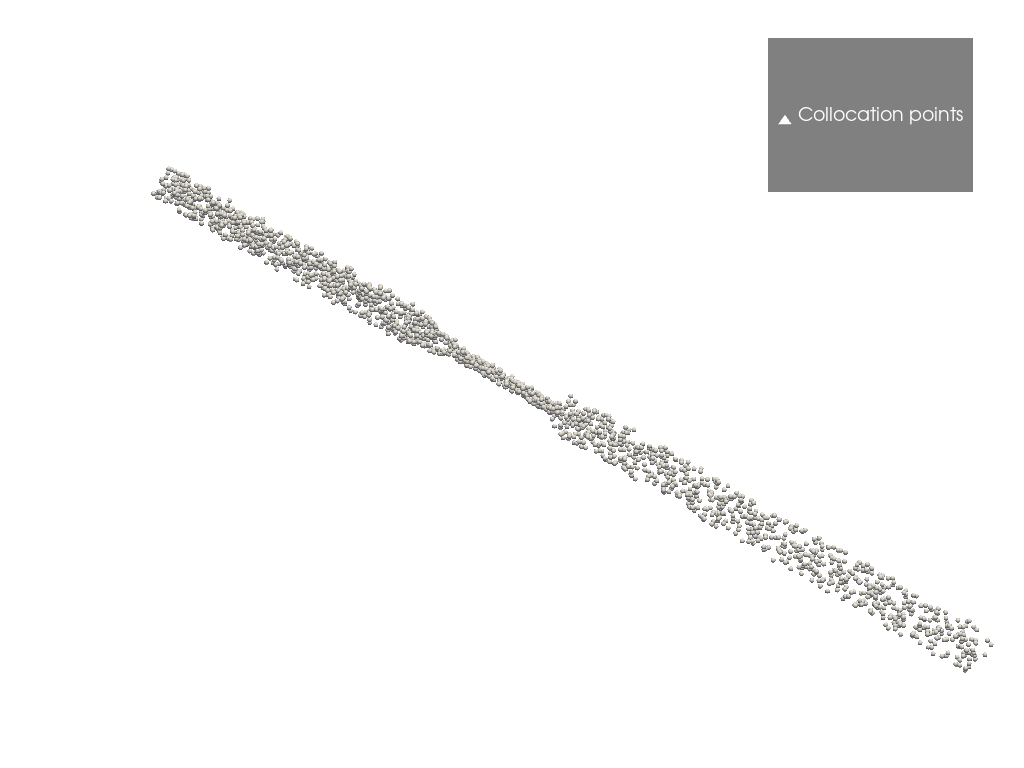}
		}
		\quad
		\subfloat[
		\label{fig:col_bet}]{
			\includegraphics[width=0.4\textwidth]{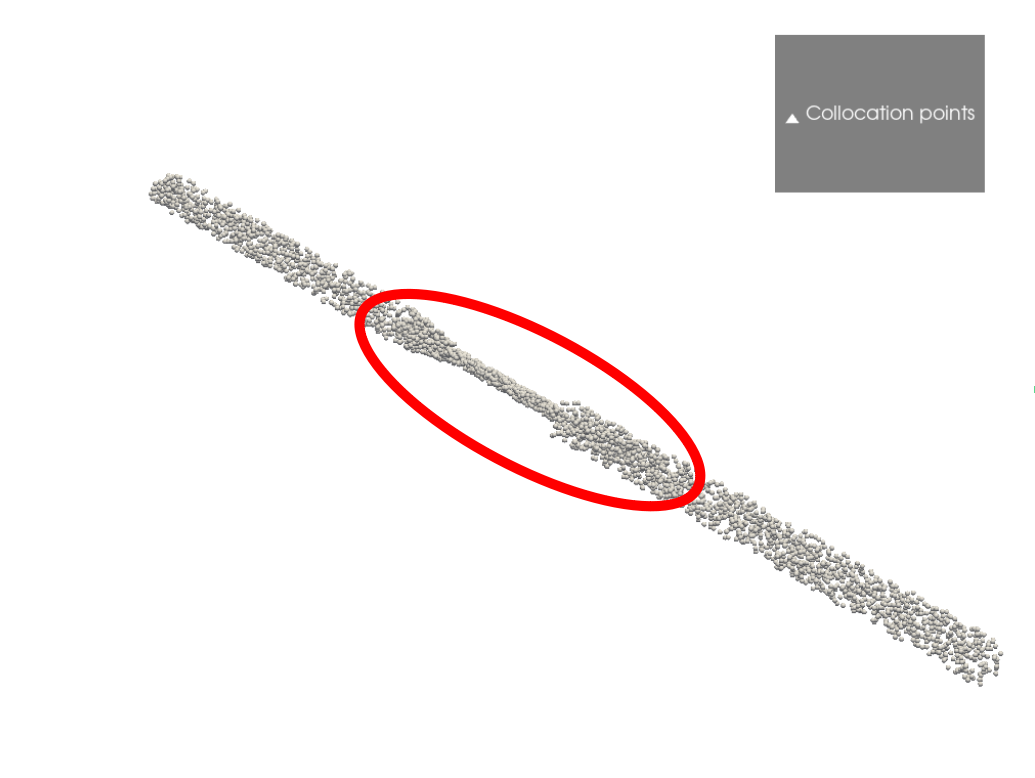}
		}
		\caption[]{Collocation points for the FDA nozzle benchmark. \textbf{(a)} Training dataset sections (red) and testing dataset sections (gray).
		\textbf{(b)} Boundary collocation points. Inlet boundary (red), outlet boundary (green) and lateral wall (gray). \textbf{(c)} PDE collocation points with $N_\mathrm{PDE}=3093$. \textbf{(d)} PDE collocation points with a higher refinement in the circled region (red), $N_\mathrm{PDE}=4341$.}
		\label{fig:col}
	\end{figure}

    We split the dataset into a training and testing subsets by selecting specific cross-sections as depicted in Figure~\ref{fig:trainPINNs}. We consider 10 training sections, with a total number of $2922$ data points. We randomly sampled $25\%$ of the data points over each section, thus obtaining a total of $N_\mathrm{train}=730$ training data points. The collocation points of the boundary portions are sketched in Figure~\ref{fig:col_bc}. We use a total of $N_{\mathrm{BC},2}=2045$ collocation points at the wall. The $N_\mathrm{PDE}$ collocation points are randomly sampled over the domain, and are depicted in Figure~\ref{fig:col_less} and Figure~\ref{fig:col_bet}, depending on the required refinement level.
	\par
    We consider different numerical setups and ANN architectures. We study the accuracy of the reconstruction of the velocity and pressure fields depending on the scaling by modifying the weights $\lambda$ associated with each loss term in Eq.~\eqref{eq:loss}. Moreover, we examine the effect of including a higher number of PDE collocation points in the convergent region, where we expect the solution to exhibit high gradients. Lastly, we study the effects of increasing the number of neurons per layer $N_\mathrm{neuron}$ on the ANN accuracy. The parameters of the loss functions used for the training of the different PINNs are reported in Table~\ref{tab:wstd}.
    \par
    All ANN initialization, training and testing was done using \texttt{TensorFlow}~\cite{tensorflow2015-whitepaper} and \texttt{Keras} \cite{chollet2015keras}. Training was performed using \num{100} epochs of the ADAM optimizer and \num{30000} epochs of L-BFGS-B for the smaller ANNs having 16 neurons per layer, whereas \num{100} epochs of the ADAM optimizer and \num{50000} epochs of L-BFGS-B were used for the larger ANNs having 32 neurons per layer. 
    \begin{table}[t]
		\centering
		\begin{tabular}{ p{5em} p{3.2em} p{3.2em} p{3.2em} p{3.2em} p{3.2em} p{3.2em} p{3.2em} }
         & \multicolumn{7}{c}{FDA-CFD}\\
        \hline
		& 1 & 2 & 3 & 4 & 5 & 6 & 7\\
		\hline
		$\lambda_{(\bm u)_1}$ & 0.91 & 100 & 100 & 100 & 100 & 0.17 & 0.17 \\
		$\lambda_{(\bm u)_2}$ & 0.91 & 100 & 100 & 100 & 100 & 0.17 & 0.17 \\
		$\lambda_{(\bm u)_3}$ &0.91& 0.91 & 0.91 & 0.91 & 0.91 & 0.17 & 0.17 \\
		$\lambda_\mathrm{BC,1}$ & 0.5  & 0.5 & 0.5 & 0.5 & 0.5 & 0.5 & 0.5 \\
		$\lambda_\mathrm{BC,2}$& 1.7  & 1.7 & 1.7 & 1.7 & 1.7 & 1.7 & 1.7 \\
		$\lambda_\mathrm{BC,3}$ & 0.5 & 0.5 & 0.5 & 0.5 & 0.5 & 0.5 & 0.5 \\
		$\lambda_\mathrm{PDE,1}$ & 2 & 2 & 2 & 3 & 3 & 3 & 3 \\
		$\lambda_\mathrm{PDE,2}$& 1.6 & 1.6 & 0.308 & 0.308 & 0.308 & 0.231 & 0.601 \\
		$\lambda_\mathrm{PDE,3}$& 1.6 & 1.6 & 0.308 & 0.308 & 0.308 & 0.231 & 0.601 \\
		$\lambda_\mathrm{PDE,4}$ & 15 & 15 & 2.88 & 2.88 & 2.88 & 2.88 & 2.88 \\
		$\lambda_{\mathrm{PDE},\mathrm{MSLE}}$ & 0 & 0 & 0 & 1 & 1& 1& 1\\
		$\lambda_{\bm u,\mathrm{MSLE}}$ & 0 & 0 & 0 & 1 & 1& 1& 1\\
		$\lambda_{\hat{\bm u}}$ & 0 & 0 & 0 & 0 & 0 & 0.1 & 0.1 \\
        $N_\mathrm{PDE}$ & 3093 & 3093 & 3093 & 3093 & 3093 & 4341 & 4341 \\
        $N_\mathrm{neuron}$ & 16 & 16 & 16 & 16 & 32 & 16 & 32 \\
		\hline
		\end{tabular}
		\caption{Parameters for the loss functions used for PINNs training in Section~\ref{sec:veriFDA}.}
		\label{tab:wstd}
	\end{table}
	  \paragraph{Test 1}
	\label{sub:test1}
    First, we compare the effects of modifying the scaling, i.e. the weights $\lambda$ associated with the data-driven loss term $\mathcal{L}_{\bm u}$ in~\eqref{eq:training_loss} on the accuracy of the reconstructed velocity. We investigate the use of larger weights for the radial velocity components with respect to the axial one in $\mathcal{L}_{\bm u}$ \eqref{eq:training_loss}, and therefore compare the two PINNs trained with the FDA-CFD-1 and the FDA-CFD-2 setup.
	\par
\begin{figure}[t]
        \centering
		\subfloat[][]{
			\includegraphics[width=0.33\textwidth]{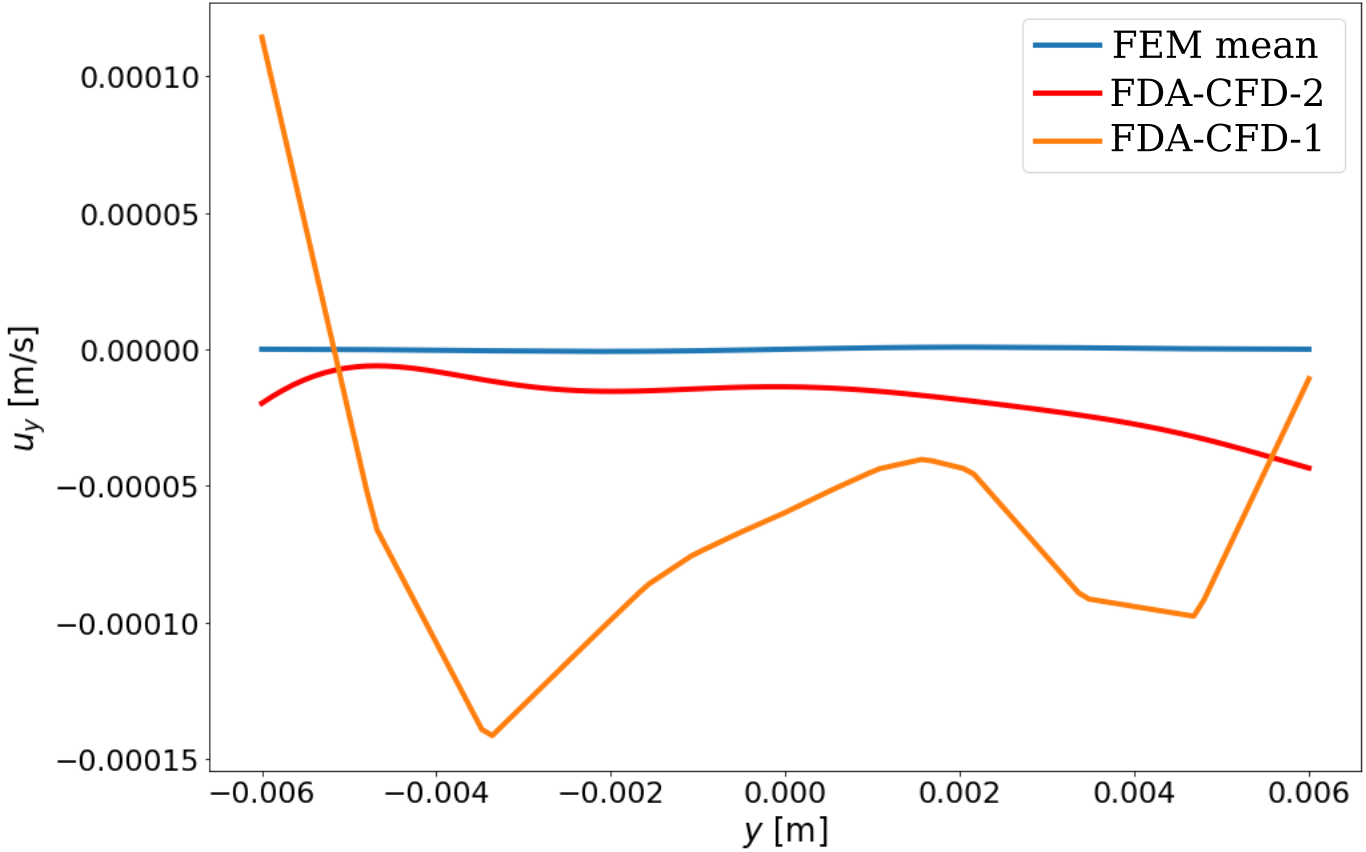}\label{fig:uycompstd}
		}\quad
		\subfloat[][]{
			\includegraphics[width=0.31\textwidth]{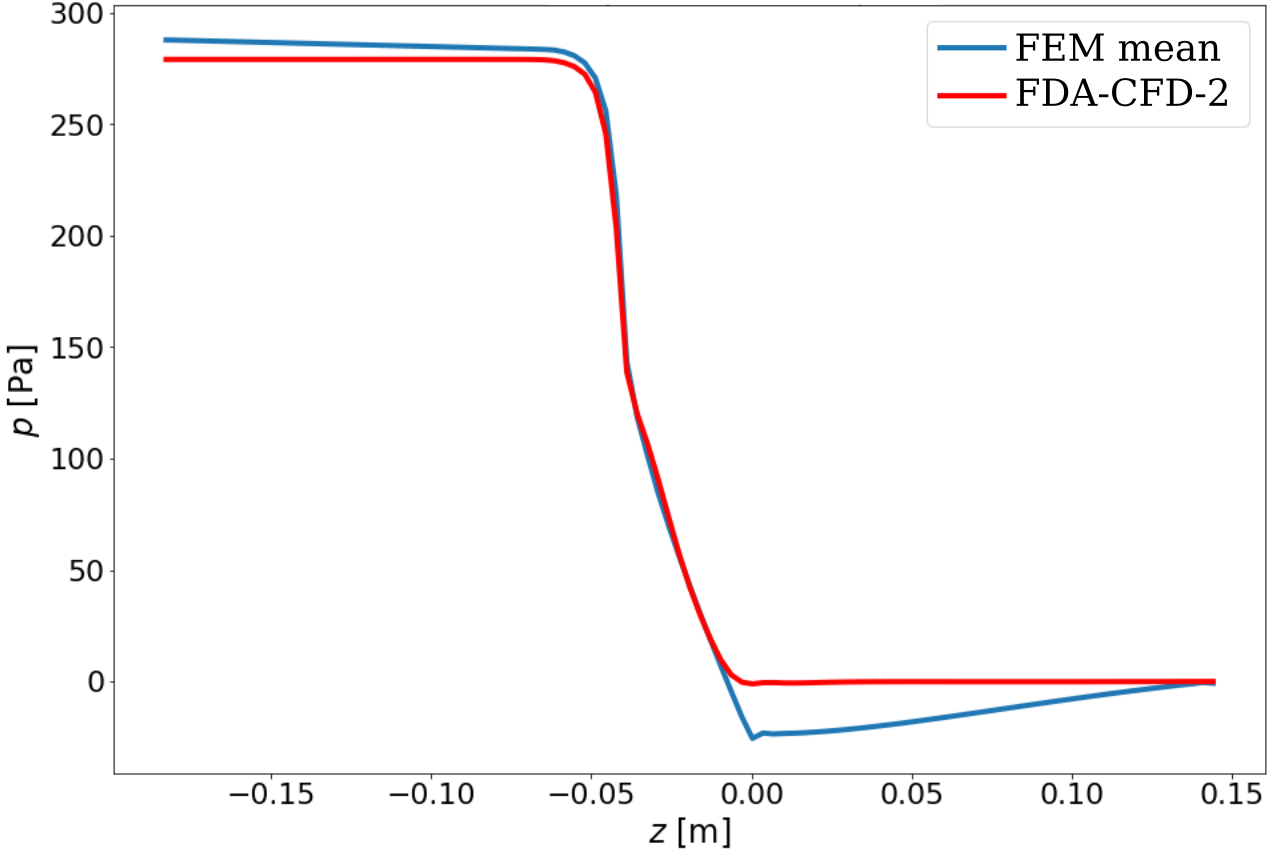}
		\label{fig:pcompstd}
		}
		\caption[]{Test 1: Comparison of approximation properties of FDA-CFD-1 and FDA-CFD-2 for low magnitude quantities and gradients. \textbf{(a)} Velocity component along the $y$-axis at $z=-0.1~\si{\meter}$. \textbf{(b)} Comparison of the centerline pressure along the $z$-axis. }
	\end{figure}

    The radial velocity components values computed by the ANN are large compared to the FEM solution which is similar to a Poiseuille flow, as reported in Figure~\ref{fig:uycompstd}. Figure~\ref{fig:uycompstd} also shows that using larger weights for the radial components in the loss function $\mathcal{L}_{\bm u}$ can mitigate this effect, as evidenced by the radial velocity reconstruction of FDA-CFD-2. We believe that the worse approximation results in regions where flow-field gradients are low are in part due to the higher relative contributions of the physics loss function $\mathcal{L}_\mathrm{PDE}$ in regions with high velocity gradients, such as the throat region.
    \par
\begin{figure}[t]
		\centering
        \subfloat[][]{
	   \includegraphics[width=0.7\textwidth]{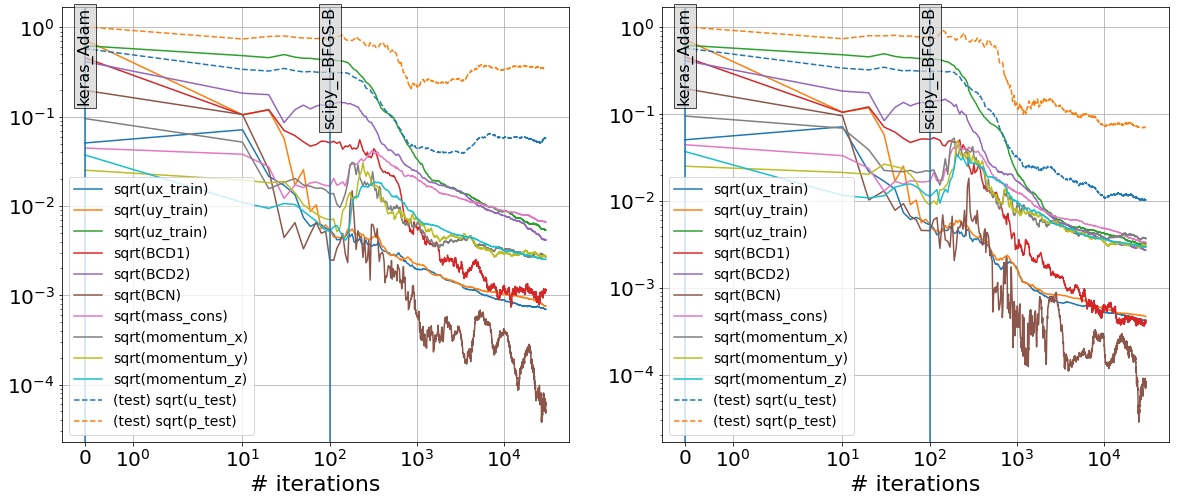}\label{fig:hist5759}
		}\\
		\subfloat[][]{
			\includegraphics[width=0.7\textwidth]{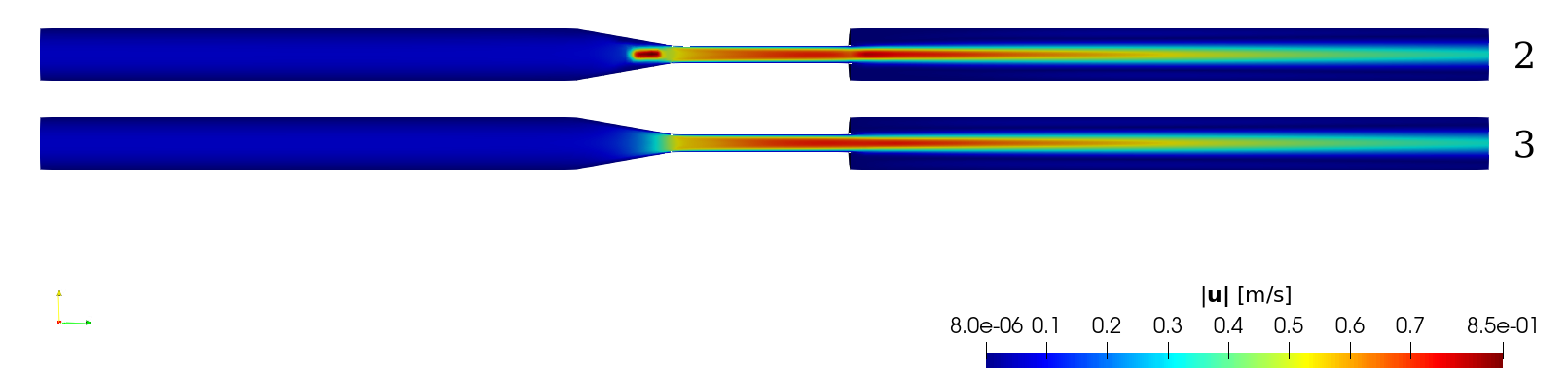}\label{fig:std_velmag}
		}
		\caption{Test 2: Comparison between FDA-CFD-2 and FDA-CFD-3 setups. \textbf{(a)} Total loss function components values during minimization. Left: FDA-CFD-2. Right: FDA-CFD-3.\textbf{(b)} Velocity magnitudes. Top: FDA-CFD-2. Bottom: FDA-CFD-3.}
		\label{fig:56comp}
	\end{figure}

    We recall that there is no training data on pressure, and it is instead reconstructed from the momentum balance in Eq.~\eqref{eq:mom}. We obtain a constant pressure before the convergent and after the diffuser, as shown in Figure~\ref{fig:pcompstd}, where we compare the pressure reconstruction obtained from the PINN with the FEM solution.
	  \paragraph{Test 2}
	\label{sub:test2}
    In this test, we compared the effects of a physically consistent scaling on the accuracy of the velocity reconstruction. We choose the weights $\left\{\lambda_{\mathrm{PDE},j} \right\}_{j=1}^4$ for FDA-CFD-3 such that the total contributions of the mass, momentum and Dirichlet no-slip residuals to $\mathcal{L}_\mathrm{PDE}$ \eqref{eq:loss_PDE} are all of the same order of magnitude. We found that this choice leads to more accurate results compared to FDA-CFD-2. We report the values of the loss functions during training for FDA-CFD-2 and FDA-CFD-3 in Figure~\ref{fig:hist5759}. In Figure~\ref{fig:hist5759}, we see that by reducing the weights for the momentum balance residuals in FDA-CFD-3, we achieve lower values of the test loss functions $\mathcal{L}_\mathrm{test}^{\bm u}$ \eqref{eq:loss_test_u} and $\mathcal{L}_\mathrm{test}^{p}$ \eqref{eq:loss_test_p} compared to FDA-CFD-2, for which the test loss functions cease to decrease after about \num{1000} iterations. In Figure~\ref{fig:std_velmag} we show the velocity magnitude for the two setups. The results of FDA-CFD-2 show a nonphysical behavior in the convergent of the nozzle. Moreover, we observed that FDA-CFD-2 has a larger dependence on the weights initialization than FDA-CFD-3. This indicates that FDA-CFD-3, for which the contributions of the physics residuals~\eqref{eq:bc_residuals} and~\eqref{eq:NS_res} are of the same order of magnitude, provides a better and more robust approximation than its counterpart, FDA-CFD-2.

\begin{figure}[t]
		\centering
		\subfloat[][]{
			\includegraphics[width=0.31\textwidth]{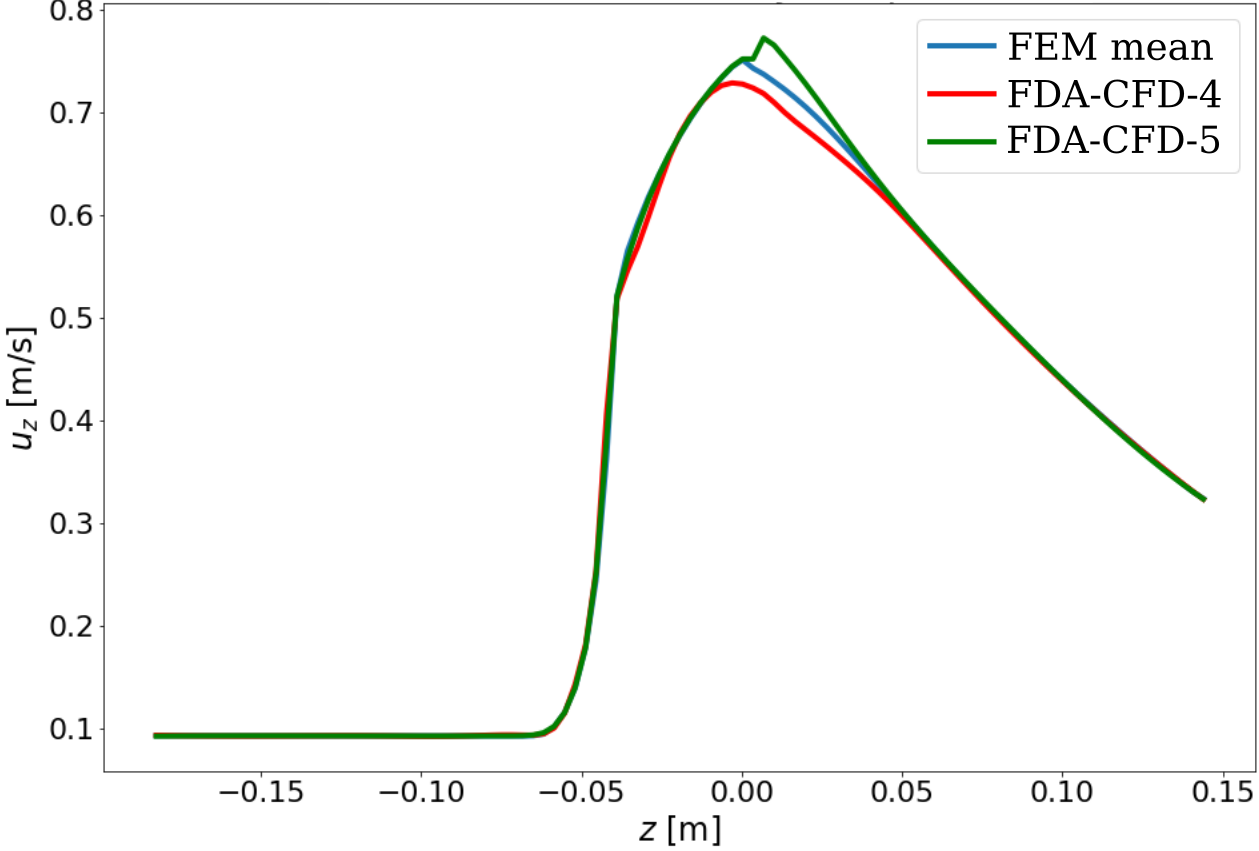}\label{fig:loguz}
		}
		\quad
        \subfloat[][]{
			\includegraphics[width=0.31\textwidth]{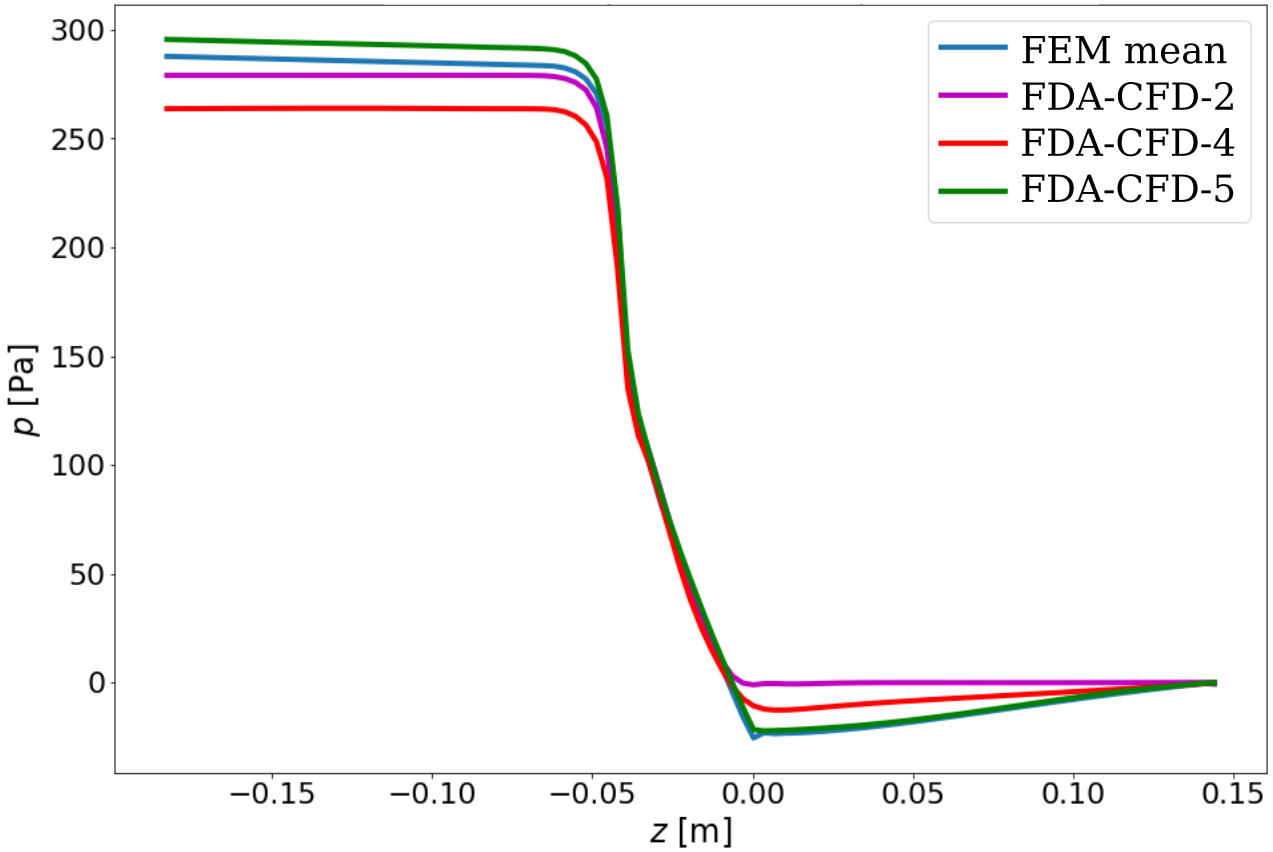}\label{fig:logp}
		}
		\caption[]{Tests 3 and 4: Comparison of computed quantities for the FEM simulation (blue) and the PINNs. \textbf{(a)} Centerline axial velocity $(\bm u)_z$. \textbf{(b)} Centerline pressure $p$.}
		\label{fig:logcomp}
	\end{figure}

	  \paragraph{Test 3}
	\label{sub:test3}
    For the FDA-CFD-4 PINN, we study the influence of the inclusion of the MSLE into the training loss function on the approximation of low magnitude data. By including the MSLE in $\mathcal{L}_\mathrm{PDE}$, we minimize more the residuals where flow-field gradients are low and we achieve a better pressure gradient reconstruction. This is shown in Figure~\ref{fig:logp}: by adding the MSLE in FDA-CFD-4, we predict a nonzero pressure gradient at the end of the duct. This is different from what we observed in FDA-CFD-2, in which we obtain a zero pressure gradient after the diffuser, which is not consistent with the FEM solution.
    
\begin{figure}[t]
		\centering
		\subfloat[][]{
			\includegraphics[width=0.9\textwidth]{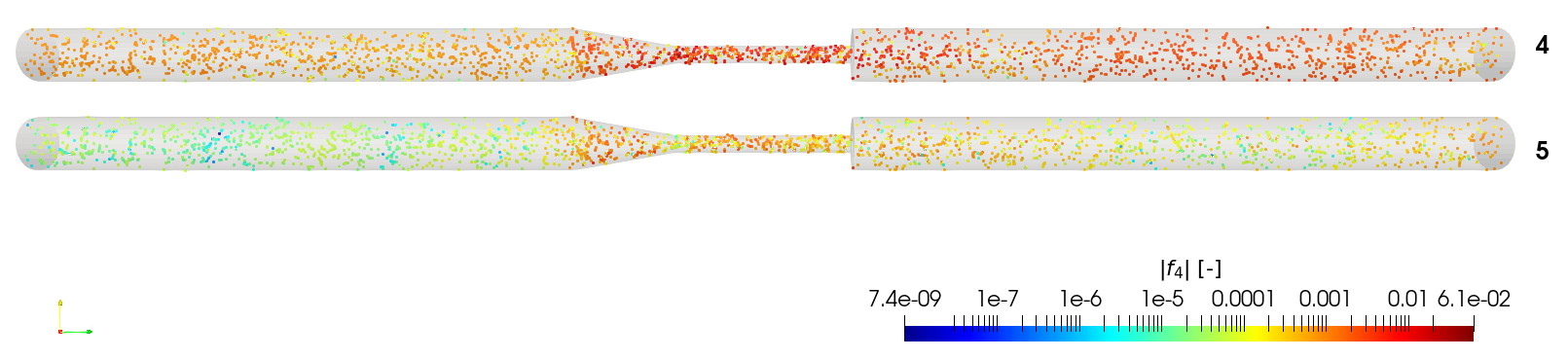}\label{fig:momzlog}
		}\quad
		\subfloat[][]{
			\includegraphics[width=0.7\textwidth]{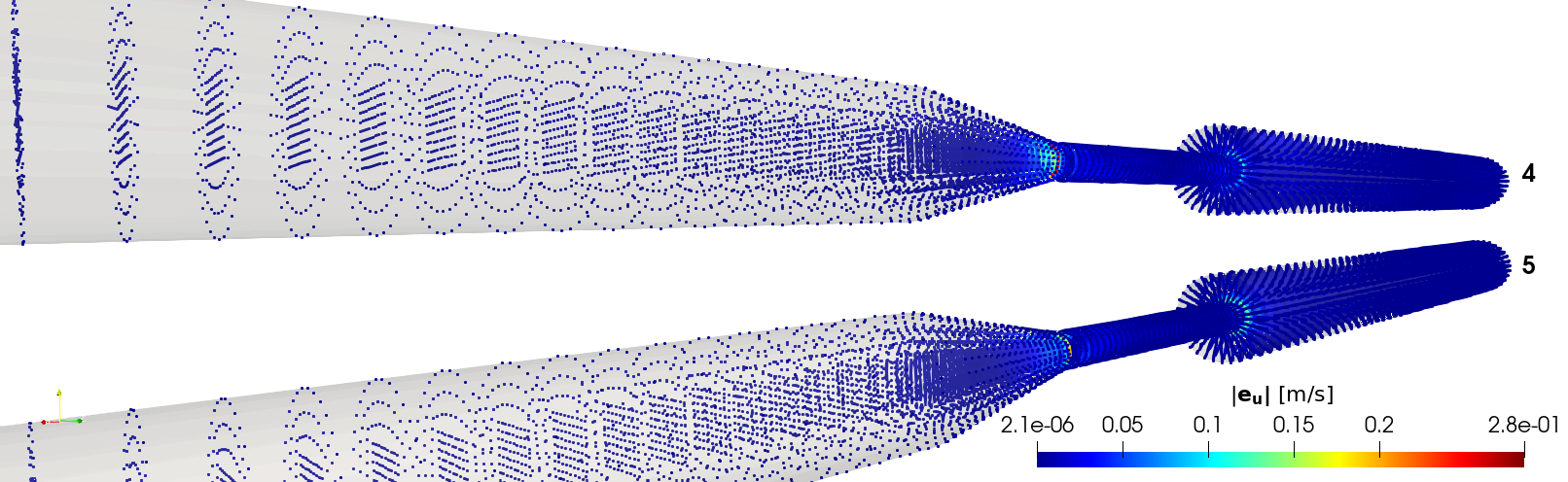}\label{fig:abslog}
		}
		\caption[]{Test 4: Pointwise comparison between FDA-CFD-4 and FDA-CFD-5. \textbf{(a)} Absolute value of the momentum balance residual along the $z$-axis, $|f_4|$, in log scale, evaluated at the collocation points. FDA-CFD-4 (top) and FDA-CFD-5 (bottom). \textbf{(b)} Velocity absolute error magnitude. FDA-CFD-4 (top) and FDA-CFD-5 (bottom).}
		\label{fig:logpts}
	\end{figure}

    \paragraph{Test 4}
	\label{sub:test4}
    We study the impact of increasing the number of neurons of the PINN on the accuracy of the pressure and velocity reconstructions. Figure~\ref{fig:logp} shows that the pressure reconstruction becomes more accurate as more neurons are used in the network. This is especially true in the region downstream the throat. In Figure~\ref{fig:momzlog}, we report the magnitude of the momentum equation residual in the $z$ direction. Indeed, for FDA-CFD-5, we obtain a lower residual compared to FDA-CFD-4, which reflects the higher accuracy when reconstructing the pressure.
    \par
\begin{figure}[t]
		\centering
		\subfloat[][]{
			\includegraphics[width=0.41\textwidth]{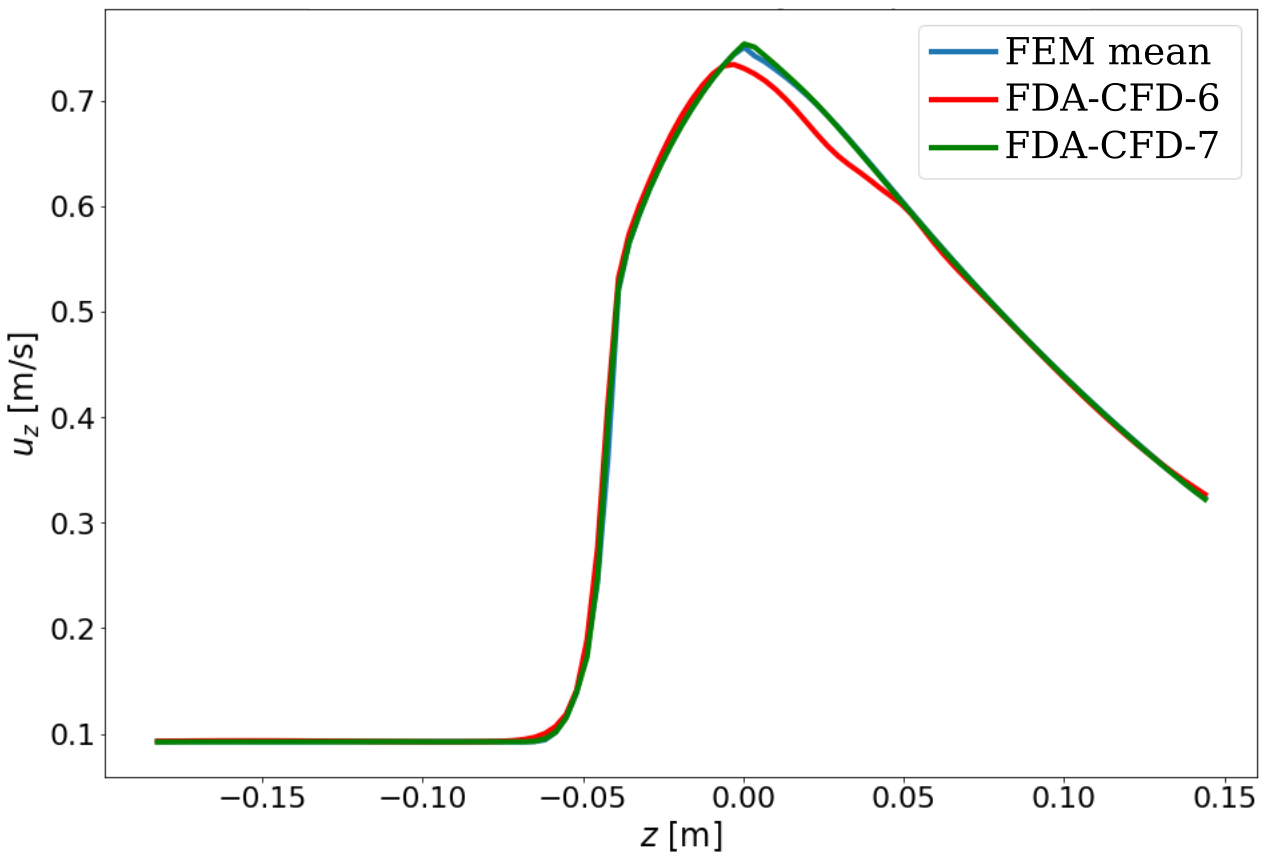}\label{fig:uzfin}
		}
		\quad
		\subfloat[][]{
			\includegraphics[width=0.41\textwidth]{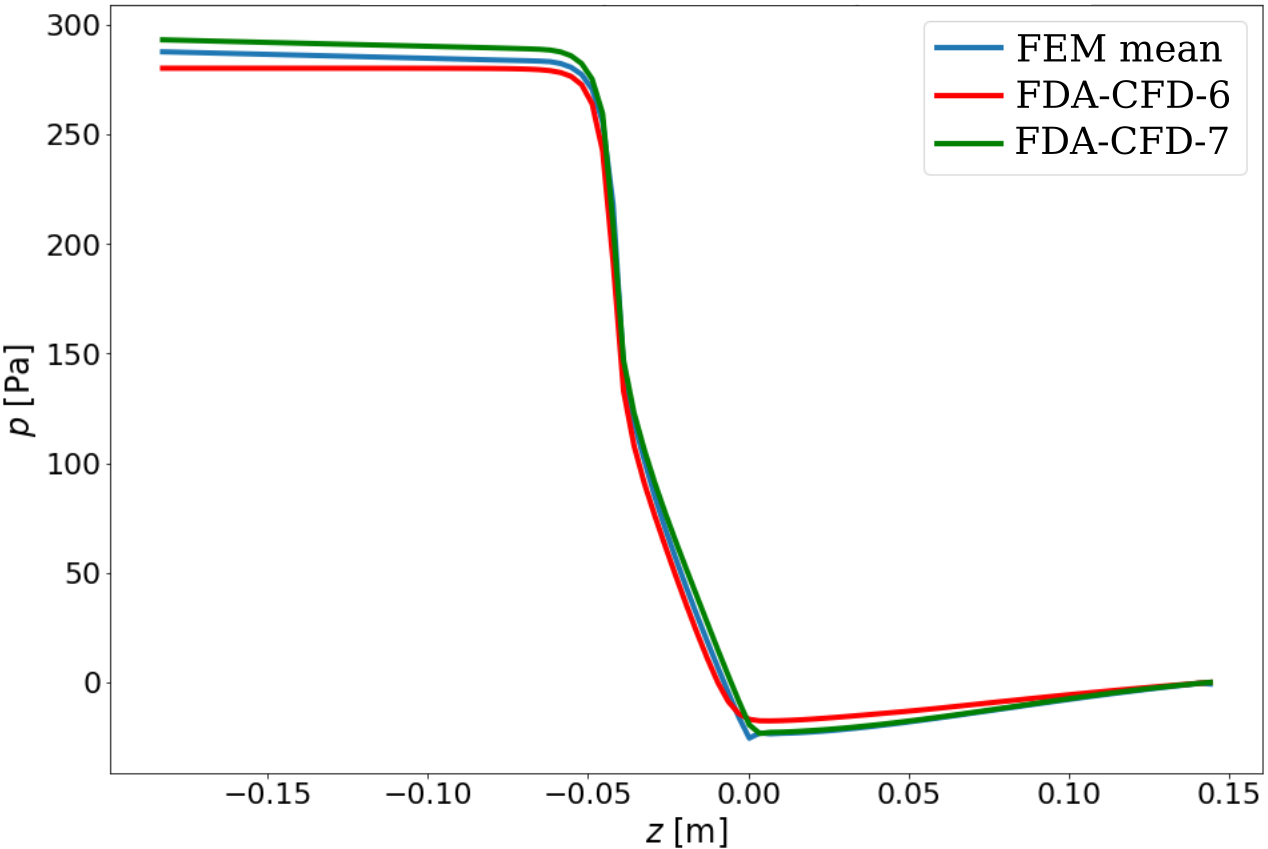}\label{fig:pfin}
		}\\
		\subfloat[][]{\includegraphics[width = 0.41\textwidth]{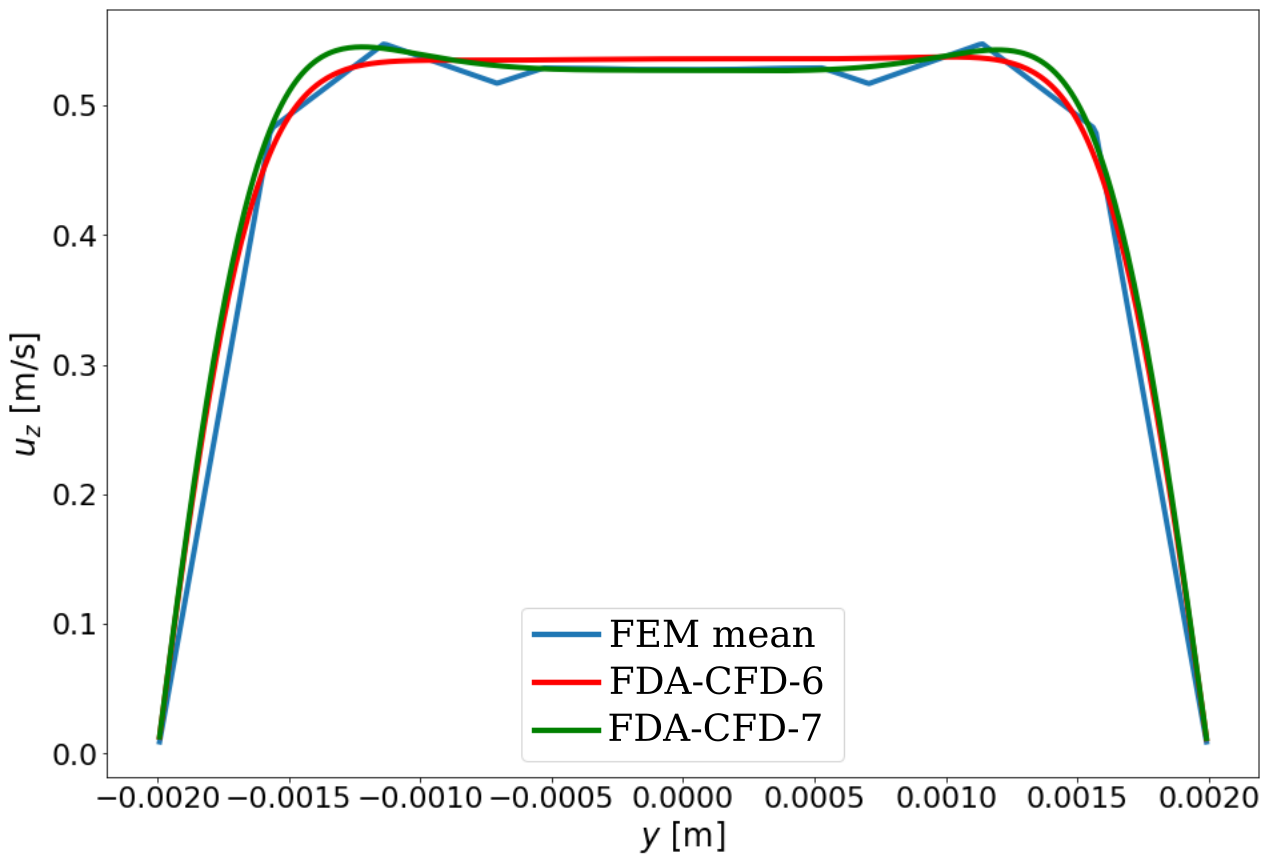}\label{fig:moroutl}}
        \quad
        \subfloat[][]{\includegraphics[width=0.41\textwidth]{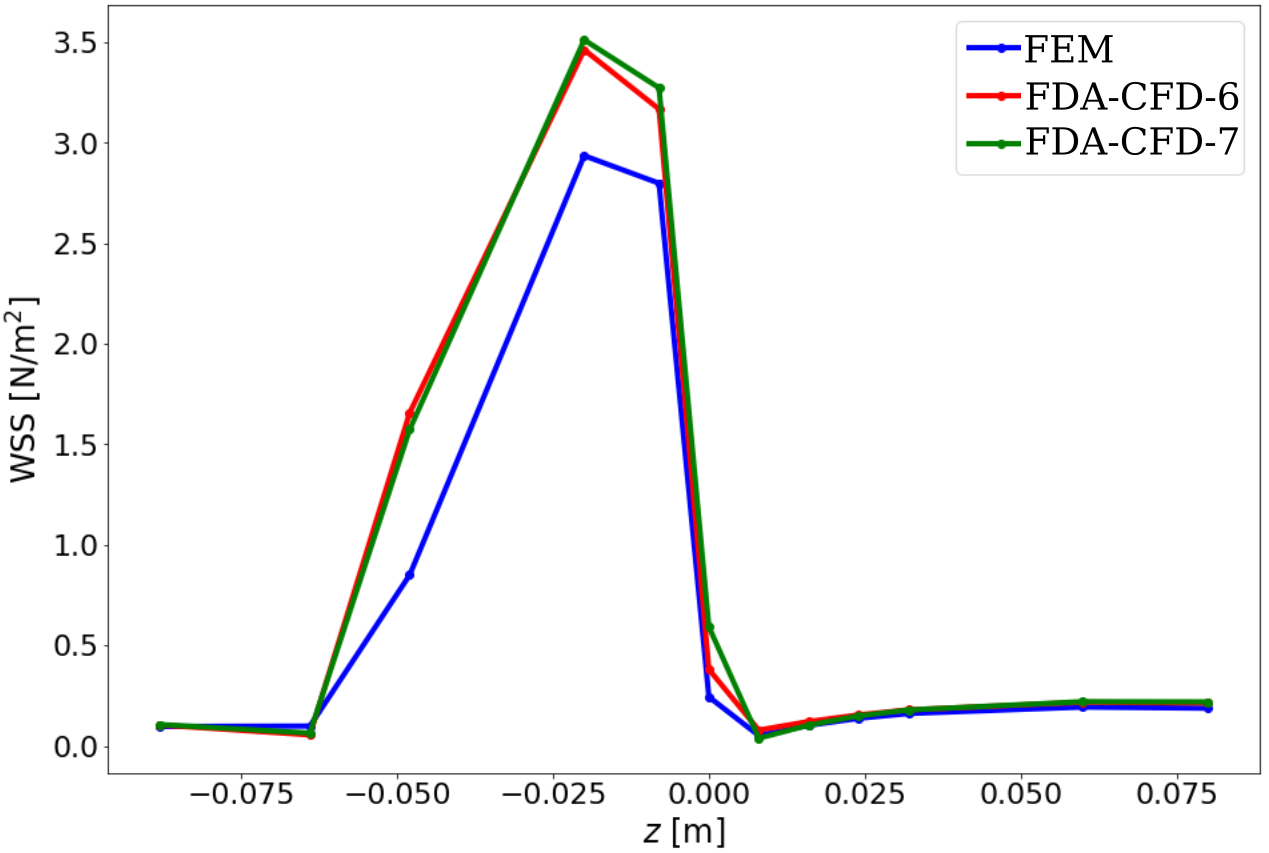}\label{fig:WSSaug}}
		\caption[]{Test 5: Comparison of computed quantities for the FEM simulation (blue) and the PINNs. \textbf{(a)} Centerline axial velocity $(\bm u)_z$. \textbf{(b)} Centerline pressure $p$. \textbf{(c)} Axial velocity  $(\bm u)_z$ distribution over training section at $z=-0.03875~\si{\meter}$. \textbf{(d)} Wall-averaged WSS magnitude comparison.}
		\label{fig:logcomp5}
	\end{figure}

    While the centerline pressure reconstruction for FDA-CFD-5 improves, the velocity $(\bm u)_z$ evaluated at the duct axis has an overshoot with respect to the FEM simulation after the diffuser, as seen in Figure~\ref{fig:loguz}. We remark that we do not give any training data section to the ANN immediately after the diffuser, as reflected in Figure~\ref{fig:trainPINNs}, so the solution there is retrieved only from the physics regularization term.
	We observe that the momentum balance residuals along the $z$-axis are larger at the end of the convergent, as shown in Figure~\ref{fig:momzlog}. The same also holds for the absolute errors of the approximated velocity with respect to the FEM simulation, reported in Figure~\ref{fig:abslog}.
    \afterpage{\clearpage}
    \paragraph{Test 5}
    \label{sub:test5}
	  In this test, we examine the effects of adding the velocity direction data in the training loss function. In Figures~\ref{fig:uzfin} and~\ref{fig:pfin} we can see the results of FDA-CFD-6 and FDA-CFD-7 for the centerline pressure and the  velocity. In Figure~\ref{fig:uzfin} we see that FDA-CFD-7 captures the peak of the centerline velocity from the FEM simulation, whereas FDA-CFD-6 underestimates it and smooths it out. We recall once more that we do not give any training data near the diffuser to the ANN, as depicted in Figure~\ref{fig:trainPINNs}. Similarly to FDA-CFD-4, FDA-CFD-6 predicts a zero pressure gradient at the beginning of the duct. Meanwhile, the pressure reconstruction in FDA-CFD-7 is coherent with the CFD simulation results. On the other hand, the higher number of parameters of the larger ANN, FDA-CFD-7, means that it overfits the numerical noise in the FEM simulations with respect to its counterpart with a lower number of parameters, FDA-CFD-6, as reflected in Figure~\ref{fig:moroutl}.
   \par
    Starting from the FEM simulation and ANN outputs, we compute the wall shear stress (WSS), a hemodynamic index commonly used for disease diagnosis~\cite{baek2009wall,caro2009discovery}, as follows:
	\begin{equation}
		\label{eq:wss}
		\epsilon = \frac{1}{2}\left(\nabla\bm{u} + \left(\nabla\bm{u}\right)^T\right),
		\qquad
		\tau = 2\mu\epsilon, \qquad
		\bm{\mathrm{WSS}} = \tau\bm{n} - (\tau\bm{n}\cdot\bm{n})\bm{n}.
	\end{equation}
	In Eq.~\eqref{eq:wss}, the second-order tensor $\epsilon$ is the symmetric part of the velocity gradient, called the shear rate tensor, the second-order tensor $\tau$ is the shear stress tensor, and $\bm{\mathrm{WSS}}$ is the vector WSS, defined on $\Gamma_\mathrm{wall}$. We compute the WSS\footnote{For the FEM simulation results we compute the WSS in postprocessing using Paraview~\cite{ahrens2005paraview}, whereas for the PINN results we use differential operators implemented in \texttt{nisaba}, that rely on automatic differentiation.} magnitude at each axial coordinate $z$ and average it over the wall. We report the computed wall-averaged WSS for FDA-CFD-6 and FDA-CFD-7 compared to the mean FEM simulation in Figure~\ref{fig:WSSaug}. Both ANNs predict a similar WSS distribution along the $z$-axis, which is larger than the one predicted using the FEM simulation. 
	\subsubsection{Validation against experimental data}
	\label{sec:valiFDA}
\begin{figure}[t]
		\centering
		\subfloat[][]{
        \includegraphics[width=0.31\textwidth]{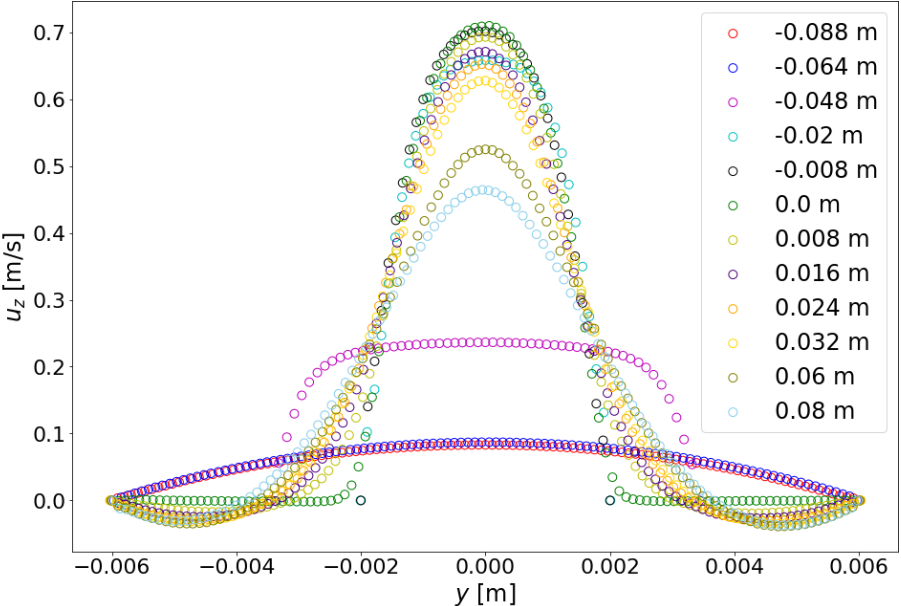}\label{fig:exp_uz}
		}
		\quad
		\subfloat[][]{
			\includegraphics[width=0.31\textwidth]{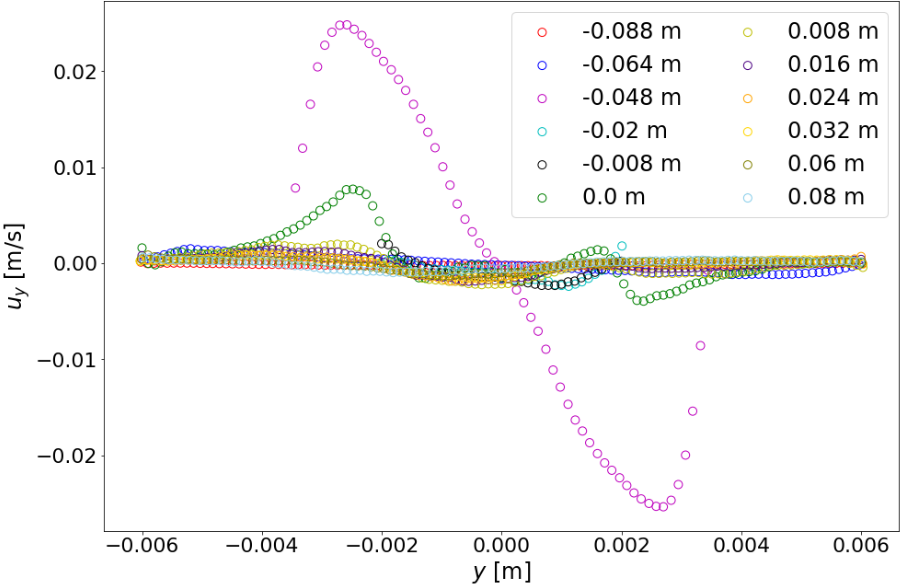}\label{fig:exp_ur}
		}\\
		\subfloat[][]{
			\includegraphics[width=0.31\textwidth]{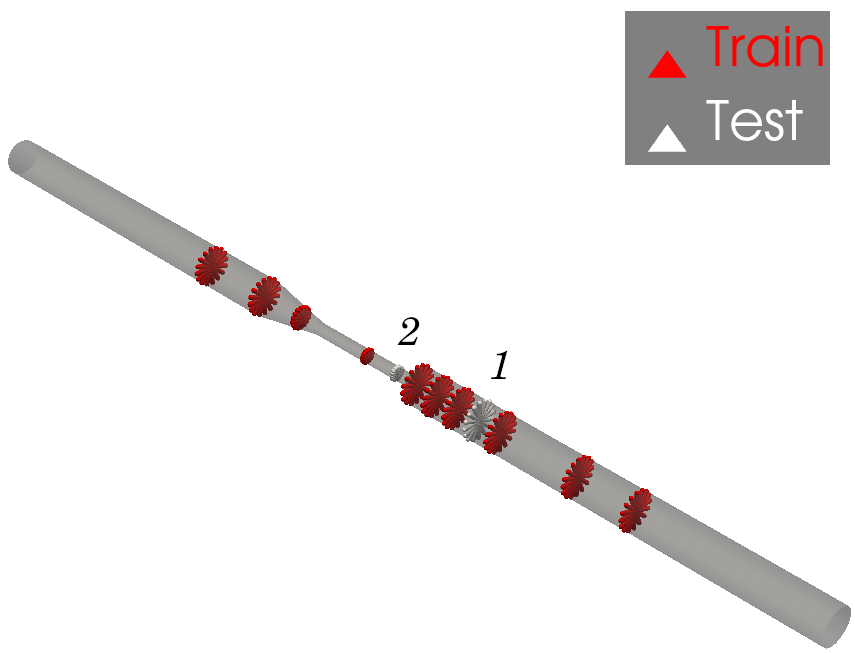}\label{fig:exp_sel}
		}
		\qquad\quad
		\subfloat[][]{
			\includegraphics[width=0.28\textwidth]{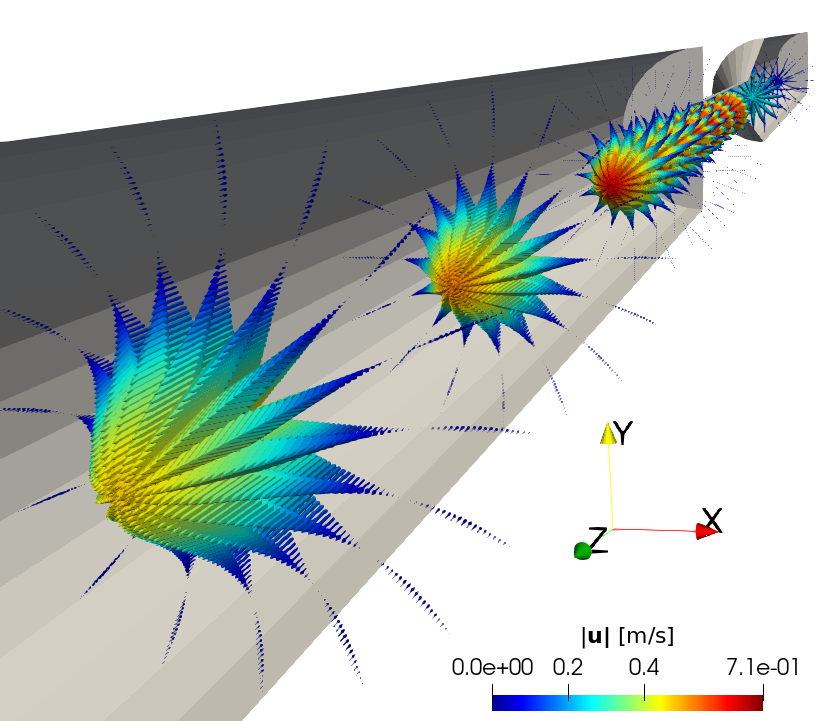}\label{fig:exp_vec}
		}
		\caption[]{PIV measurements from the \texttt{999} dataset~\cite{hariharan2011multilaboratory}. \textbf{(a)} Axial velocity $(\bm u)_z$ profile along the $y$-axis at different $z$ values. \textbf{(b)} Radial velocity $(\bm u)_x$ and $(\bm u)_y$ along the $y$-axis at different $z$ values. \textbf{(c)} Training and testing data sections partition. Test sections are: \num{1} at $z_1=-0.008\ \si{\meter}$ and \num{2} at $z_2=0.024\ \si{\meter}$. \textbf{(d)} Vector plot of the rotated data.}
		\label{fig:exp_dat}
	\end{figure}

	We train and test the ANN using solely the experimental data obtained from PIV measurements~\cite{hariharan2011multilaboratory} for the FDA nozzle benchmark dataset~\cite{NCIHub43}. When using the in silico dataset, the fluid properties were available in all points of the computational domain. Differently, in this case, the velocity is available only on certain cross-sections. First, in Section~\ref{sub:valiSetFDA}, we describe the setup. Subsequently, in Test 6, we test the FDA-EXP networks by progressively removing data sections on which the PINNs are trained.
	\paragraph{Setup}
	\label{sub:valiSetFDA}
	We describe the problem setup using the PIV experimental data. The data consist of axial and radial velocity samples on 12 different cross-sections, depicted in Figure~\ref{fig:exp_sel}. Five different data acquisitions are available for the FDA nozzle benchmark. We employ the dataset labeled \texttt{999} for the training, as it provides the most extensive results~\cite{hariharan2011multilaboratory}. The axial and radial velocity samples are represented in Figures~\ref{fig:exp_uz} and \ref{fig:exp_ur}.
    \par
    The training dataset is obtained starting from the $8992$ data points by removing select cross-sections and randomly sampling the data points over each remaining cross-section. We consider three different setups whose sampling rates and removed cross-sections are reported in Table~\ref{tab:setexp}.  
    \begin{table}[t]
		\centering
		\begin{tabular}{ p{8em} p{5.5em} p{5.5em} p{5.5em} }
			\hline
			& FDA-EXP-1 & FDA-EXP-2 & FDA-EXP-3 \\
			\hline
			Removed section & None & 1 & 1 and 2  \\
			Section sampling & $15\%$ & $20\%$ & $20\%$ \\
			\hline	
		\end{tabular}
		\caption{Setup for the different PINNs in Section~\ref{sec:valiFDA}. For section numbering we refer to Figure~\ref{fig:exp_sel}.}
		\label{tab:setexp}
	\end{table}
    To distribute the samples circumferentially over each section, we rotate the points and the velocity samples on the $y$-axis with respect to the $z$-axis by $22.5^\circ$ increments, resulting in angles ranging from $0^\circ$ to $157.5^\circ$. The resulting velocity vector data alignment is shown in Figure~\ref{fig:exp_vec}.
    \par
    The ANN architecture, minimization algorithm and loss function parameters are the same as for FDA-CFD-7 in Section~\ref{sec:veriSetFDA}, with the parameters reported in Table~\ref{tab:wstd}, the only difference being we set $\lambda_{\hat{\bm{u}},\mathrm{MSLE}} = \lambda_{{\bm{u}},\mathrm{MSLE}}=0$. We disregard the MSLE in the training loss function since we observed no major differences by incorporating it along with the unit velocity vector loss function $\mathcal{L}_{\hat{\bm u}}$.
    
	\paragraph{Test 6}
	\label{sub:valiResFDA}
	\begin{figure}[t]
		\centering
		\subfloat[][]{
			\includegraphics[width=0.31\textwidth]{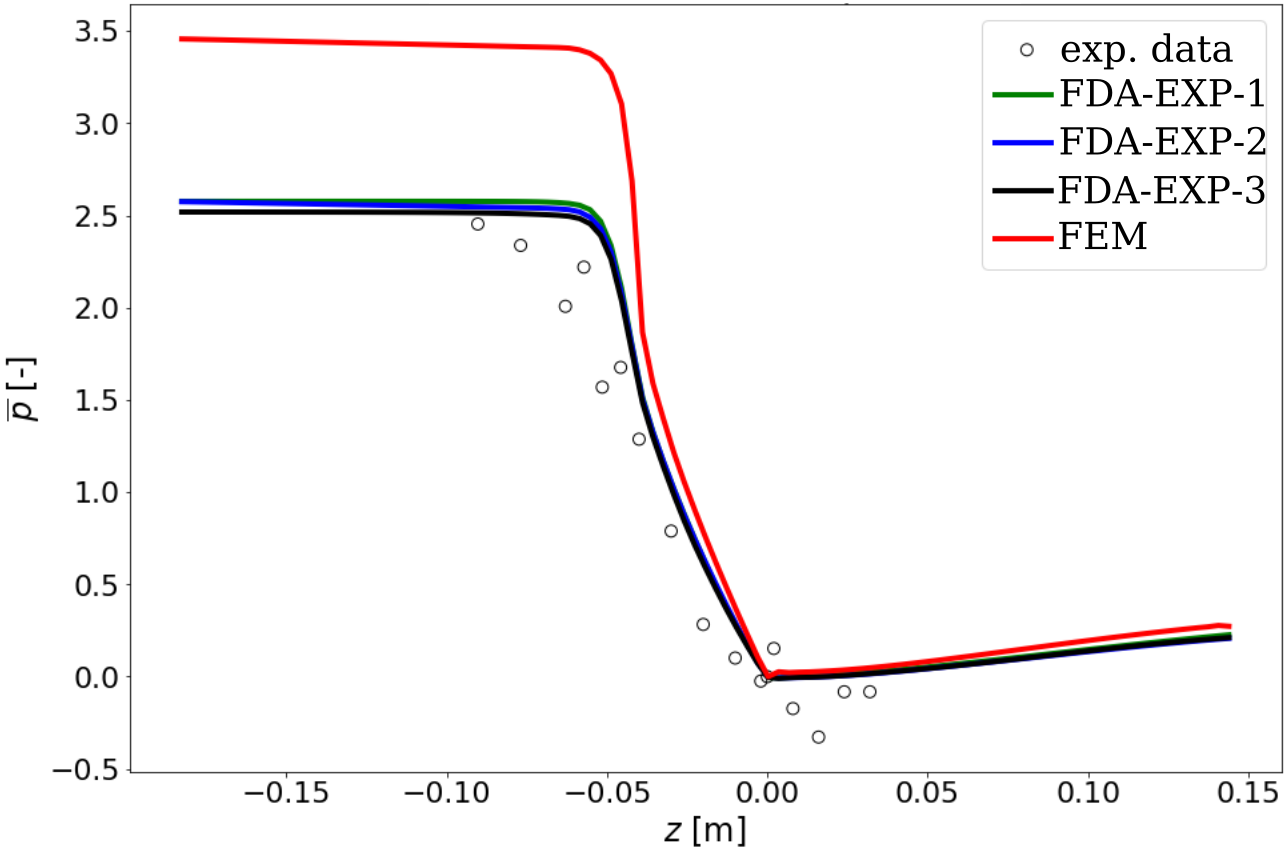}\label{fig:exp_pax}
		}
        \subfloat[][]{
			\includegraphics[width=0.31\textwidth]{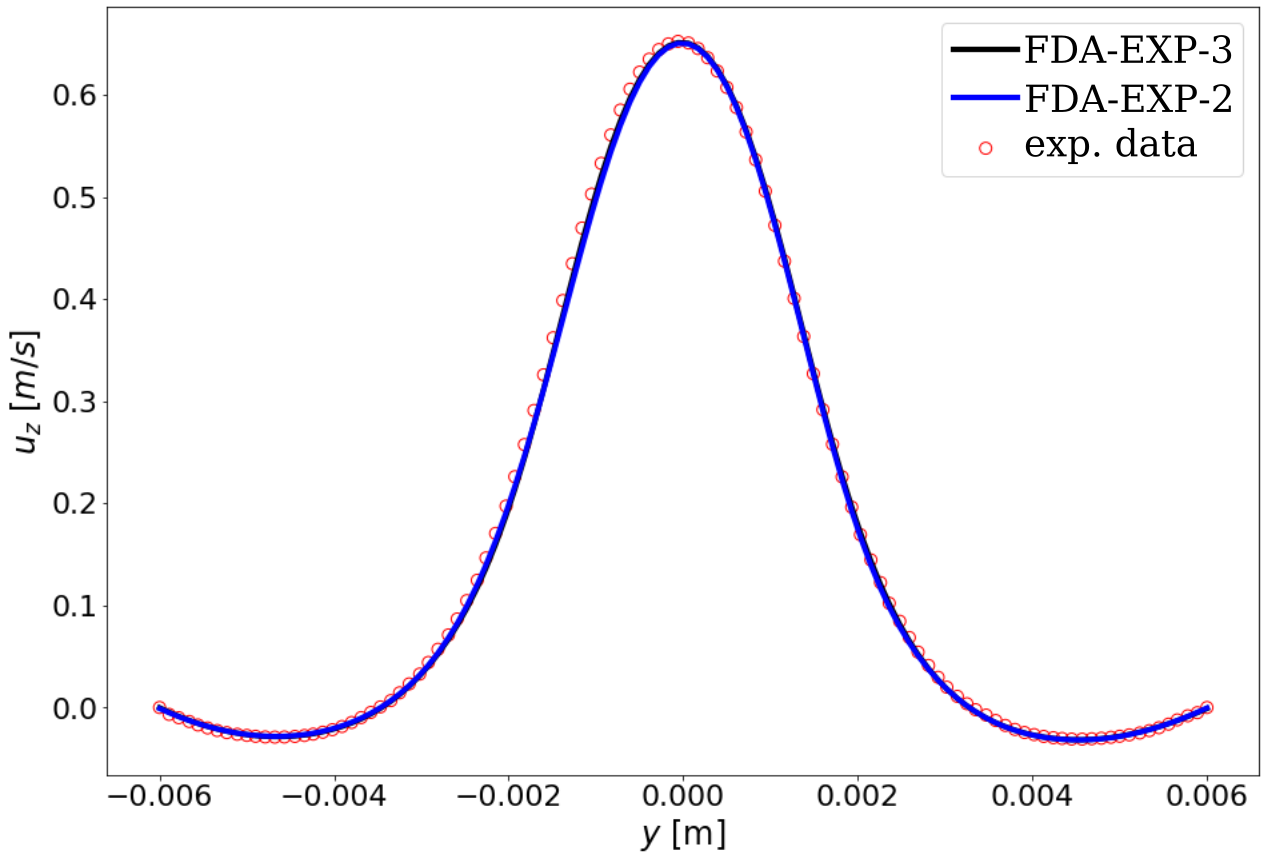}\label{fig:exp_024}
		}
		\subfloat[][]{
			\includegraphics[width=0.31\textwidth]{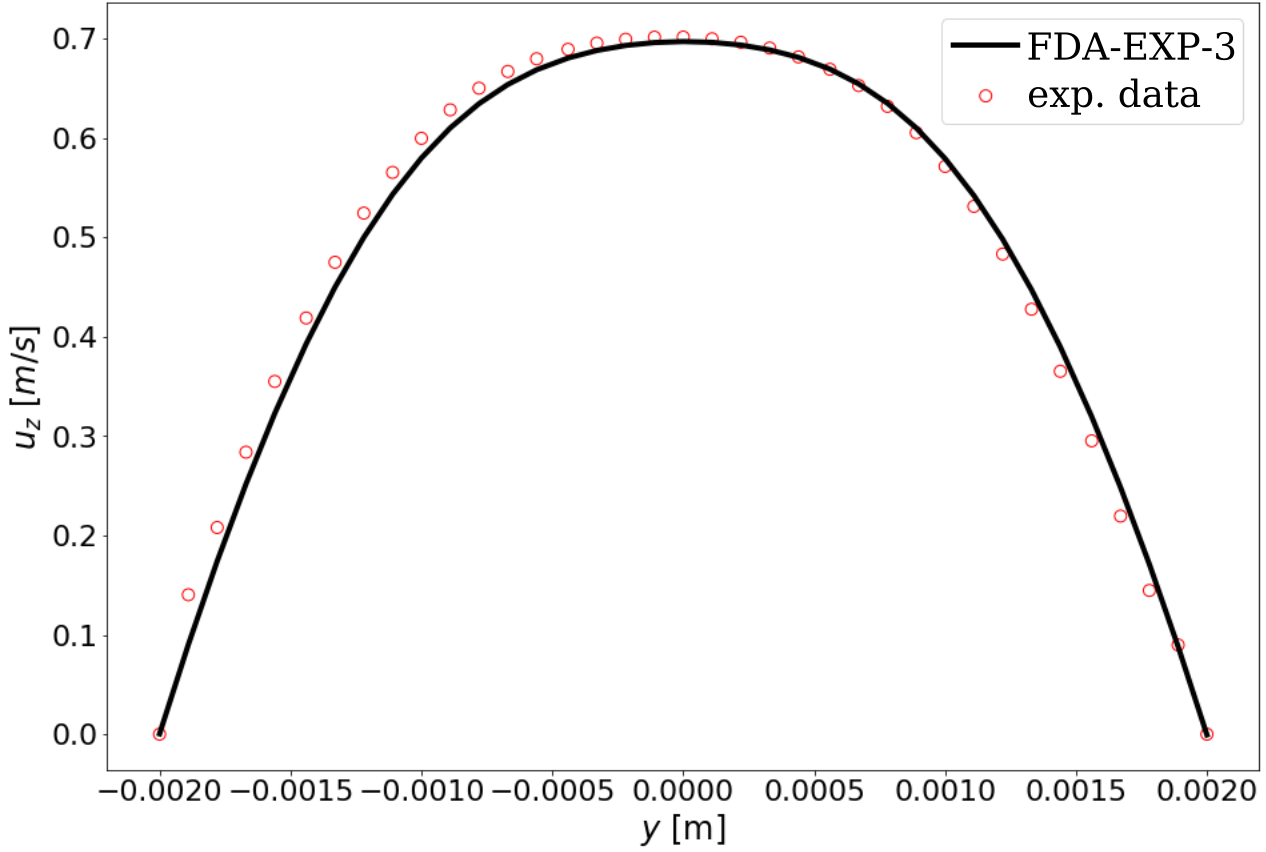}\label{fig:exp_m008}
		}
		\caption[]{Test 6: Comparison of the PINN results with the \texttt{999} experimental data~\cite{hariharan2011multilaboratory}.\textbf{(a)} Normalized center-line pressure $\overline{p}$. \textbf{(b)}  Axial velocity profile at $z = 0.024\ \si{\meter}$. \textbf{(c)} Axial velocity profile at $z = - 0.008\ \si{\meter}$.}
		\label{fig:exp_cent}
	\end{figure}
	
	\begin{figure}[t]
		\centering
		\subfloat[][]{
			\includegraphics[width=0.31\textwidth]{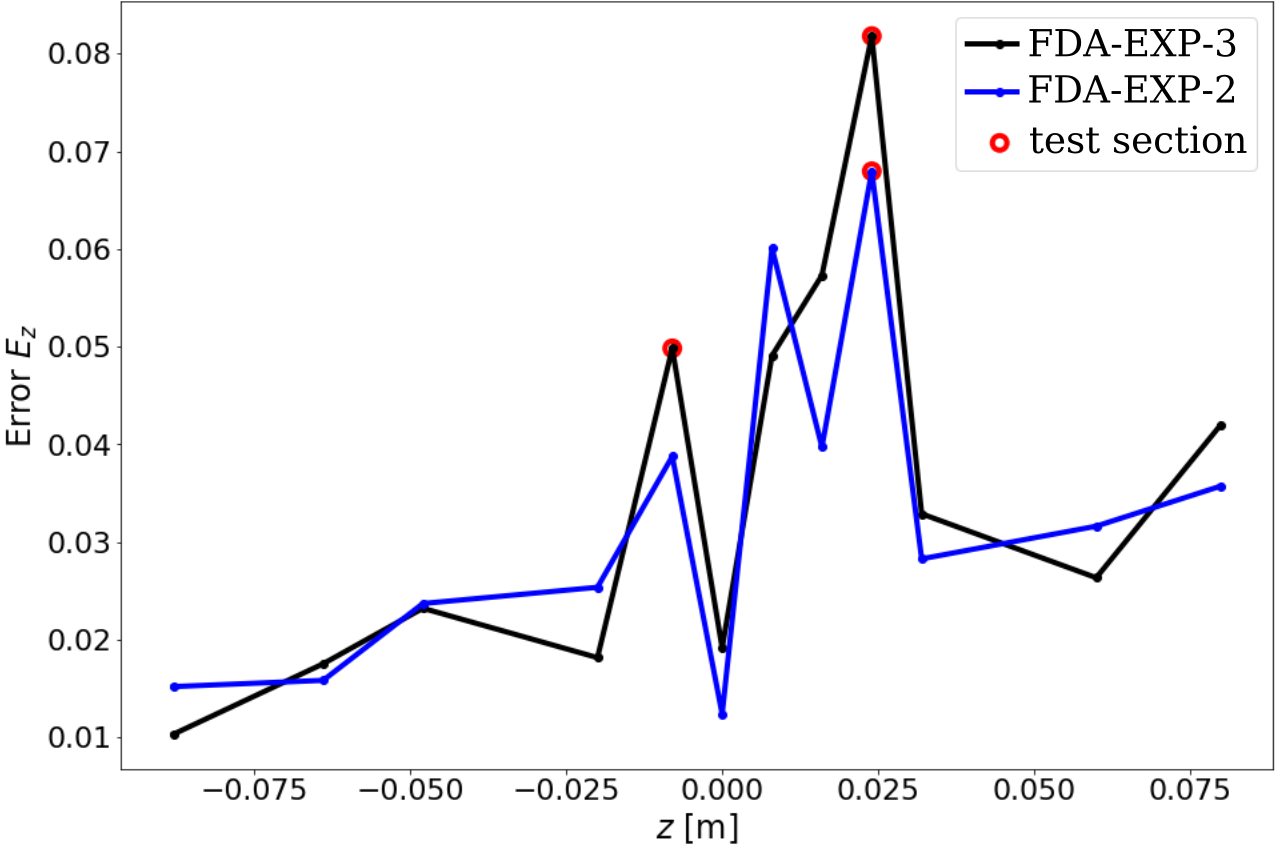}\label{fig:err_exp1}
		}
		\subfloat[][]{
			\includegraphics[width=0.31\textwidth]{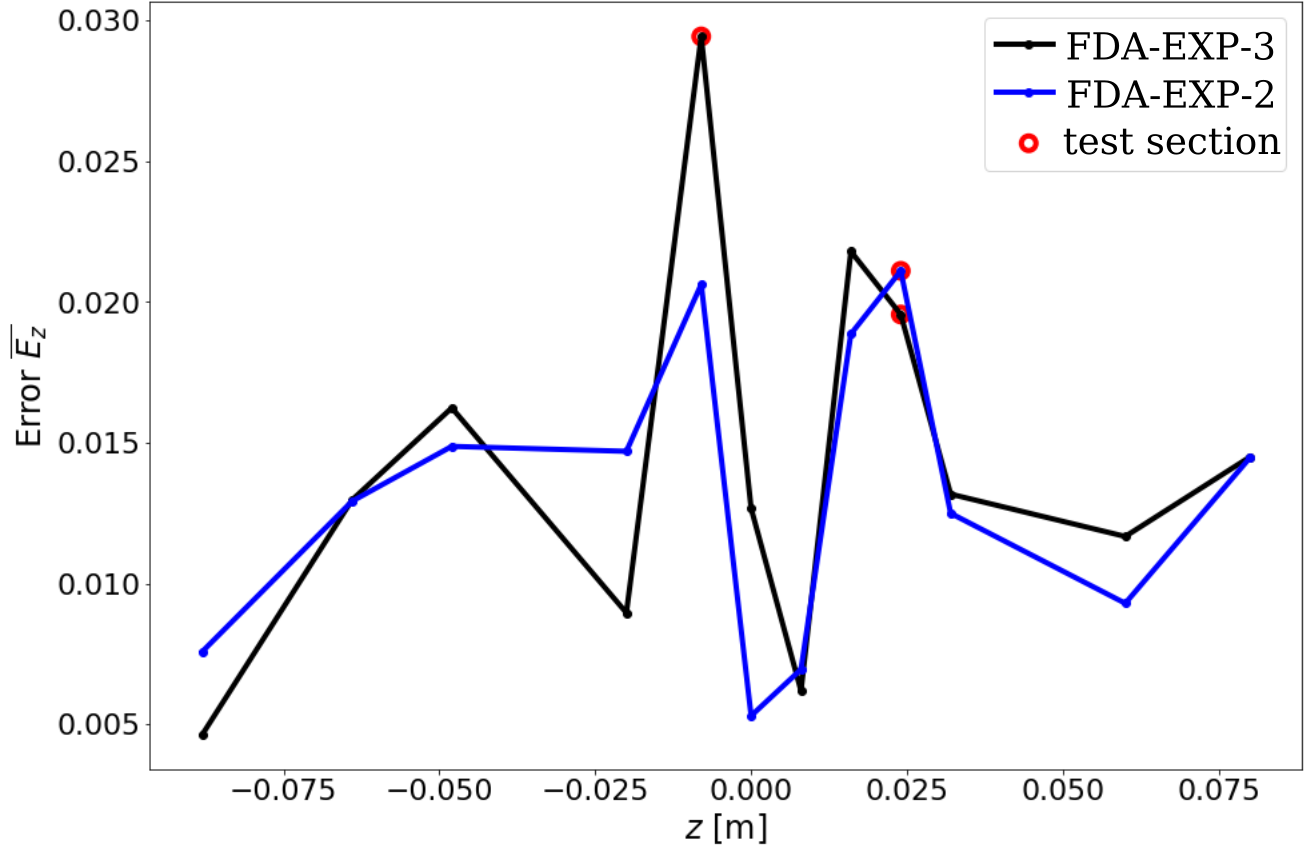}\label{fig:err_exp2}
		}
		\subfloat[][]{
			\includegraphics[width=0.31\textwidth]{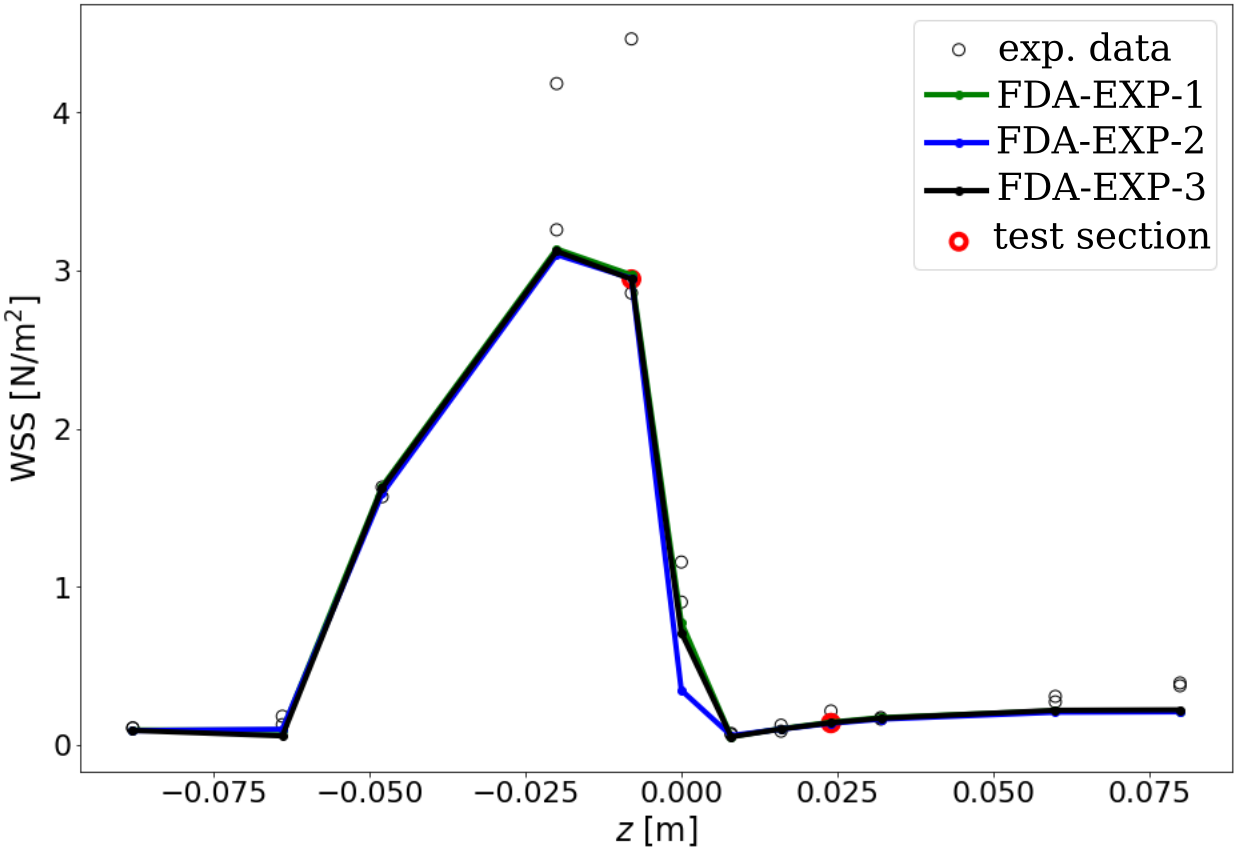}\label{fig:WSSexp}
		}
		\caption[]{Test 6: Comparison of errors computed from the PINNs with the ones computed using the experimental data on the sections. Red circles represent values at test sections. \textbf{(a)} Cross-section average of the relative error $E_z$. \textbf{(b)} Error of the average cross-section velocity $\overline{E_z}$. \textbf{(c)} Wall-averaged WSS magnitude comparison.}
		\label{fig:exp_err}
	\end{figure}

	\begin{figure}[t]
		\centering
		\subfloat[][]{
			\includegraphics[width=0.31\textwidth]{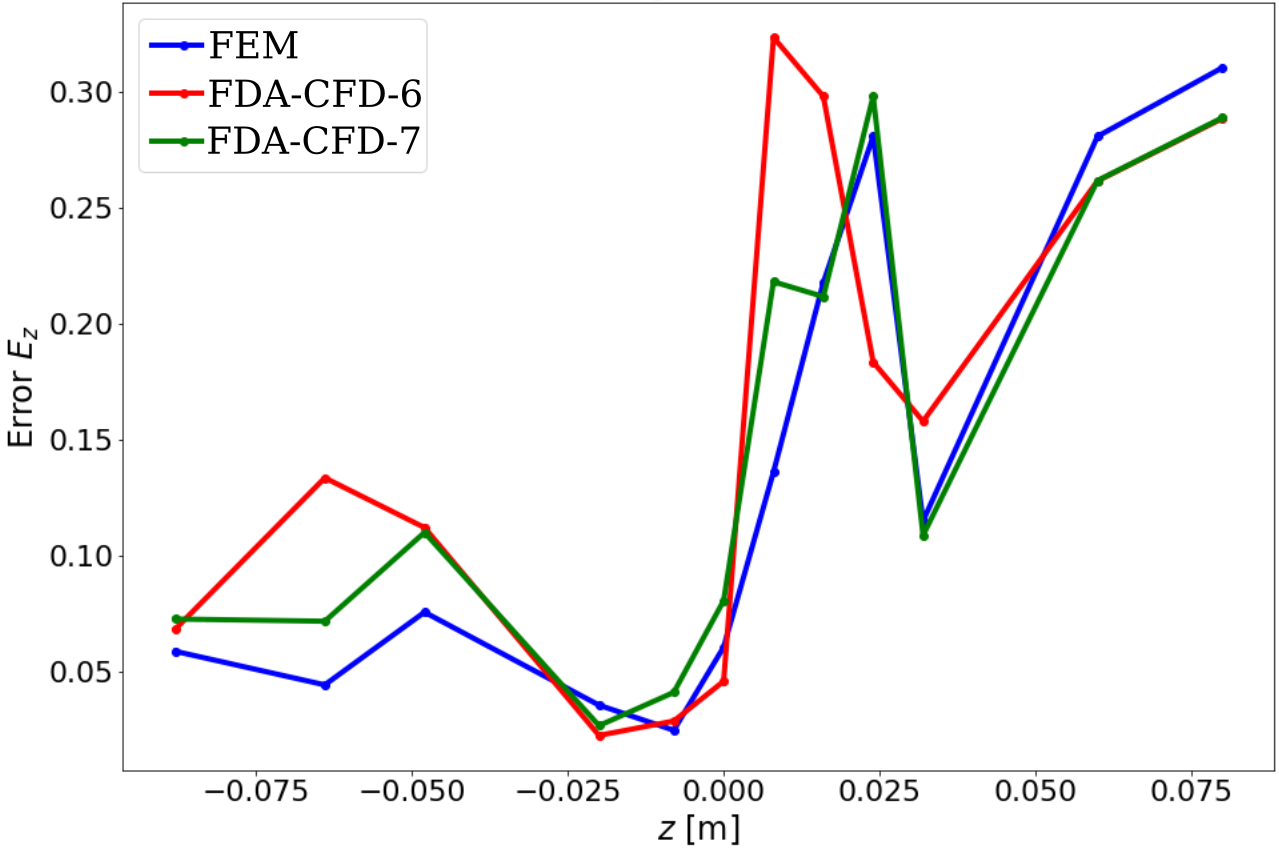}\label{fig:ezfin}
		}\quad
		\subfloat[][]{
			\includegraphics[width=0.31\textwidth]{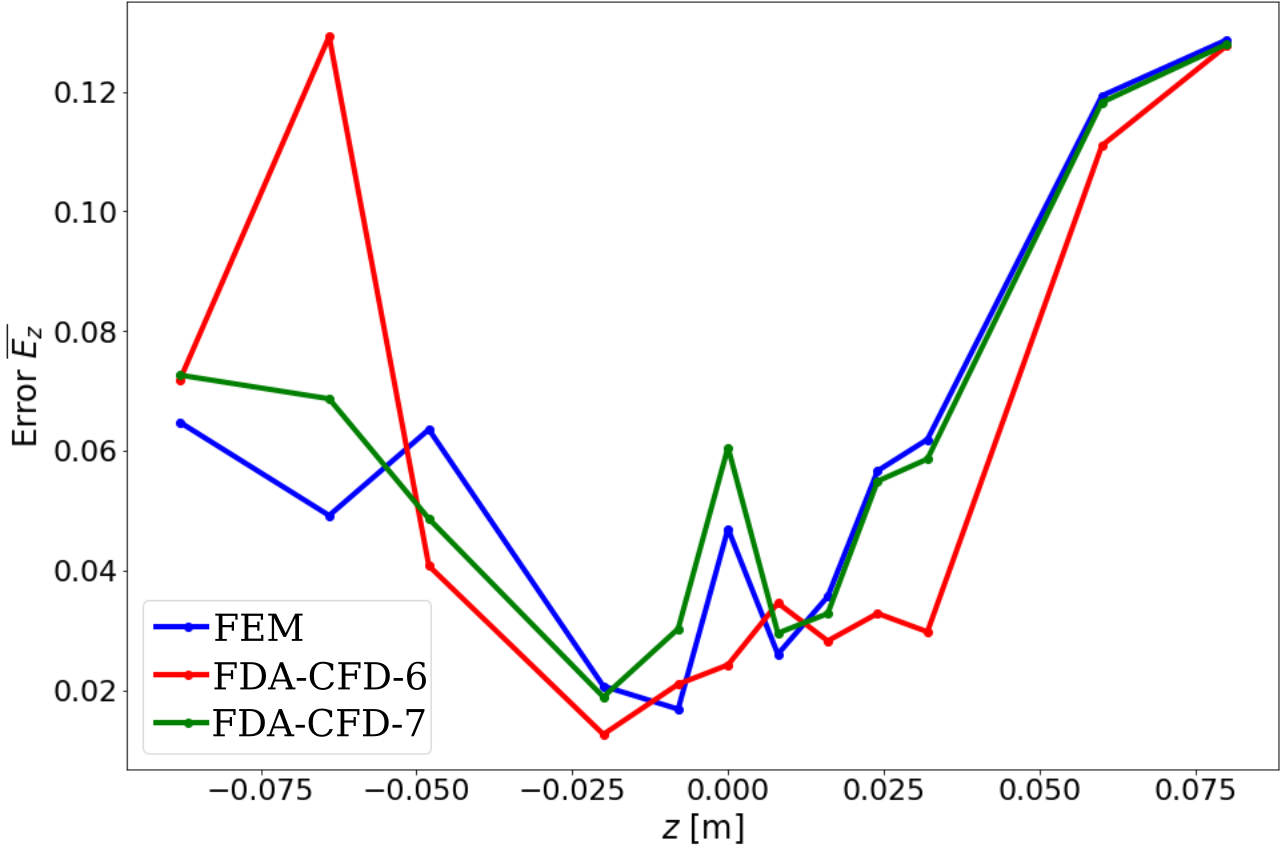}\label{fig:eztfin}
		}
		\caption[]{Test 5: Comparison of errors with respect to the experimental data computed from the PINNs with the ones computed using FEM data. \textbf{(a)} Cross-section average of the relative error ${E_z}$. \textbf{(b)} Error of the average cross-section velocity $\overline{E_z}$.}
		\label{fig:errfin}
	\end{figure}
    
    We compare the three different PINNs, trained with different sampling rates and by removing different cross-sections, as reported in Table~\ref{tab:setexp}. In Figure~\ref{fig:exp_pax} we report the reconstructed normalized pressure along the duct axis, where we see that all PINN setups performed far better at replicating experimental data than the FEM simulation.
    \par
    By comparing the reconstructed and experimental axial velocity for FDA-EXP-2 at the removed section at $z = 0.024~\si{\meter}$, in Figure~\ref{fig:exp_024}, we see that both PINNs compensate for the slight asymmetry in the experimental data. This also holds for FDA-EXP-3, where we have removed two sections from the training dataset, the first one at $z = 0.024\ \si{\meter}$ as in FDA-EXP-2, and the second one at $z = -0.008\ \si{\meter}$. While the first section is in a region where there are several other training sections, the second section is in a region where data are scarce. Moreover, it is located near the discontinuity due to the diffuser. Even so, the obtained velocity profile for the second section is very close to the experimental data, as evidenced by Figure~\ref{fig:exp_m008}. 
    \par
    In order to quantify the accuracy of the PINN approximations an to examine the effects of adding the experimental data in the loss function, we define errors between the axial velocity over each cross-section from the FEM simulation $(\bm u)_z$ and an average among the axial velocities of the five available experimental datasets $(\bm u_{
    \mathrm{exp}})_z$.
    \par
    We define the errors $E_z$ and $\overline{E_z}$ for a generic cross-section with axial coordinate $z$ as:  
	\begin{equation}
		E_z = \frac{1}{N_{z,\mathrm{sec}}}\sum_{i = 1}^{N_{z,\mathrm{sec}}} \left|\frac{(\bm u_{\mathrm{exp},i})_z- (\bm u_{i})_z }{(\bm u_{\mathrm{exp},i})_z}\right|, \qquad
		\overline{E_z} = \frac{1}{N_{z,\mathrm{sec}}}\sum_{i = 1}^{N_{z,\mathrm{sec}}} \left|\frac{(\bm u_{\mathrm{exp},i})_z-(\bm u_{i})_z }{(\bar{\bm u}_\mathrm{exp})_z}\right|,
		\label{eq:err}
	\end{equation}
    where ${N_{z,\mathrm{sec}}}$ is the total number of points on the given cross-section and $\bar{\bm u}_{\mathrm{exp}}$ is the average of the derived experimental values ${\bm u}_{\mathrm{exp},i}$ over the cross-section,
    \begin{equation*}
        \bar{\bm u}_{\mathrm{exp}} = \frac{1}{N_{z,\mathrm{sec}}}\sum_{i = 1}^{N_{z,\mathrm{sec}}}\bm u_{\mathrm{exp},i}.
    \end{equation*}
    In other words, in Eq.~\eqref{eq:err}, $E_z$ is the average of the axial velocity relative errors on the cross-section, while in $\overline{E_z}$ the relative error is computed with respect to the average experimental axial velocity on the cross-section. In particular, by considering $\overline{E_z}$ we avoid divisions by small values of ${u}_{z,\mathrm{exp}}$ near the wall that might cause $E_z$ to assume very large values.
	\par
    We compute the errors with respect to the experimental dataset using Eq.~\eqref{eq:err}, and report them in Figures~\ref{fig:err_exp1} and~\ref{fig:err_exp2}. Both FDA-EXP-2 and FDA-EXP-3 display a similar error trend. In addition to the approximated velocities, we compute the wall-averaged WSS using Eq.~\eqref{eq:wss} and compare it to the experimental data in Figure~\ref{fig:WSSexp}. The computed wall-averaged WSS approximates the experimental data well in the low-stress region, even though the data are quite dispersed. However, in the convergent region, where the WSS data are highest, all ANNs underestimate the WSS.
    \par
    In order to examine the effect of adding the experimental dataset to the loss function, we compute the errors using Eq.~\eqref{eq:err} of the PINNs in Section~\ref{sec:veriFDA}, that are trained only using in silico data. We report the computed errors in Figures~\ref{fig:ezfin} and~\ref{fig:eztfin}. From Figures~\ref{fig:ezfin} and \ref{fig:eztfin} we can see that the FDA-CFD-6 setup leads to larger errors in the first part of the duct, whereas the errors obtained using the FDA-CFD-7 setup follow more closely the error of the FEM simulation. This is coherent with the fact that the training data are taken from the FEM simulation and that FDA-CFD-7 has a larger number of neurons and therefore parameters. Overall, we see that PINNs trained using only the FEM data exhibit overall larger errors than the PINNs trained using the experimental data.    
	\subsection{Application to blood flows in aneurysm}
	\label{sec:aneurysm}
	We apply the PINN method to describe the flow in a realistic aneurysm geometry. First, in Section~\ref{sec:aneVeri}, we verify the PINN method using FEM simulation data. Then, in Section~\ref{sec:aneVali}, we show results for the validation of our method using 4D flow MRI data~\cite{pravdivtseva20213d}.
	\subsubsection{Numerical verification against in silico data}
	\label{sec:aneVeri}
	We perform PINN validation with in silico data obtained from FEM simulations. In Section~\ref{sec:veriSet} we describe the FEM simulation setup and results used for the training. Next, in Section \ref{sec:veriRes} we discuss the results of the PINN training.
	\paragraph{Setup}
	\label{sec:veriSet}
	
	\begin{figure}[t]
		\centering
		\hspace{-1cm}\subfloat[][]{
			\includegraphics[width=0.68\textwidth]{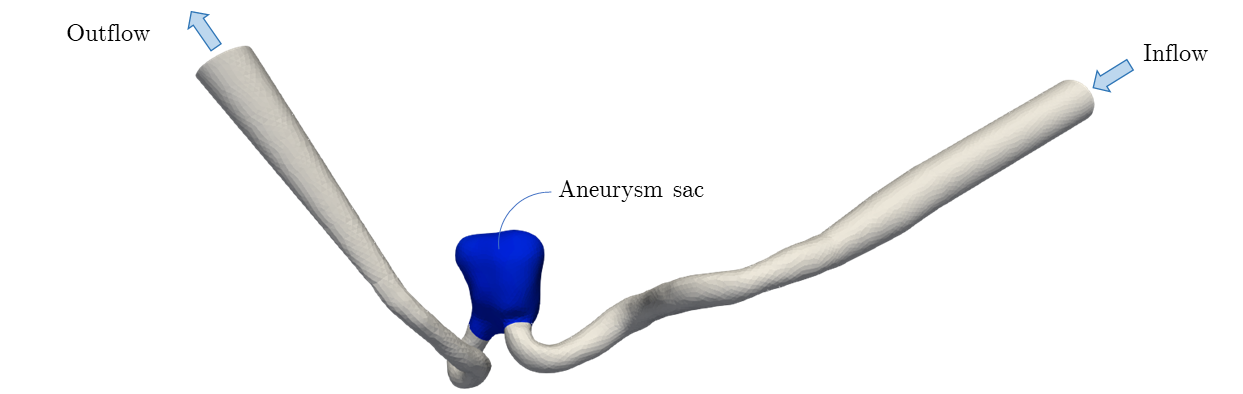}\label{fig:modane}
		}
		\subfloat[][]{
			\includegraphics[width=0.28\textwidth]{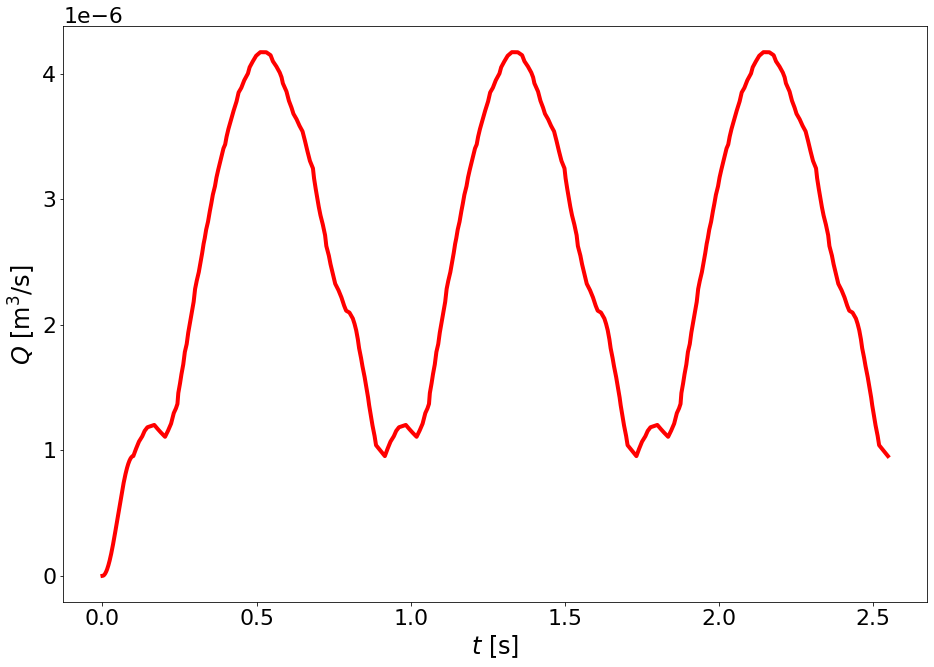}\label{fig:inl_ane}
		}
		\caption[]{Aneurysm geometry and simulation features. \textbf{(a)} Aneurysm model, provided by~\cite{pravdivtseva20213d}. The region where 4D flow MRI data are available is in blue. \textbf{(b)} Imposed inlet flow rate over time.}
		\label{fig:setber}
	\end{figure}
    We consider the aneurysm model represented in Figure~\ref{fig:modane}, provided in~\cite{pravdivtseva2022influence}. In their study, the geometry is obtained from an in vitro patient-derived 3D-printed aneurysm model. It is segmented from an intraophthalmic, extradural internal carotid artery aneurysm of a patient (52 years, female)~\cite{pravdivtseva20213d}.
    \par
    We carry out the FEM simulations using \texttt{life$^\texttt{X}$}~\cite{africa2022lifex}. We set the flow properties, reported in Table~\ref{tab:fluid_flow_properties}, second column, in accordance with the experimental data~\cite{pravdivtseva20213d}. The Reynolds number is $\mathbb{R}\mathrm{e} \approx 1500$ at the inlet of the aneurysm sac. We impose a parabolic velocity profile at the inlet with a volumetric flow rate as shown in Figure~\ref{fig:inl_ane}. We apply a homogeneous Neumann boundary condition at the outlet, and a Dirichlet no-slip boundary condition at the lateral wall. As in Section~\ref{sec:veriSetFDA}, we apply the Backward Euler scheme with a semi-implicit treatment of the nonlinear term along with the VMS-LES method. We consider simulations on two different mesh refinement levels and simulate a total of three heartbeats. We report the simulation parameters in Table~\ref{table:FEM}. Simulations were ran on 32 cores on the \texttt{gigatlong} cluster of the Department of Mathematics, Politecnico di Milano~\cite{polimiHardware}. 
    \par
\begin{figure}[t]
		\centering
		\subfloat[][]{
			\includegraphics[width=0.24\textwidth]{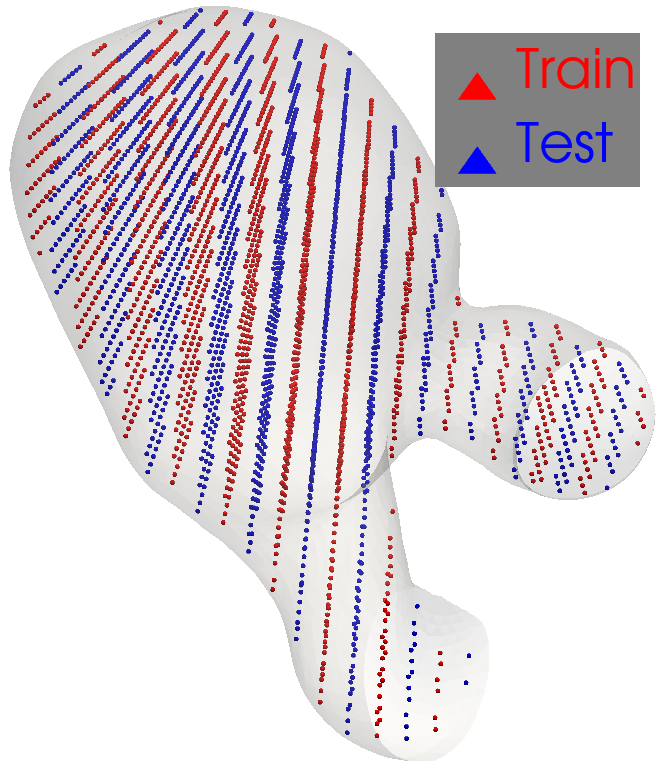}\label{fig:FEManetrain}
		}
		\qquad \qquad
		\subfloat[][]{
			\includegraphics[width=0.24\textwidth]{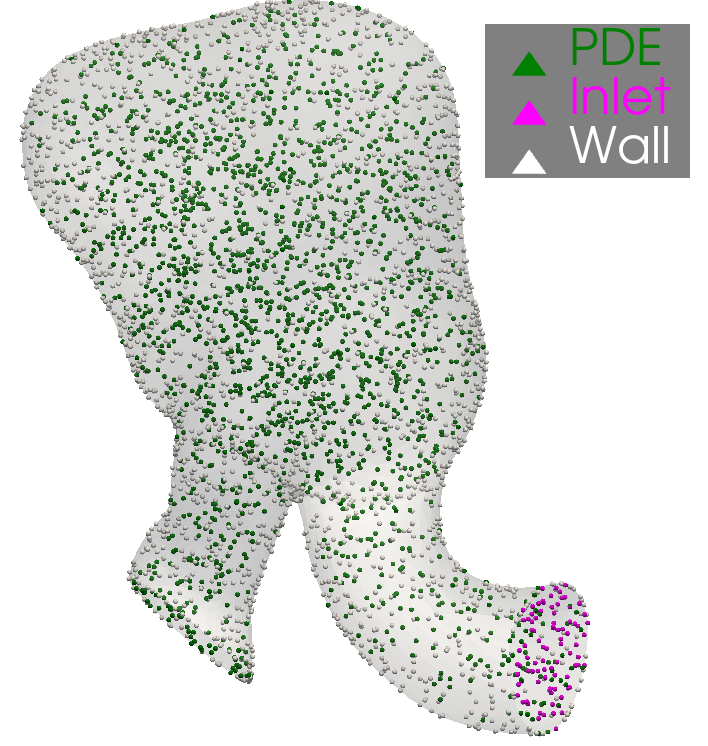}\label{fig:FEManePDE}
		}
		\caption[]{Collocation points. \textbf{(a)} Training and testing dataset partition. \textbf{(b)} Collocation points for the PDE regularization.}
		\label{fig:FEManecol}
	\end{figure}

    For training/test data generation we consider a part of the domain in the proximity of the aneurysm sac. We sample the velocity field over a uniform grid with a resolution of $0.7 \times 0.7 \times 0.7 \ \si{\milli\meter\cubed}$ along the Cartesian directions, which yields a velocity distribution coherent with the 4D flow MRI data~\cite{pravdivtseva2022influence}. We obtain the training dataset by removing half of the cross-sections along the $z$ direction, according to Figure~\ref{fig:FEManetrain}.
    \par
    Regarding the loss function associated to the PDE, we use the same definition of residual as for the previously analyzed steady-state case in Section~\ref{sec:FDAbench}, by fixing the time in Eq.~\eqref{eq:NS_res} to $t=2.024~\si{\second}$. Hence, we avoid the computational cost of adding the physical time as an input variable to the PINNs. Some information on the time-dependent terms is nonetheless recovered through the data-driven loss function~$\mathcal{L}_\mathrm{train}$. We define the total loss function~$\mathcal{L}$ as:
    \begin{equation}
        \label{eq:loss_ane}
        \mathcal{L}(\Tilde{\bm{\theta}}) = \mathcal{L}_{\mathrm{PDE}}(\Tilde{\bm{\theta}}) + \mathcal{L}_{\mathrm{BC}}(\Tilde{\bm{\theta}}) +
		\mathcal{L}_{\bar p}(\Tilde{\bm{\theta}}) +
		\mathcal{L}_{\mathrm{train}}(\Tilde{\bm{\theta}}),
    \end{equation}
    where $\mathcal{L}_{\bar p}$ represents the zero mean constraint penalty on the pressure:
	\begin{equation}
		\label{eq:pmean}
        \mathcal{L}_{\bar{p}}(\Tilde {\bm{\theta}}) = \frac{\lambda_{\bar{p}}}{N_\mathrm{PDE}}\left|\sum_{i=1}^{N_\mathrm{PDE}} p_i(\Tilde {\bm{\theta}})\right|^2,
	\end{equation}
    where $\lambda_{\bar{p}}$ is a positive parameter and the reconstructed pressure $p$ is evaluated at the same $N_\mathrm{PDE}$ collocation points as $\mathcal{L}_\mathrm{PDE}$. The role of the zero mean pressure loss function $\mathcal{L}_{\bar{p}}$ in Eq.~\eqref{eq:pmean} is akin to the zero mean constraint used in FEM simulations to guarantee pressure uniqueness when no Neumann boundary condition is applied. We introduce \eqref{eq:pmean} since we observed that using it instead of the Neumann boundary condition residual~\eqref{eq:NS_N}, by setting $\lambda_{\mathrm{BC},3}=0$,
    we improve the accuracy of the pressure reconstruction.
    \par
	We use $N_\mathrm{train} = 1420$ points for training and $N_\mathrm{test} = 1436$ points for testing. We randomly seed $N_\mathrm{PDE} = 1959$ collocation points for the PDE loss function and $N_\mathrm{BC,2} = 2057$ collocation points at the wall, which we report in Figure~\ref{fig:FEManePDE}.
    \par
	For the ANNs we use three hidden layers of 32 neurons each, with batch normalization before each layer. We solve the minimization problem related to the training using \num{100} epochs of the ADAM optimizer, with a learning rate of $10^{-2}$, and \num{14000} epochs of L-BFGS-B. The total time to perform the numerical minimization\footnote{Computational resources used for training are CPU 11th Gen Intel(R) Core(TM) i7-11800H @ 2.30GHz, GPU NVIDIA GeForce RTX 3070} is approximately 3 hours.
    \par
\begin{figure}[t]
		\centering
		\subfloat[][]{
        \includegraphics[width=0.22\textwidth]{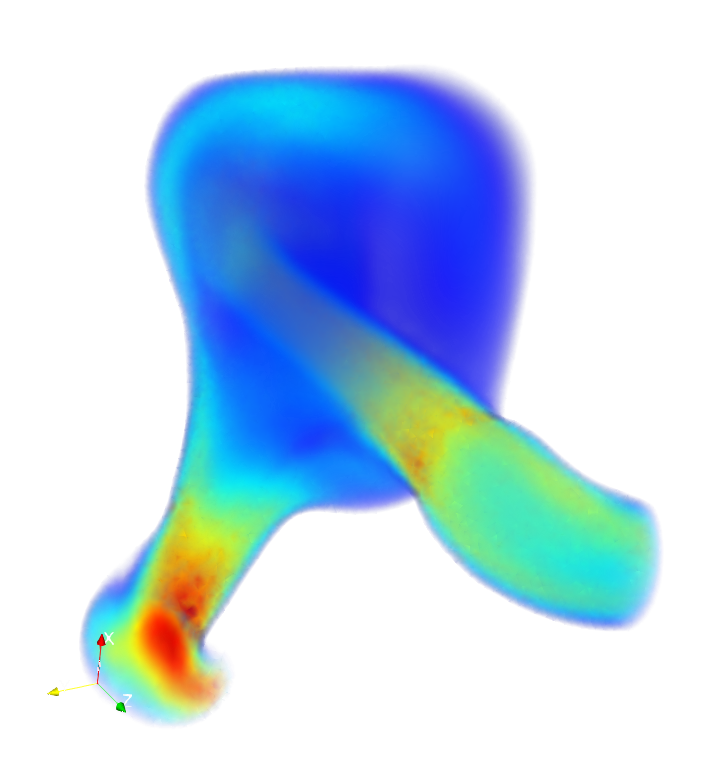}\label{fig:FEManev}
		}
		\subfloat[][]{
			\includegraphics[width=0.22\textwidth]{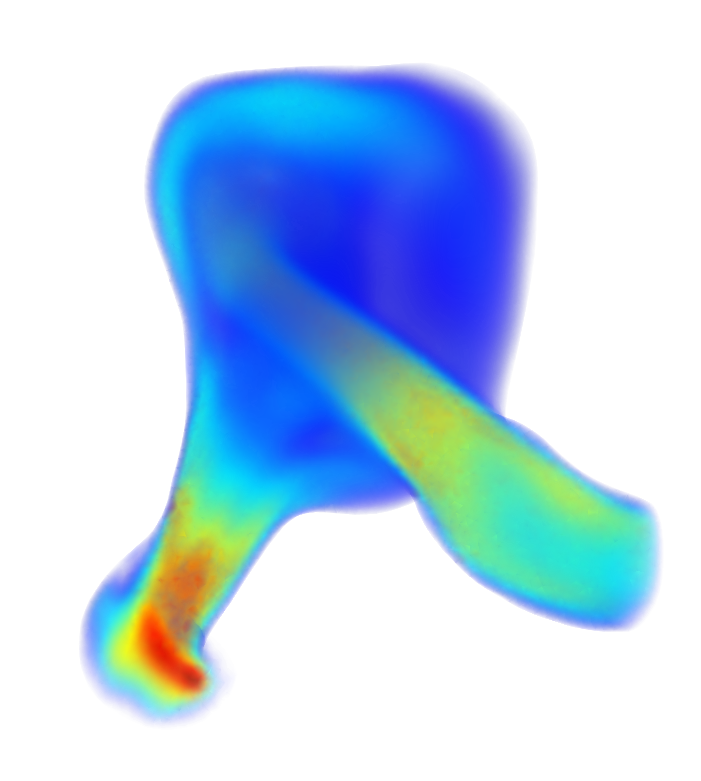}\label{fig:FEManeu1}
		}
		\subfloat[][]{
			\includegraphics[width=0.22\textwidth]{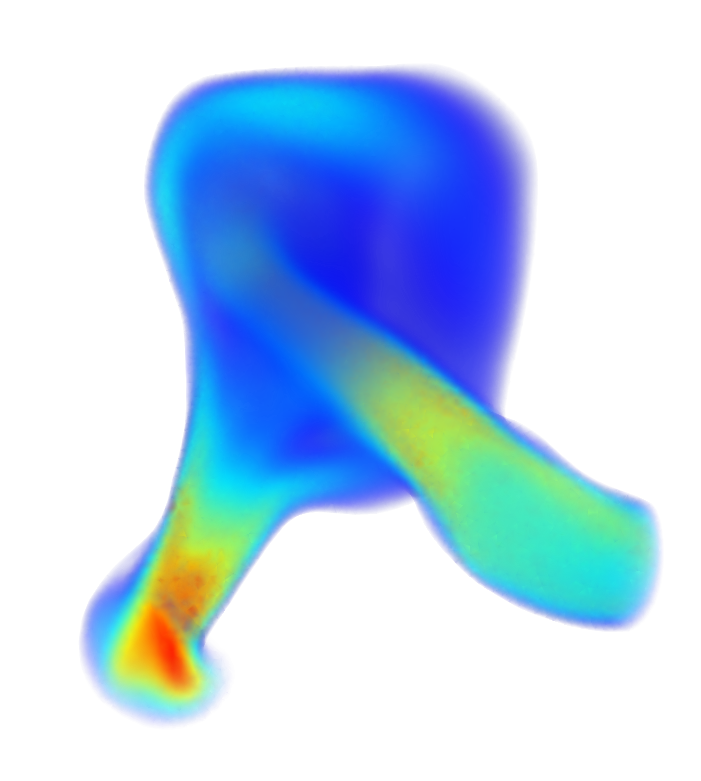}\label{fig:FEManeu2}
		}
		\subfloat[][]{
			\includegraphics[width=0.22\textwidth]{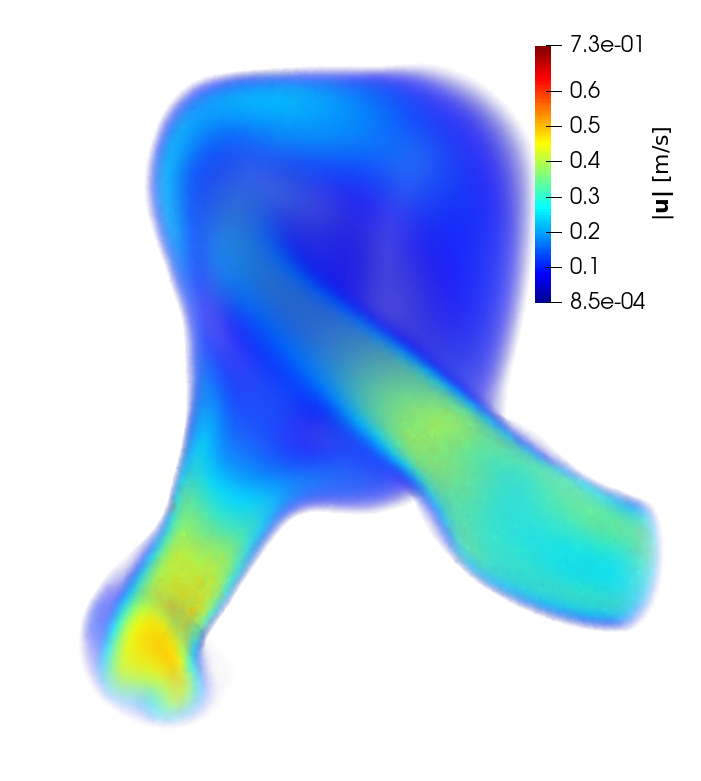}\label{fig:FEManeu3}
		}
		\caption[]{Computed magnitude of the velocity field. \textbf{(a)} FEM. \textbf{(b)} AN-CFD-1. \textbf{(c)} AN-CFD-2. \textbf{(d)} AN-CFD-3.}
		\label{fig:FEManevel}
	\end{figure}
\begin{figure}[t]
		\centering
		\subfloat[][]{
			\includegraphics[width=0.22\textwidth]{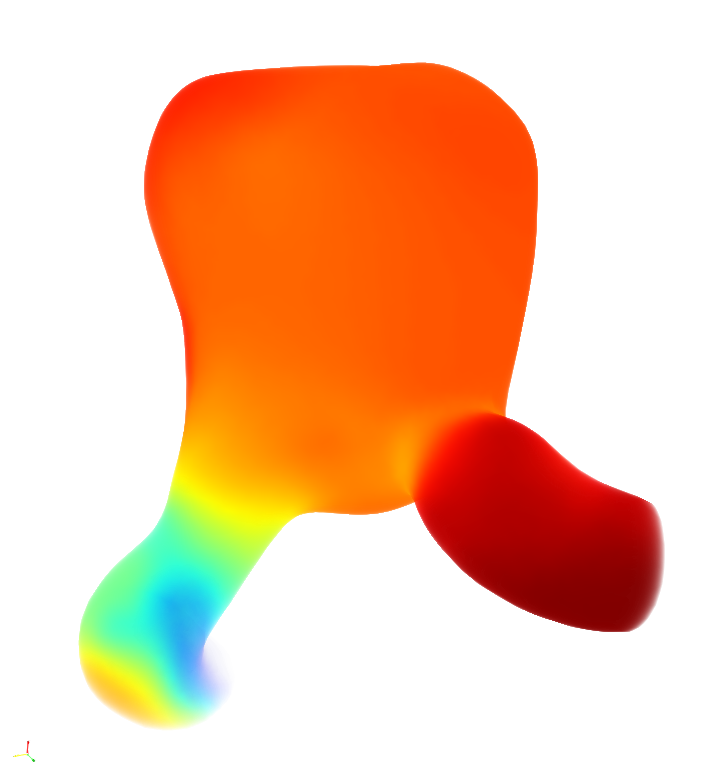}\label{fig:FEManep}
		}
		\subfloat[][]{
			\includegraphics[width=0.22\textwidth]{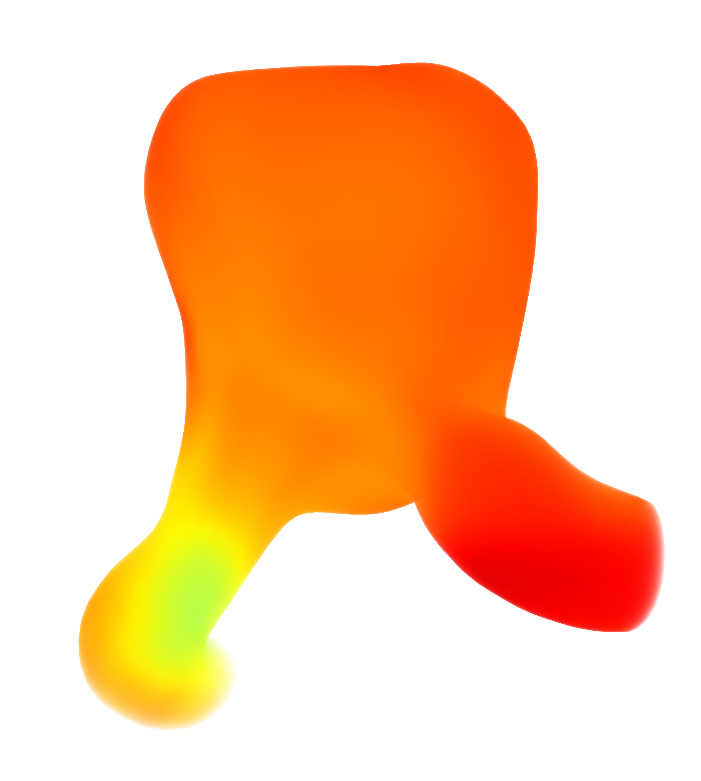}\label{fig:FEManep1}
		}
		\subfloat[][]{
			\includegraphics[width=0.22\textwidth]{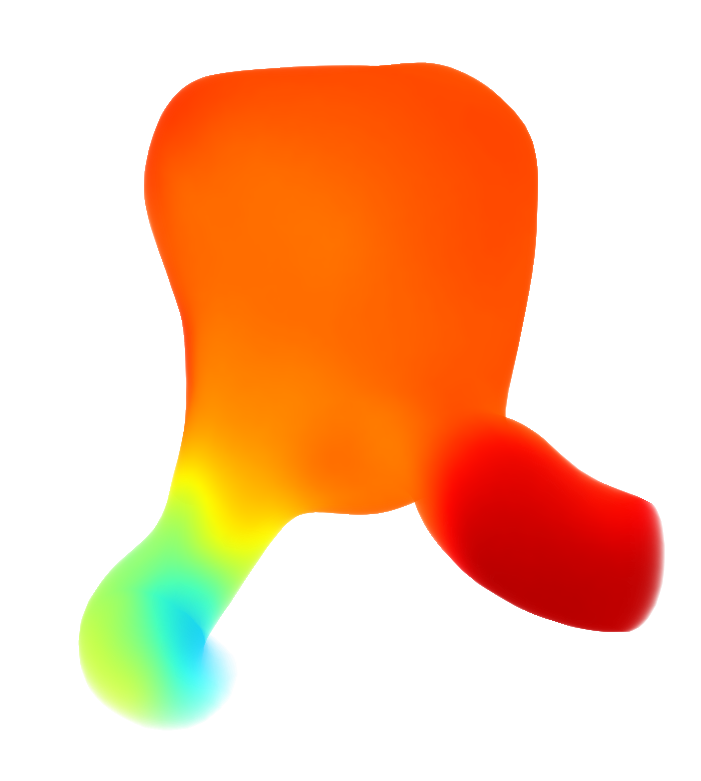}\label{fig:FEManep2}
		}
		\subfloat[][]{
			\includegraphics[width=0.22\textwidth]{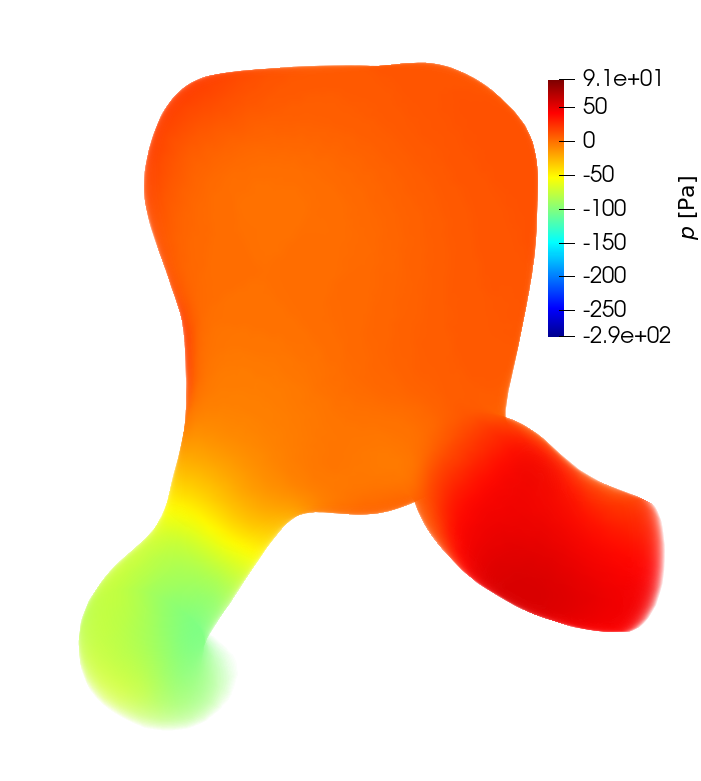}\label{fig:FEManep3}
		}
		\caption[]{Reconstructed pressure field. \textbf{(a)} FEM. \textbf{(b)} AN-CFD-1. \textbf{(c)} AN-CFD-2. \textbf{(d)} AN-CFD-3.}
		\label{fig:FEManepres}
	\end{figure}

    We train three different PINNs with decreasing weights on the training data loss terms with respect to the physics loss term. We summarize the weights for the three different configurations in Table~\ref{tab:wFEMane}.
    \begin{table}[t]
		\centering
		\begin{tabular}{p{4em} p{5em} p{5em} p{5em}}
			\hline
			& AN-CFD-1 & AN-CFD-2 & AN-CFD-3 \\
			\hline
			$\lambda_{u_1}$ & 170 & 17 & 1.7 \\
            $\lambda_{u_2}$ & 170 & 17 & 1.7 \\
            $\lambda_{u_3}$ & 170 & 17 & 1.7 \\
            $\lambda_{\mathrm{BC,1}}$ & 0.5 & 0.5 &0.5\\
			$\lambda_{\mathrm{BC,2}}$ & 21.6 & 21.6 & 21.6 \\
            $\lambda_{\mathrm{BC,3}}$ & 0 & 0 & 0 \\
			$\lambda_{\mathrm{PDE,1}}$ & 5.76 & 5.76 & 5.76 \\
			$\lambda_{\mathrm{PDE,2}}$ & 1 & 1 & 1 \\
			$\lambda_{\mathrm{PDE,3}}$ & 1 & 1 & 1 \\
			$\lambda_{\mathrm{PDE,4}}$ & 1.56 & 1.56 & 1.56 \\
			$\lambda_{\hat{\bm{u}}}$ & 17 & 1.7 & 0.17  \\
			$\lambda_{\bar p}$ & 0.5 & 0.5 & 0.5  \\
			\hline
		\end{tabular}
		\caption{Parameters for the loss functions used for PINNs training in Section~\ref{sec:aneVeri}.}
		\label{tab:wFEMane}
	\end{table}
	\paragraph{Results}
	\label{sec:veriRes}
	We report the results of the PINNs trained using the FEM simulation output. In Figure~\ref{fig:FEManevel}, we show a comparison between the velocity magnitude field obtained from the FEM simulation and the velocity magnitude obtained from the PINNs. We observe that by decreasing the weight of the data in the loss function the peak velocity magnitude decreases.
	Next, in Figure~\ref{fig:FEManepres}, we report the reconstructed pressure fields, where we zero average the  FEM simulation pressure over the domain. We observe that the main features of the pressure field are matched by the PINN output. However, AN-CFD-1 and AN-CFD-3 predict lower pressure gradients with respect to the FEM simulation in the two vessel branches. On the other hand, AN-CFD-2 matches the FEM simulation best, with some discrepancy close to the outlet.

	\begin{figure}[t]
		\centering
		
		\subfloat[][]{
			\includegraphics[width=0.31\textwidth]{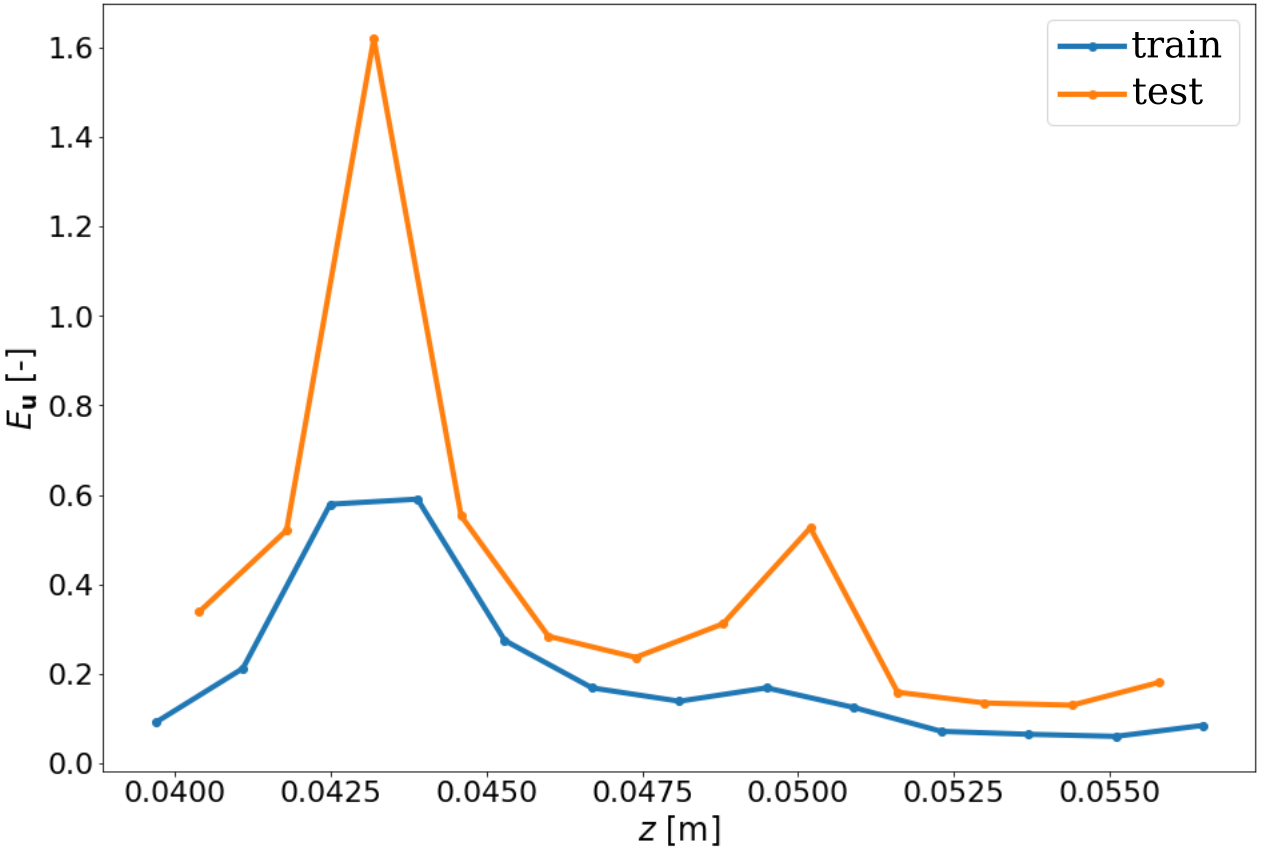}\label{fig:FEManeer1}
		}
		\subfloat[][]{
			\includegraphics[width=0.31\textwidth]{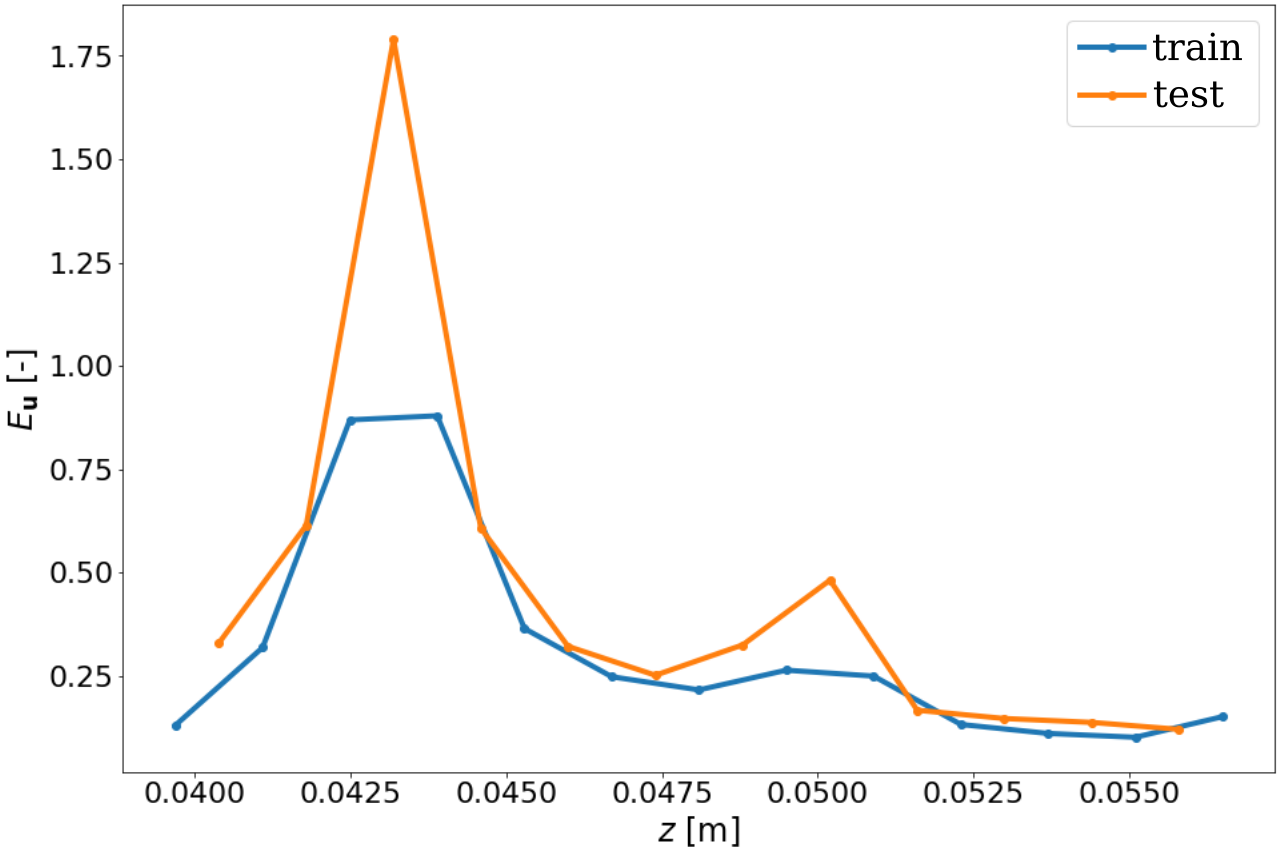}\label{fig:FEManeer2}
		}
		\subfloat[][]{
        \includegraphics[width=0.31\textwidth]{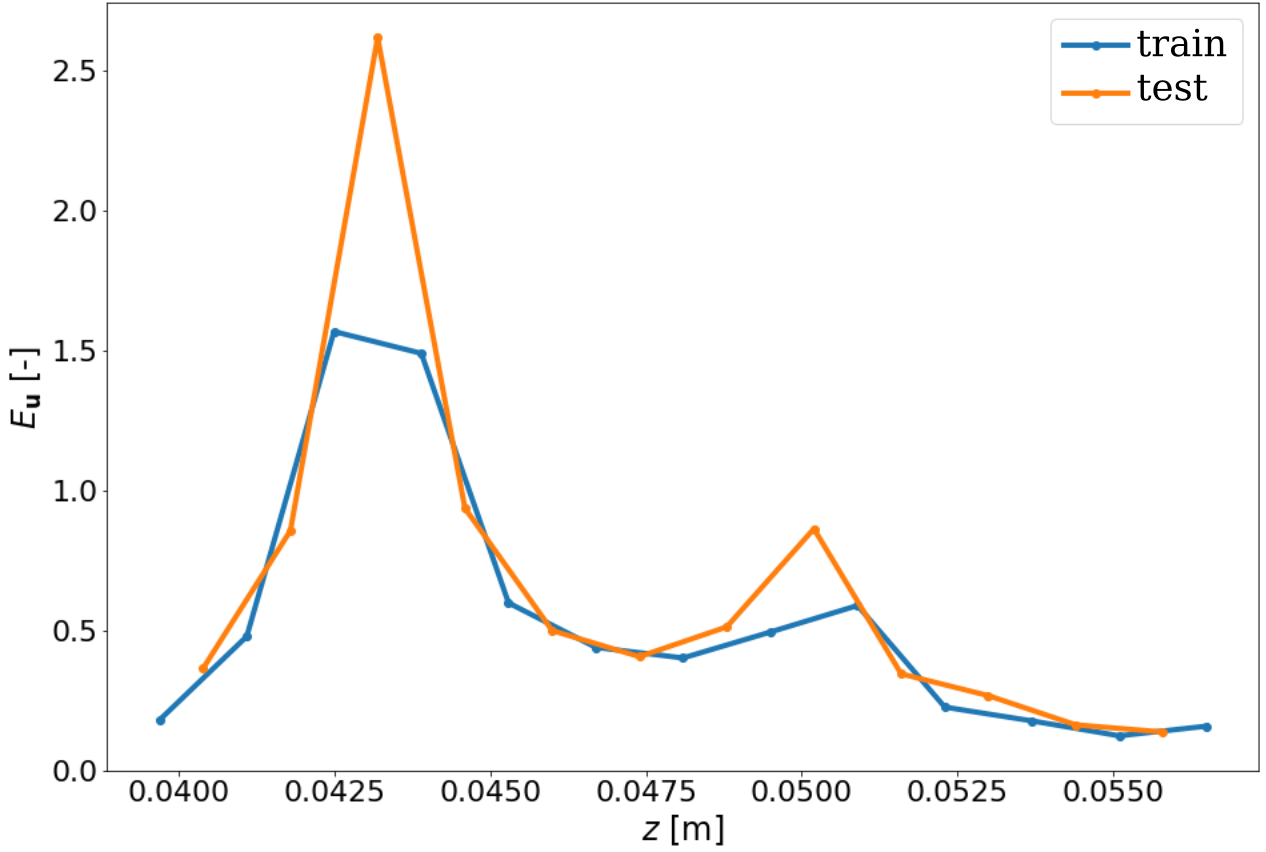}\label{fig:FEManeer3}
		}
		\caption[]{Normalized error $E_{\bm{u}}$ computed on different sections. \textbf{(a)} AN-CFD-1. \textbf{(b)} AN-CFD-2. \textbf{(c)} AN-CFD-3.}
		\label{fig:FEManeerrs}
	\end{figure}
	To avoid problems with the division by near-zero values, we compute the normalized error with respect to the average velocity magnitude over each train and test section using the following expression:
	\begin{equation}
		E_{\bm{u}} = \begin{dcases}
		    &\frac{|\bm{u}-\bm{u}_\mathrm{ref}|}{|\overline{\bm{u}}|},\quad\text{if} \ |\overline{\bm{u}}| > 10^{-3},\\
            & |\bm{u}-\bm{u}_\mathrm{ref}|, \quad \text{otherwise}.
		\end{dcases}
	\end{equation}
    We report the computed errors in Figure~\ref{fig:FEManeerrs}. The computed errors for AN-CFD-1 are substantially larger on the test data than on the training data, whereas the errors on the test data for AN-CFD-3 follow more closely the errors on the training data. Even so, their overall magnitude is significantly larger than for the other two cases. On the other hand, AN-CFD-2, provides a good compromise between the two, since the error on the test data follows the training data error trend and both errors are of a smaller magnitude than the ones for AN-CFD-3.
	\subsubsection{Validation against experimental data}
	\label{sec:aneVali}
	We validate the PINNs using 4D flow MRI data. First, in Section~\ref{sec:valiSet}, we provide a description of the 4D flow MRI data from~\cite{pravdivtseva2022influence} and the setup for the PINNs training. Next, in Section~\ref{sec:valiRes}, we examine the obtained results.
	\paragraph{Setup}
	\label{sec:valiSet}
	The 4D flow MRI data provided by~\cite{pravdivtseva2022influence} are available only in the aneurysm sac, as shown in Figure~\ref{fig:modane}. We report the fluid flow properties in Table~\ref{tab:fluid_flow_properties}, third column. The flow is pulsatile, \color{black}with the period corresponding to a single heartbeat of a duration of $0.816~\si{\second}$. The three velocity components are available in each point of the spatial domain in 24 time samples, with a time interval between each sample of $0.034~\si{\second}$.
    \par
	Following the error analysis in~\cite{pravdivtseva2022influence}, we use the the lower resolution $1.5~\si{\milli\meter\cubed}$ voxel dataset for training and the higher resolution $1.0~\si{\milli\meter\cubed}$ voxel dataset for comparison. We use the datasets with a CS factor of \num{2.5}. For the two datasets, \num{144} and \num{224} samples are given in the $x$ direction, \num{144} and \num{224} in the $y$ direction, and \num{40} and \num{60} samples in the $z$ direction. 
	\par
\begin{figure}[t]
		\centering
		\subfloat[][]{
			\includegraphics[width=0.3\textwidth]{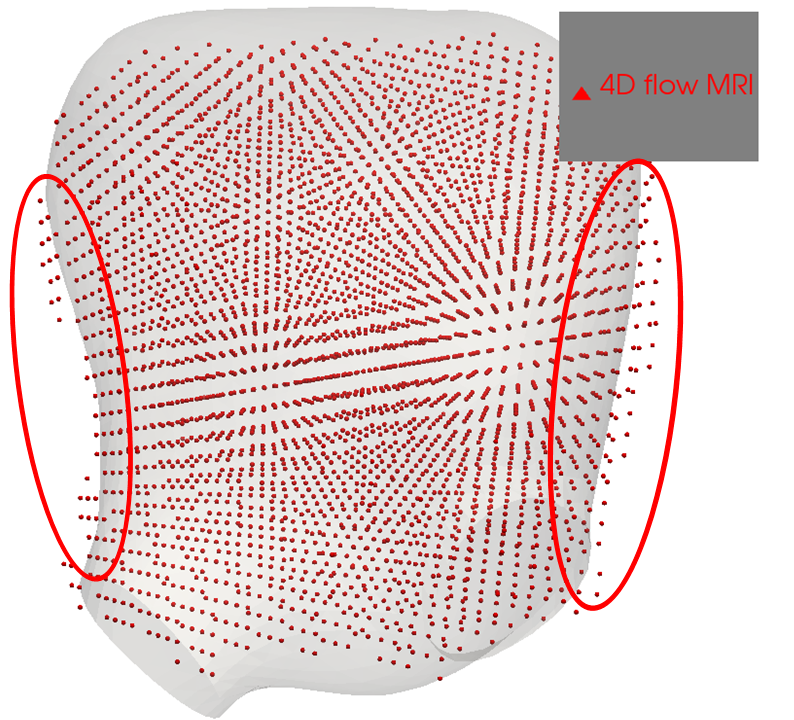}\label{fig:bcdisp}
		}
		\quad
		\subfloat[][]{
			\includegraphics[width=0.4\textwidth]{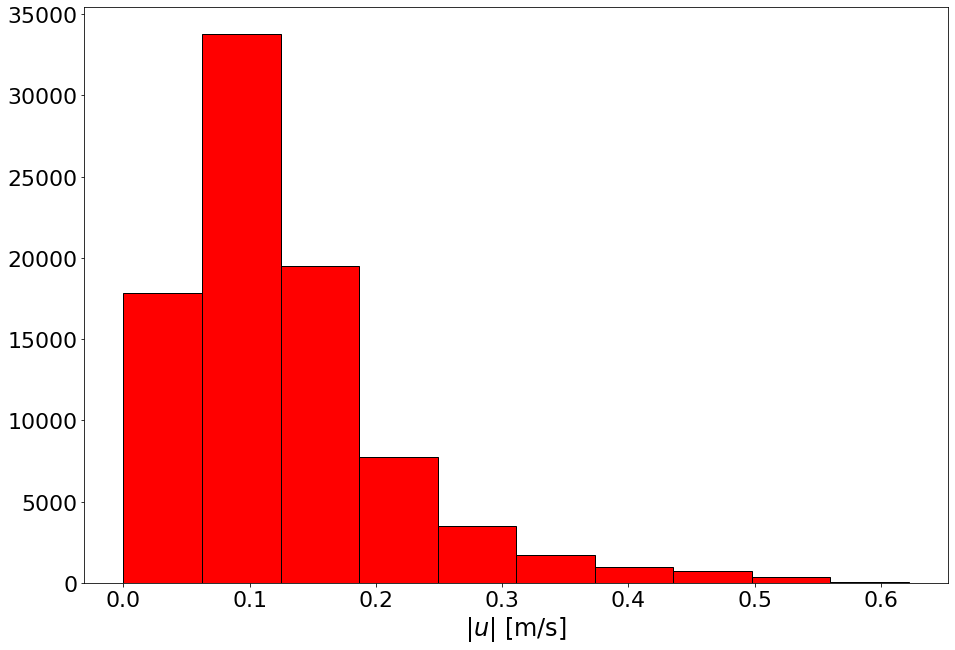}\label{fig:vdistr}
		}
		\caption[]{Collocation points and velocity distribution for the $1.5~\si{\milli\meter\cubed}$ voxel dataset. \textbf{(a)} Measurement points (red) and aneurysm domain (gray). Measurement points that cross the boundary are circled in red. \textbf{(b)} Recorded velocity magnitudes distribution.}
		\label{fig:filt4D}
	\end{figure}

    Due to the data having a different orientation from the model depicted in Figure~\ref{fig:modane}, we recenter the measurements with respect to the model. In Figure~\ref{fig:bcdisp} we depict the measurement points with respect to the aneurysm domain. The 4D flow MRI dataset is very noisy and dispersed at the boundary. The measurement points cross the lateral boundaries, while at the bottom of the aneurysm sac there is a lack of data, making the boundary very difficult to identify. We filter out the points outside of the domain and we report the resulting distribution of the velocity magnitudes over all time samples in Figure~\ref{fig:vdistr}. Most of the velocity magnitude samples have low values, close to $0.1~\si[per-mode=reciprocal]{\meter\per\second}$.
	\par
	\par
	We train the PINNs following what we did in Section~\ref{sec:veriSet}, using the loss function \eqref{eq:loss_ane}, where we neglect the inlet boundary loss term by setting $\lambda_{\mathrm{BC},1}=0$, since no information is available on its spatiotemporal distribution. We use ANNs made of three hidden layers, with batch normalization before each layer.
    \par
    We test four different PINN configurations at a single timestep corresponding to $t=0.34~\si{\second}$. Configurations AN-4DMRI-1, 2 and 3 correspond to decreasing weights on the data in the definition of the total loss function, whereas AN-4DMRI-4 disregards the regularization with the PDE residual $\bm{f}_\mathrm{PDE}$, and is used to assess the contribution of the wall no-slip boundary condition residual $\bm f_{\mathrm{BC},2}$. The parameters of the different PINN configurations are summarized in Table~\ref{tab:wane}.
	\begin{table}[t]
		\centering
			\begin{tabular}{ p{4em} p{6em} p{6em} p{6em} p{6em} }
				\hline
				& AN-4DMRI-1 & AN-4DMRI-2 & AN-4DMRI-3 & AN-4DMRI-4\\
				\hline
				$\lambda_{u_1}$& 170 & 17 & 1.7 & 1.7 \\
				$\lambda_{u_2}$& 170 & 17 & 1.7 & 1.7 \\
				$\lambda_{u_3}$& 170 & 17 & 1.7 & 1.7 \\
				$\lambda_{\hat{\bm{u}}}$& 17 & 1.7 & 0.17 & 0.17 \\
				$\lambda_{\mathrm{BC},1}$ & 0 & 0 & 0 & 0 \\
				$\lambda_{\mathrm{BC},2}$ & 21.6 & 21.6 & 21.6 & 36\\
			    $\lambda_{\mathrm{BC},3}$ & 0.5 & 0.5 & 0.5 & 0\\
				$\lambda_{\mathrm{PDE},1}$ & 5.76 & 5.76 & 5.76 & 0\\
				$\lambda_{\mathrm{PDE},2}$ & 1 & 1 & 1 & 0\\
				$\lambda_{\mathrm{PDE},3}$ & 1 & 1 & 1 & 0\\
				$\lambda_{\mathrm{PDE},4}$ & 1.56 & 1.56 & 1.56 & 0\\
				$\lambda_{\bar p}$ & 0.5 & 0.5 & 0.5 & 0 \\$N_{\mathrm{neuron}}$ & 32 & 32 & 32 & 16 \\
				\hline
			\end{tabular}
		\caption{Parameters of different loss functions used for PINNs training in Section~\ref{sec:valiSet}.}
		\label{tab:wane}
	\end{table}
    \par
\begin{figure}[t]
		\centering
		\subfloat[][]{
			\includegraphics[width=0.24\textwidth]{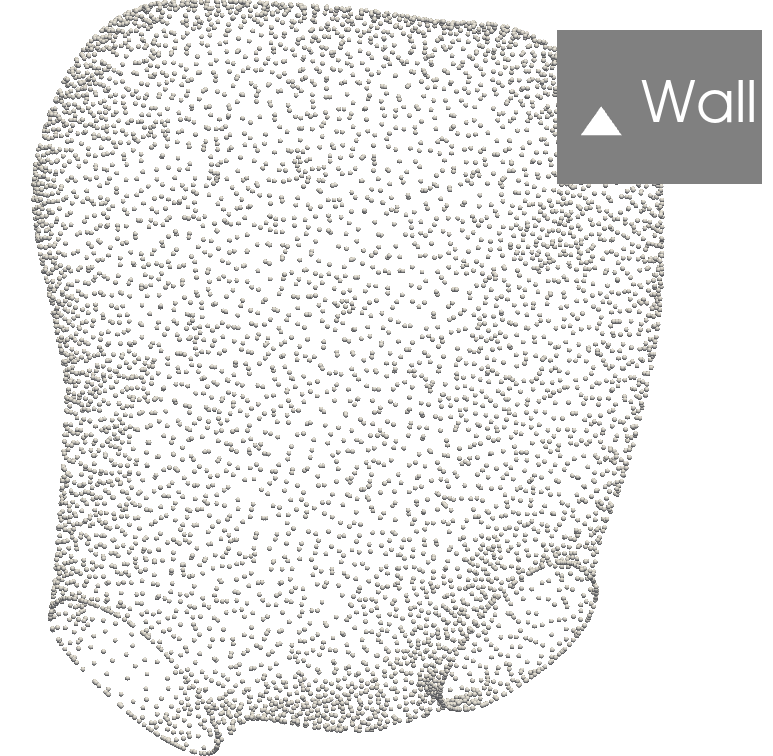}\label{fig:bcane}
		}
		\qquad
		\subfloat[][]{
			\includegraphics[width=0.24\textwidth]{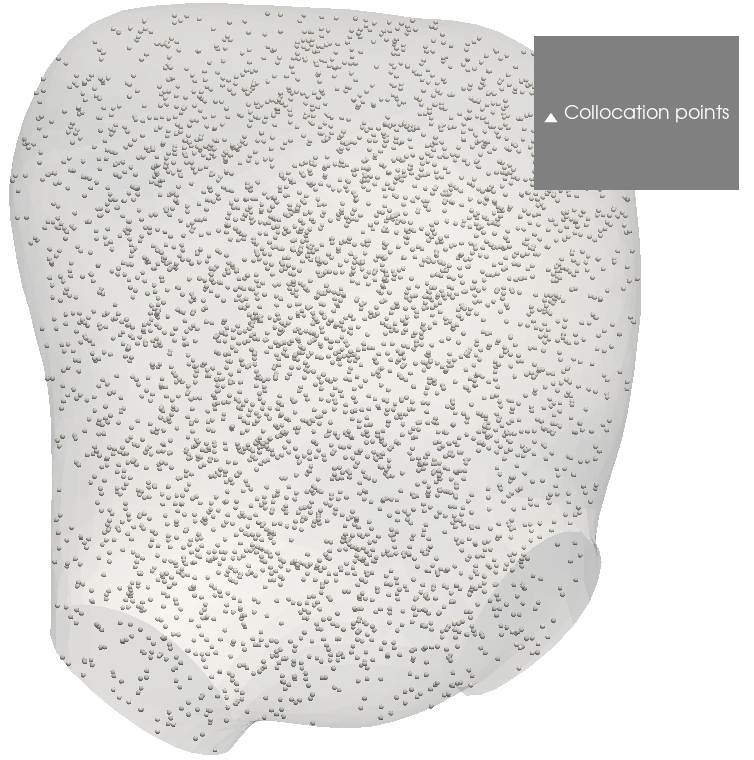}\label{fig:PDEane}
		}
		\caption[]{Collocation points. \textbf{(a)} Collocation points at wall. \textbf{(b)} PDE collocation points.}
		\label{fig:colane}
	\end{figure}

	For AN-4DMRI-1, 2 and 3, after filtering out the points outside of the domain, we use $N_\mathrm{train} = 1107$ points for training and $N_\mathrm{test} = 2249$ points for testing. We randomly seed $N_\mathrm{PDE} = 4420$ collocation points for the PDE loss term and $N_{\mathrm{BC},2} = 3198$ for collocation points at the wall. Wall boundary condition and PDE collocation points are depicted in Figures~\ref{fig:bcane} and~\ref{fig:PDEane}, respectively. For the fourth configuration, we use 80\% of the total dataset for training, and the remaining for testing purposes.
    
    \par
    For the training, we use $100$ epochs of the ADAM optimizer, followed by \num{14000} epochs of L-BFGS-B for the first three configurations and \num{7000} epochs for the fourth configuration.
    
    \afterpage{\clearpage}
    
\begin{figure}[t]
		\centering
		\subfloat[][]{
			\includegraphics[width=0.22\textwidth]{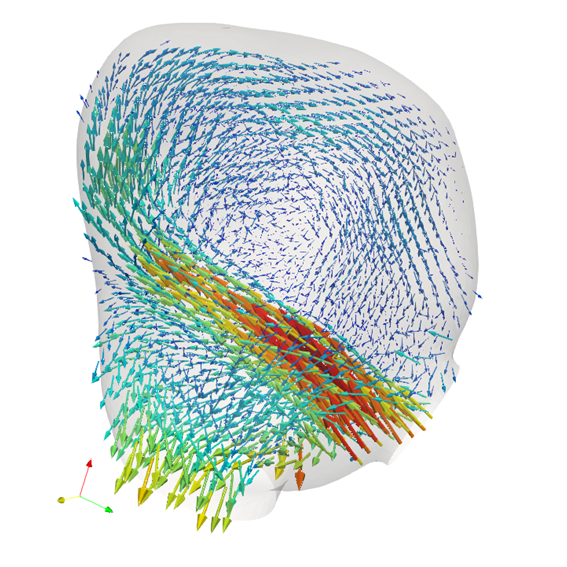}\label{fig:u04D}
		}
		\subfloat[][]{
			\includegraphics[width=0.22\textwidth]{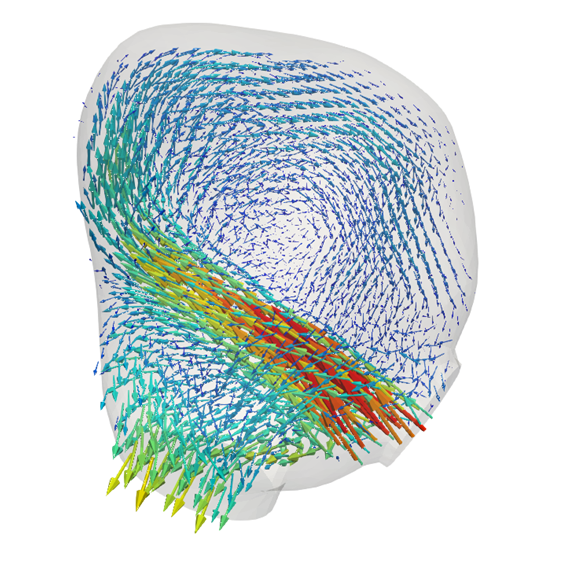}\label{fig:u14D}
		}
		\subfloat[][]{
			\includegraphics[width=0.22\textwidth]{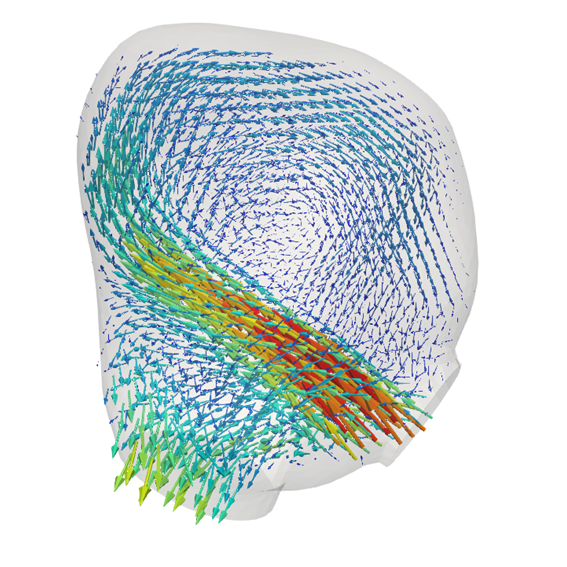}\label{fig:u24D}
		}
		\subfloat[][]{
			\includegraphics[width=0.22\textwidth]{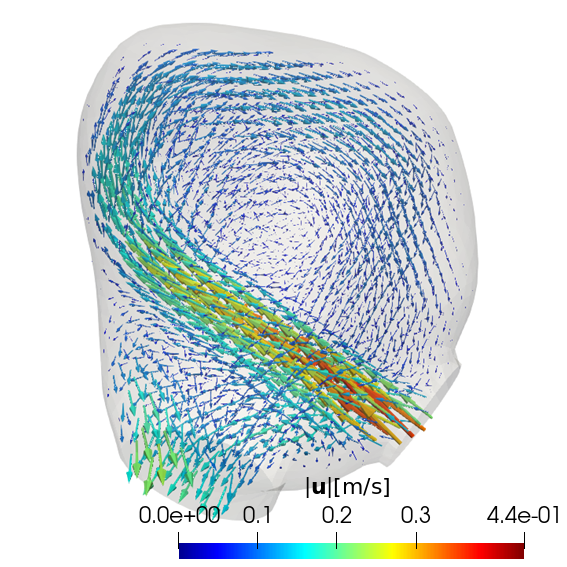}\label{fig:u34D}
		}
		\caption[]{Velocity field evaluated at the collocation points of the $1.5~\si{\milli\meter\cubed}$ voxel dataset. \textbf{(a)} 4D flow MRI $1.5~\si{\milli\meter\cubed}$ data. \textbf{(b)} AN-4DMRI-1. \textbf{(c)} AN-4DMRI-2. \textbf{(d)} AN-4DMRI-3.}
		\label{fig:velane}
	\end{figure}

\begin{figure}[t]
		\centering
		\subfloat[][]{
			\includegraphics[width=0.22\textwidth]{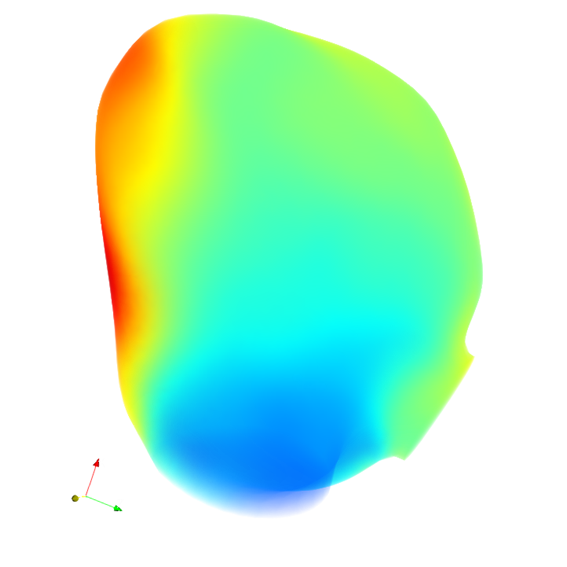}\label{fig:p14D}
		}
		\subfloat[][]{
			\includegraphics[width=0.22\textwidth]{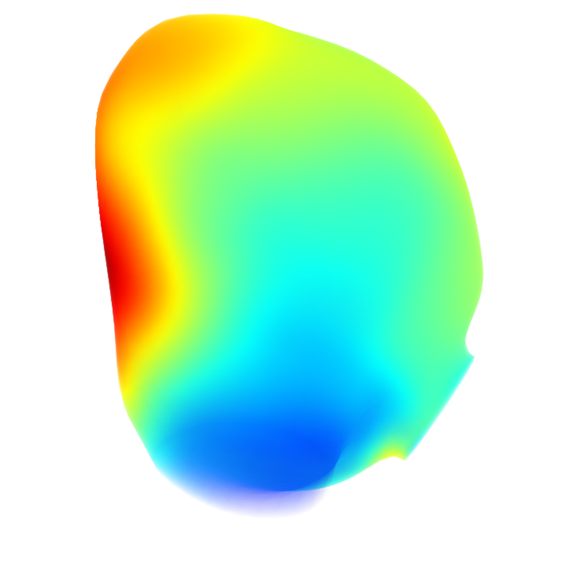}\label{fig:p24D}
		}
		\subfloat[][]{
			\includegraphics[width=0.22\textwidth]{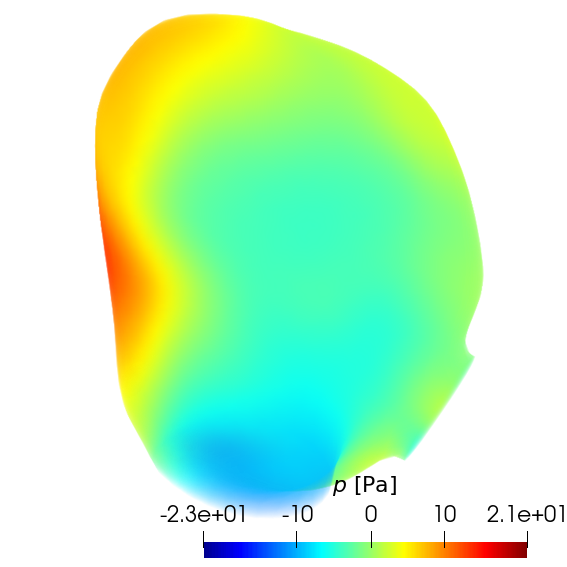}\label{fig:p34D}
		}
		\caption[]{Pressure field reconstruction. \textbf{(a)} AN-4DMRI-1. \textbf{(b)} AN-4DMRI-2. \textbf{(c)} AN-4DMRI-3.}
		\label{fig:press_ane}
	\end{figure}

	\paragraph{Results}
	\label{sec:valiRes}
	We report the reconstructed velocity fields from the PINNs along with the experimentally observed velocity field in Figure~\ref{fig:velane}. All PINNs retrieved the main velocity field features, however, by decreasing the the weights on the data, the peak velocity decreases and the result matches the boundary conditions better, particularly near the inlet, where the data are highly dispersed. We report the pressure reconstruction obtained from Eq.~\eqref{eq:mom} in Figure~\ref{fig:press_ane}.
    From Figure~\ref{fig:p24D} we see that all three pressure reconstructions have a similar pressure distribution. We remark that these represent non-invasive pressure measurements, as no pressure experimental data were introduced.
	\begin{figure}[t]
		\centering
		\subfloat[][]{
			\includegraphics[width=0.22\textwidth]{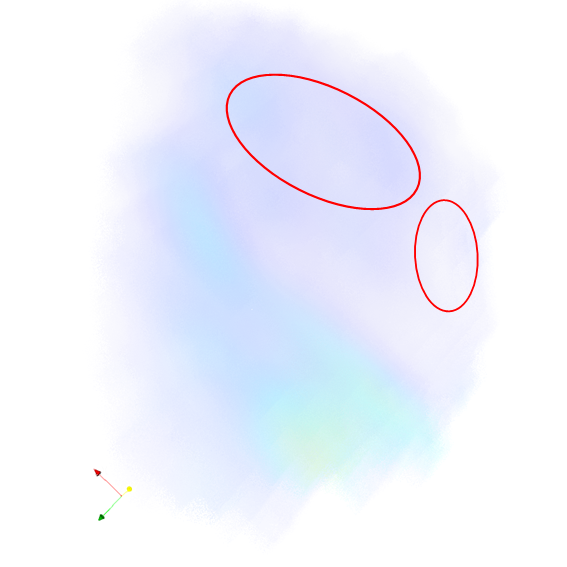}\label{fig:om154D}
		}
		\subfloat[][]{
			\includegraphics[width=0.22\textwidth]{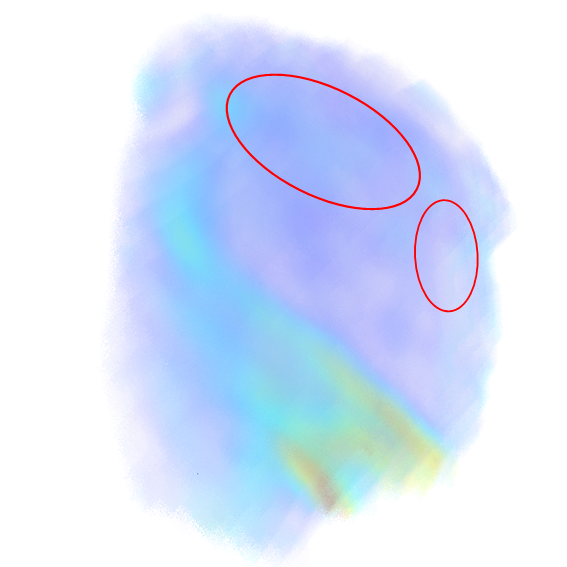}\label{fig:om104D}
		}
		\subfloat[][]{
			\includegraphics[width=0.22\textwidth]{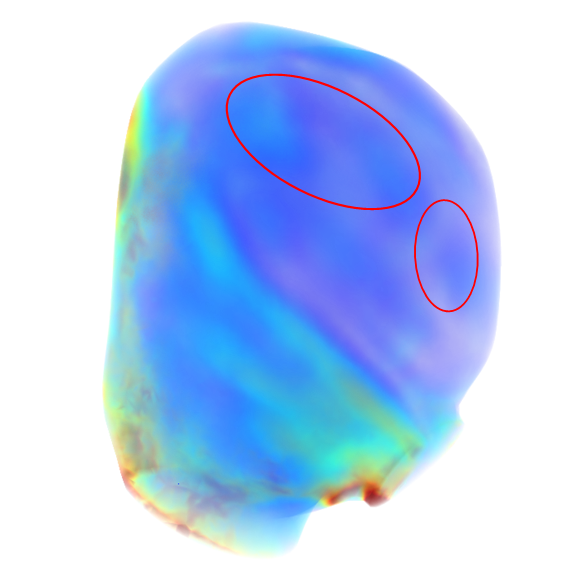}\label{fig:om14D}
		}\\
		\subfloat[][]{
			\includegraphics[width=0.22\textwidth]{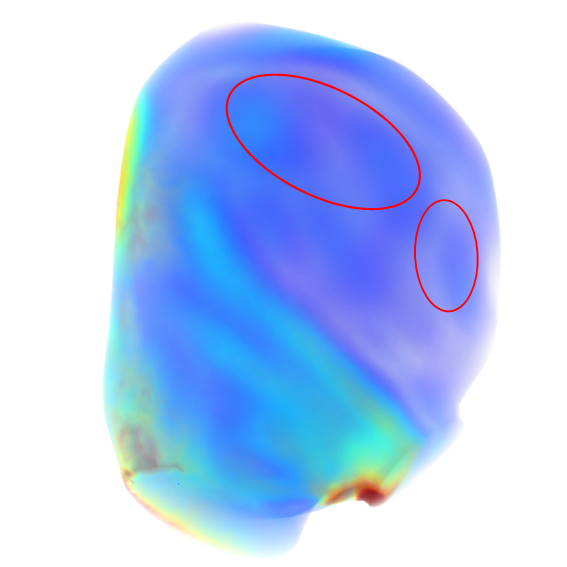}\label{fig:om24D}
		}
		\subfloat[][]{
			\includegraphics[width=0.22\textwidth]{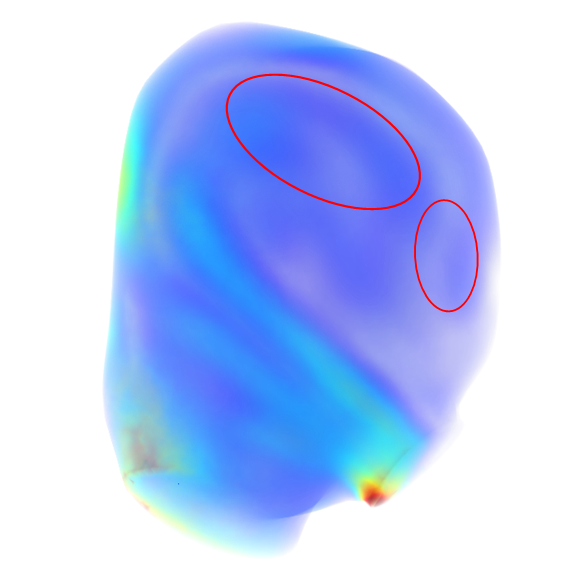}\label{fig:om34D}
		}
		\subfloat[][]{
			\includegraphics[width=0.22\textwidth]{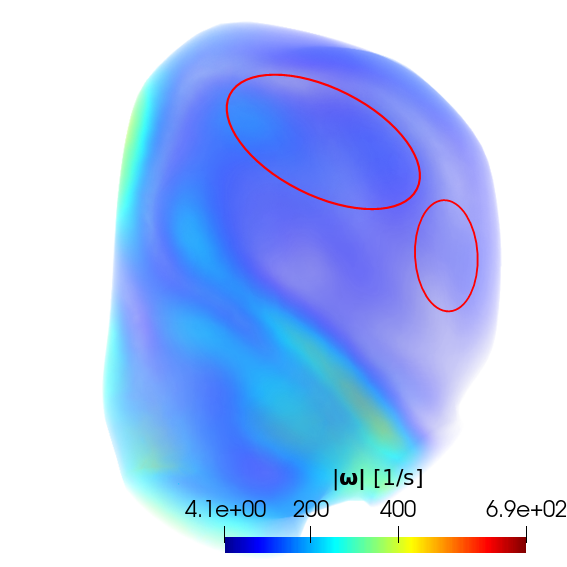}\label{fig:om44D}
		}
		\caption[]{Computed vorticity field for the 4D flow MRI data and for the PINN outputs. Recovered vorticity paths are circled in red. \textbf{(a)} 4D flow MRI $1.5~\si{\milli\meter\cubed}$ data. \textbf{(b)} 4D flow MRI $1.0~\si{\milli\meter\cubed}$ data. \textbf{(c)} AN-4DMRI-1. \textbf{(d)} AN-4DMRI-2. \textbf{(e)} AN-4DMRI-3. \textbf{(f)} AN-4DMRI-4.}
		\label{fig:vort_plot}
	\end{figure}
	Next, we compute the vorticity $\bm\omega=\nabla\cross\bm u$, and the wall shear stress $\bm{\mathrm{WSS}}$ for the reconstructed solutions. In Figure~\ref{fig:vort_plot}, we report the computed vorticities for the 4D flow MRI data and for the PINN outputs. AN-4DMRI-1, AN-4DMRI-2 and AN-4DMRI-3 all have a high value of the vorticity close to the inlet, the outlet and the left wall. Moreover, they all appear to have an overshoot of the vorticity near the inlet. As already discussed, the loss of accuracy at the inlet may be due to the fact that the Dirichlet inlet boundary condition is not available and the data at the inlet is scarce. We also see that using PINNs we are able to identify some vorticity paths that are not visible in the initial training data. Indeed, the ANN was trained using the $1.5\ \si{\milli\meter\cubed}$ voxel dataset, whereas the retrieved vorticity field is closer to the $1.0\ \si{\milli\meter\cubed}$ voxel data. This can be better appreciated for the AN-4DMRI-1 and AN-4DMRI-2 configurations.
    \par
\begin{figure}[t]
		\centering
		\subfloat[][]{
			\includegraphics[width = 0.4\textwidth]{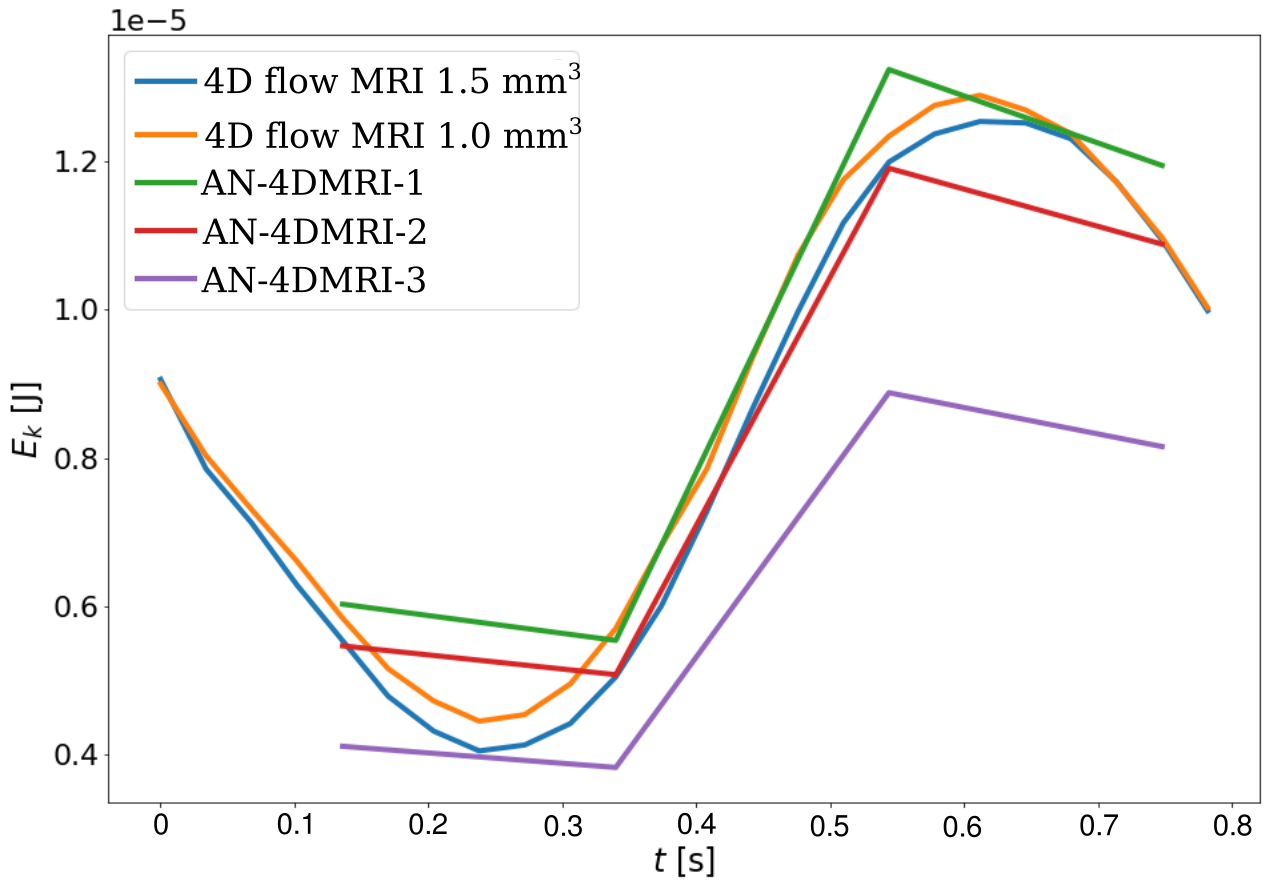}\label{fig:kinetic_energy}
		}
		\quad
		\subfloat[][]{
			\includegraphics[width = 0.4\textwidth]{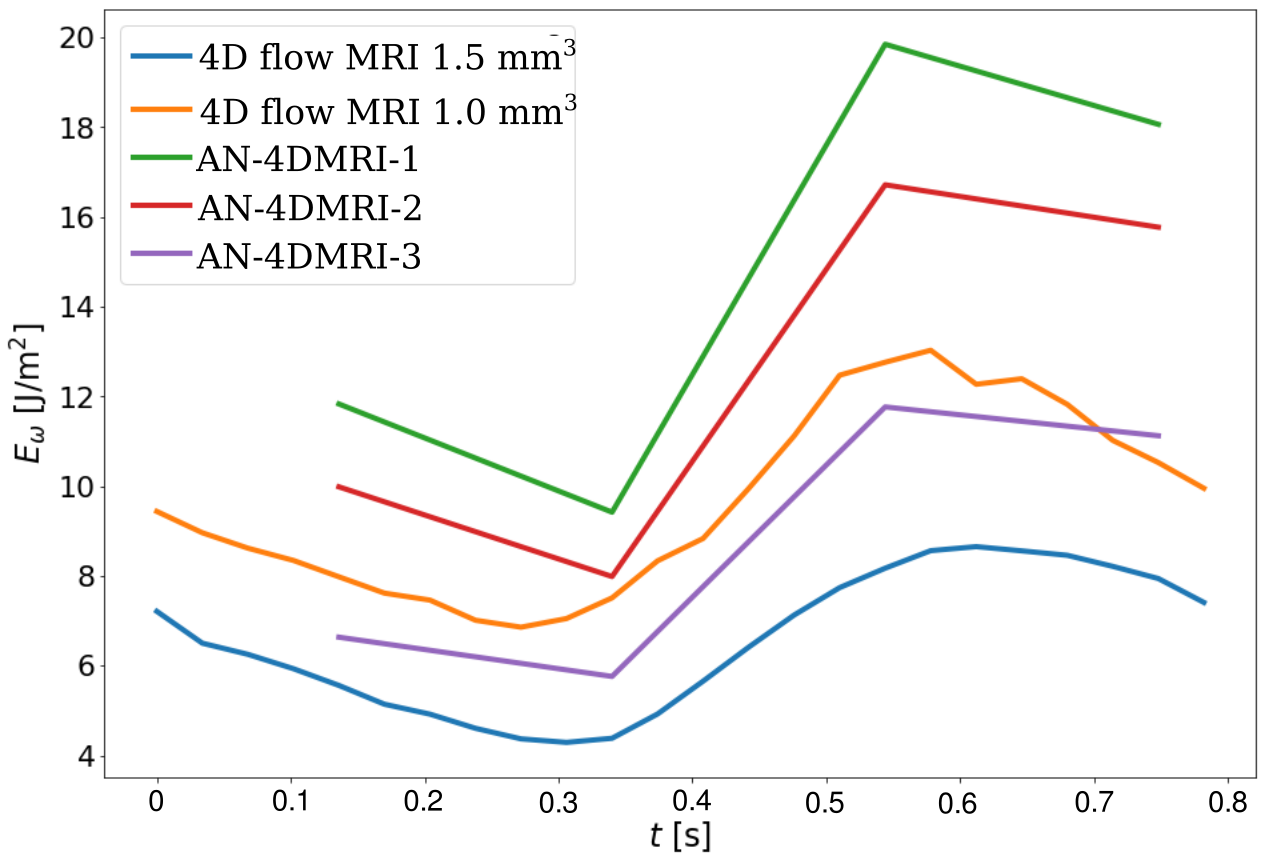}\label{fig:enstrophy}
		}
		\caption{Integral quantities during a single heartbeat. \textbf{(a)} Kinetic energy. \textbf{(b)} Enstrophy.}
		\label{fig:kinenstr}
	\end{figure}

	Our analysis can be extended to different timesteps. To quantify the variation of the flow properties during a single heartbeat, we introduce the kinetic energy $E_k$ and the enstrophy $E_{\bm\omega}$, both dimensional integral quantities computed as:
    \begin{equation*}
        E_k = \int_\Omega \frac{1}{2}\rho |\bm u|^2, \qquad E_{\bm\omega} = \int_\Omega \frac{1}{2}\rho |\bm\omega|^2,
    \end{equation*}
    in dimensional units. In Figure~\ref{fig:kinenstr}, we report the variation in time of the integral quantities associated to the flow. The flow field close to the inlet is disregarded during the integration for the enstrophy computation, due to the strong overshoots in the observed vorticity. The ANN models were trained at time instances of $t=0.136~\si{\second}$, $t=0.34~\si{\second}$, $t=0.544~\si{\second}$ and $t=0.748~\si{\second}$. In Figure~\ref{fig:kinetic_energy}, we show the variation of the kinetic energy during the course a single heartbeat. While AN-4DMRI-1 gives larger estimates of the kinetic energy, and AN-4DMRI-3 strongly underestimates it, AN-4DMRI-2 provides a good compromise between the two and overall a better estimate over the duration of the whole heartbeat. In Figure~\ref{fig:enstrophy} we report the enstrophy, which increases when going from the $1.5\ \si{\milli\meter\cubed}$ to the $1.0\ \si{\milli\meter\cubed}$ voxel dataset. On the other hand, the enstrophy computed from the PINN outputs decreases with decreasing weights on the data, i.e. from AN-4DMRI-1 to AN-4DMRI-3, in agreement with the distribution of the vorticity reported in Figure~\ref{fig:vort_plot}. We remark that the high dispersion of the enstrophy is in part due to errors in the computation of the vorticity, since it is performed with discrete derivatives.
    \par
\begin{figure}[t]
		\centering
		\subfloat[][]{
			\includegraphics[width=0.24\textwidth]{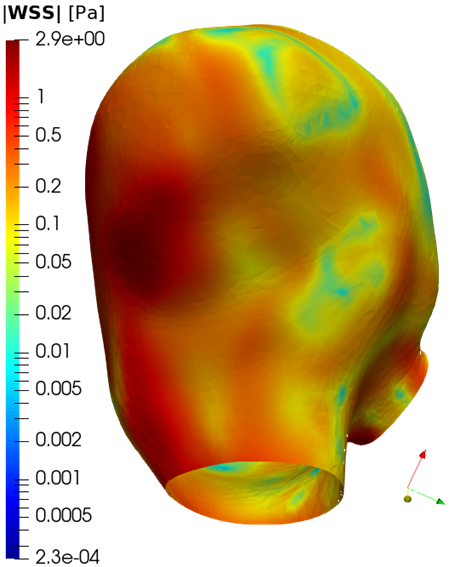}\label{fig:wss14D}
		}
		\subfloat[][]{
			\includegraphics[width=0.24\textwidth]{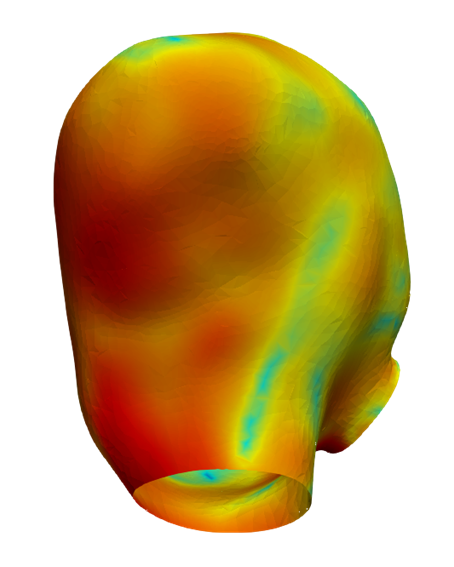}\label{fig:wss24D}
		}
		\subfloat[][]{
			\includegraphics[width=0.24\textwidth]{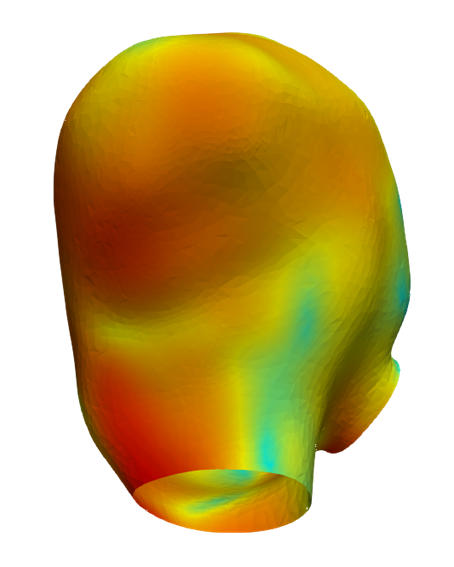}\label{fig:wss34D}
		}
		\subfloat[][]{
			\includegraphics[width=0.24\textwidth]{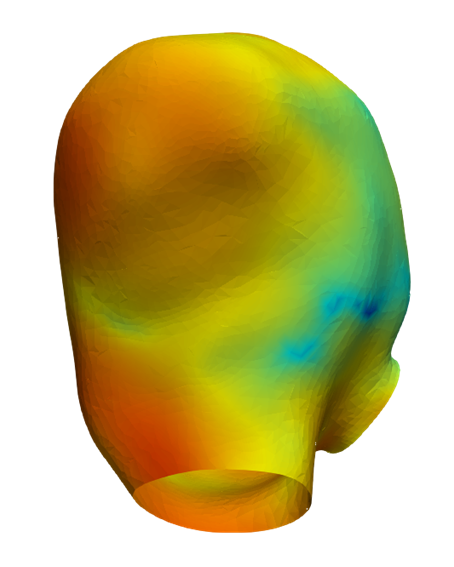}\label{fig:wss44D}
		}
		\caption[]{Computed wall shear stress magnitude in log scale. \textbf{(a)} AN-4DMRI-1. \textbf{(b)} AN-4DMRI-2. \textbf{(c)} AN-4DMRI-3. \textbf{(d)} AN-4DMRI-4.}
		\label{fig:WSS}
	\end{figure}

    We report in Figure~\ref{fig:WSS} the computed WSS. In all the networks, the highest WSS is concentrated in the same wall region. However, AN-4DMRI-4 strongly underestimates the wall shear stress, with an order of magnitude of difference. The wall shear stress prediction decreases from AN-4DMRI-1 to AN-4DMRI-3. Moreover, an overshoot of the wall shear stress appears near the inlet. As discussed for the vorticity, this may be due to the scarcity of data. The results of the wall shear stress are poor and difficult to assess, since so much dispersion arises in our study.
    \par
	Lastly, we compute the errors as the square root of the training and test loss function for the velocity at $t=0.34~\si{\second}$. We compute an average absolute error for the $1.5~\si{\milli\meter\cubed}$ voxel dataset. The results are summarized in Table~\ref{tab:errdataA}.
	
\begin{table}[t]
		\centering
		\begin{tabular}{ p{5em} p{6em} p{6em} p{6em} p{6em} }
        \hline
        & AN-4DMRI-1 & AN-4DMRI-2 & AN-4DMRI-3 & AN-4DMRI-4\\
			\hline
			 $\sqrt{\mathcal{L}_\mathrm{train}^{\bm{u}}}$ & $2.10\cdot10^{-2}$ & $4.20\cdot10^{-2}$ & $7.35\cdot10^{-2}$ & $6.76\cdot10^{-2}$\\
			$\sqrt{\mathcal{L}_\mathrm{test}^{\bm{u}}}$ & $4.21\cdot10^{-2}$ & $4.83\cdot10^{-2}$ & $6.40\cdot10^{-2}$ & $5.39\cdot10^{-2}$ \\
			$e_{\bm{u}}~[\si[per-mode=reciprocal]{\meter\per\second}]$ & $8.07\cdot10^{-3}$ & $1.11\cdot10^{-2}$ & $1.63\cdot10^{-2}$ & $1.46\cdot10^{-2}$ \\
			\hline	
		\end{tabular}
		\caption{Computed errors for the four PINN configurations in Section~\ref{sec:aneVali}. 
        }
		\label{tab:errdataA}
	\end{table}

\section{Discussion}
\label{sec:discussion}
This study demonstrates the potential of PINNs to integrate sparse, noisy, or under-resolved experimental data with physical modeling, enabling more accurate and spatially resolved reconstruction of hemodynamic quantities. Our focus was on enhancing 4D flow MRI data, which, while non-invasive, suffers from limited spatial and temporal resolution, especially near vessel walls and in regions of complex flow. By embedding the incompressible Navier–Stokes equations into the PINNs, we leveraged both data and physical laws to improve the reconstruction the velocity field.
\par
The velocity, pressure, and WSS reconstruction have shown strong agreements both with FEM simulations and with experimental data for the FDA nozzle benchmark. Initially, training the ANN using only data without physics-based regularization yielded poor results, particularly when data were limited or spatially under-resolved. This limitation was addressed by incorporating a the Navier-Stokes residuals into the loss function, which acted as a regularizer and compensated for the sparse nature of the data. The physics-based regularization not only improved the accuracy of the reconstructed fields but also enabled the estimation of quantities not directly measured, such as the pressure and WSS.
\par
One of the main challenges in PINN training was the a priori selection of optimal hyperparameters for the numerical optimization. This process can be automated using software libraries for hyperparameter tuning, such as \texttt{scikit-learn}’s \texttt{GridSearchCV}~\cite{JMLR:v12:pedregosa11a}. Another significant challenge was the minimization of loss terms spanning multiple orders of magnitude. We mitigated this by introducing a logarithmic transformation of the residuals in the loss function, which improved the network’s ability to capture both high- and low-magnitude flow features.
\par
The placement and density of collocation points also played a crucial role. Increasing the number of collocation points in regions of high flow gradients, such as the convergent section of the nozzle, enhanced the training process and the accuracy of the results. Furthermore, incorporating velocity unit vectors and increasing the number of neurons in the ANN improved the reconstruction without inducing overfitting. When trained with PIV measurements, PINNs provided results that more closely matched experimental observations compared to classical CFD simulations, highlighting the robustness of the physics-informed approach.
\par
For the aneurysm case, PINNs demonstrated the ability to compute key hemodynamic quantities of interest, such as velocity, pressure and WSS, from 4D flow MRI data. The integration of physical modeling was particularly valuable in regions where data were sparse or noisy, such as near vessel walls. The reconstructed velocity fields were more accurate and the vorticity distributions better matched high-resolution 4D flow MRI data compared to pure data-driven approaches. While the estimation of WSS remains challenging due to its high sensitivity to boundary conditions and data quality, the physics-based regularization was essential for obtaining non-invasive pressure information and improving the overall fidelity of the hemodynamic indicators.
\subsection{Limitations}
\label{sec:limitations}
While the proposed procedure can be easily extended to other applications, the present analysis is limited to two specific test cases: the FDA nozzle benchmark and an in vitro aneurysm model.
\par
Regarding the FDA nozzle benchmark, only the laminar flow regime was considered. However, transitional and turbulent flow regimes, which are common in physiological scenarios, present more complex fluid flow structures and were not explored in this work. Future studies should investigate the PINNs performance under these conditions.
\par
The present analysis also relies on a quasi-steady flow assumption to reduce the number of parameters in the PINNs. Although the aneurysm case study incorporates pulsatile flow, the framework’s ability to accurately capture rapid, unsteady flow features, such as those occurring during systolic peaks, has not been thoroughly validated. The quasi-steady assumption may not hold in highly dynamic scenarios, potentially limiting the accuracy of the results.
\par
The aneurysm case study relied on in vitro 4D flow MRI data which are inherently less noisy and more controlled than in vivo acquisitions. In vivo 4D flow MRI data are often affected by significant noise, motion artifacts and measurement errors. The robustness of the proposed approach in handling such noisy, real-world data remains to be thoroughly tested. Additionally, the lack of a ground truth for the aneurysm model limited the ability to validate the reconstructed quantities, particularly the wall shear stress (WSS), which exhibited high variability and dispersion in the results.
\par
The accuracy of the PINN also depends on the correct specification of the boundary conditions. In the aneurysm case study, the lack of precise inlet velocity profiles or outlet pressure conditions introduced uncertainties, particularly near the boundaries. This limitation is exacerbated in patient-specific geometries, where boundary conditions are often approximated or unknown.
\par
The performance of the PINN framework is highly sensitive to the selection of hyperparameters, including the weights assigned to different terms in the loss function and the number and placement of collocation points for the PDE regularization. Balancing data-driven and physics-based loss terms is non-trivial and often requires extensive trial and error. Suboptimal weights can lead to poor convergence, nonphysical solutions, or overfitting, as observed in some of the test cases. Due to the manual optimization of the PINN hyperparameters, there is no guarantee of their optimality. Further works could explore the integration of automated optimization tools into the current framework and strategies for optimal collocation point placement, particularly in complex geometries where flow features are not known a priori.
\section{Conclusion}
\label{sec:conclusions}
In this work, we developed and tested a PINN-based method to enhance measurement-based data, with a primary emphasis on improving data from 4D flow MRI. Our findings highlight the value of incorporating physical knowledge with sparse or noisy measurements, enabling the reconstruction of high-fidelity hemodynamic fields and estimation of quantities not directly measured, such as the pressure and the WSS. 
\par
For the FDA nozzle benchmark, we successfully combined experimental observations with the Navier-Stokes equations, demonstrating that physics-based regularization improves the reconstruction of velocity and pressure fields from localized samples. Training with PIV measurements confirmed that PINNs can achieve good agreement with experimental data, outperforming classical CFD approaches.
\par
In the aneurysm case, PINNs improved the reconstruction of velocity fields and vorticity distributions, and enabled non-invasive pressure estimation from 4D flow MRI data. While the WSS estimation remains challenging, the integration of physical modeling was crucial for obtaining meaningful results in data-sparse regions in the study of the aneurysm.
\par
Advancements in machine learning techniques, such as adaptive activation functions~\cite{jagtap2020adaptive}, or enhanced 4D flow MRI preprocessing, as for instance using a singular value decomposition~\cite{henry19928}, could further enhance accuracy. Overall, we have shown PINNs to be a valuable tool for enhancing data -- particularly for emerging technologies such as 4D flow MRI -- by effectively combining sparse measurements with physical modeling.

\section*{Acknowledgements}
LD and IR acknowledge the support by the FAIR (Future Artificial Intelligence Research) project, funded by the NextGenerationEU program within the PNRR-PE-AI scheme (M4C2, investment 1.3, line on Artificial
Intelligence), Italy. LD, SP, and IR acknowledge their membership to INdAM GNCS -- Gruppo Nazionale per il Calcolo Scientifico (National Group for Scientific Computing), Italy. The research of LD, SP, and IR is part of the activities of ``Dipartimento di Eccellenza 2023–2027", MUR, Italy, Dipartimento di Matematica, Politecnico di Milano. 
\afterpage{\clearpage}
 
	\appendix
	\section{Appendix A}
	
	\subsection{Training sections coordinates}
	We report the $z$ coordinates of the training cross-sections for the FDA nozzle benchmark training by in silico data in Section~\ref{sec:veriFDA} and tests in Sections~\ref{sub:test1}, ~\ref{sub:test2}, ~\ref{sub:test3}, ~\ref{sub:test4}, and~\ref{sub:test5}.
	\begin{table}[h]
		\centering
		\begin{tabular}{c c}
			\hline
			\multirow{2}*{Figure~\ref{fig:trainPINNs} $z$ coords [m]} & -0.179, -0.145, -0.108, -0.07, -0.057,\\
			& -0.039, -0.016, 0.058, 0.099, 0.141\\
			\hline
		\end{tabular}
		\caption{$z$ coordinates for the training sections from FEM simulation.}
		\label{tab:trainz}
	\end{table}

    \subsection{FEM simulation parameters}
    We report the FEM simulation parameters from Sections~\ref{sec:veriFDA} and~\ref{sec:aneVeri}.
    \begin{table}[H]
		\centering 
		\begin{tabular}{p{8em} p{6em} p{6em} p{6em} p{6em}}
			\hline
			& FDA-coarse & FDA-fine & ANE-coarse & ANE-fine \\
			\hline
			Time step~$[\si{\second}]$ & $10^{-3}$ & $10^{-3}$ & $10^{-3}$ & $10^{-3}$  \\
			FE spaces & $\mathbb{Q}_1-\mathbb{Q}_1$ & $\mathbb{Q}_1-\mathbb{Q}_1$ & $\mathbb{P}_1-\mathbb{P}_1$  & $\mathbb{P}_1-\mathbb{P}_1$ \\
			Number of cells & \num{7325} & \num{36128} & \num{133166} & \num{1242082}\\
			Velocity DOFs & \num{24696} & \num{116442} & \num{78621} &  \num{632682} \\
			Pressure DOFs & \num{8232} & \num{38814} & \num{26207} & \num{210894} \\
			Total DOF & \num{32928} & \num{155256} & \num{104828} & \num{843576}\\
			$h_{\mathrm{min}}~[\si{\meter}]$ & $2.31\cdot10^{-3}$ & $1.23\cdot10^{-3}$ & $3.87\cdot10^{-4}$ & $1.88\cdot10^{-4}$ \\
			$h_{\mathrm{max}}~[\si{\meter}]$ & $8.38\cdot10^{-3}$ & $4.10\cdot10^{-3}$ & $1.14\cdot10^{-3}$ & $6.32\cdot10^{-4}$ \\
			$h_\mathrm{mean}~[\si{\meter}]$ & $6.20\cdot10^{-3}$ & $3.10\cdot10^{-3}$ & $7.93\cdot10^{-4}$ & 
			$3.76\cdot10^{-4}$\\
			\hline
		\end{tabular}
		\\[10pt]
		\caption{Simulation parameters for the FDA nozzle benchmark problem for the coarse and fine hexahedral meshes.}
		\label{table:FEM}
	\end{table}

	
	\bibliographystyle{elsarticle-num-names} 
	\bibliography{bibliography.bib}
    \biboptions{numbers,sort&compress}
	
	
	
	
	
\end{document}